\documentclass[11pt,a4paper]{article}

\usepackage{iftex}
\ifPDFTeX
  \usepackage[utf8]{inputenc}
  \usepackage[T1]{fontenc}
  \usepackage{lmodern}
\fi

\usepackage[a4paper,hmargin=3.5cm,vmargin=3cm]{geometry}
\usepackage{amsmath,amssymb,amsthm,mathtools}
\usepackage{graphicx}
\usepackage{booktabs}
\usepackage{array}
\usepackage[table]{xcolor}
\usepackage{enumitem}
\usepackage[font=small,labelfont=bf]{caption}
\usepackage{microtype}
\usepackage[colorlinks=true,linkcolor=blue!55!black,citecolor=blue!55!black,urlcolor=blue!55!black,pdftitle={Conjugation Diameters and Covering Numbers within the General Linear Family over Arbitrary Fields. I. Special Linear Diameter and Determinant Fibres},pdfauthor={Juliet Cooke}]{hyperref}
\usepackage[capitalise,noabbrev]{cleveref}

\definecolor{tabband}{HTML}{F6F3ED}
\definecolor{tabfocal}{HTML}{FAF3E3}
\definecolor{tabmuted}{HTML}{595145}

\DeclareMathOperator{\im}{im}

\theoremstyle{plain}
\newtheorem{theorem}{Theorem}[section]
\newtheorem{proposition}[theorem]{Proposition}
\newtheorem{lemma}[theorem]{Lemma}
\newtheorem{corollary}[theorem]{Corollary}

\theoremstyle{definition}
\newtheorem{definition}[theorem]{Definition}
\newtheorem{example}[theorem]{Example}

\theoremstyle{remark}
\newtheorem*{remark}{Remark}

\newenvironment{introresult}[2]{%
  \par\medskip\noindent\textbf{\Cref{#1}\ (#2).}\ \itshape
}{\par\medskip}

\title{\textbf{Conjugation Diameters and Covering Numbers within the General Linear Family over Arbitrary Fields}\\[0.4em]
\large I. Special Linear Diameter and Determinant Fibres}
\author{Juliet Cooke\thanks{Email: \href{mailto:juliet@julietcooke.net}{\texttt{juliet@julietcooke.net}}\hspace{0.6em}ORCID: \href{https://orcid.org/0000-0003-3131-3862}{\texttt{0000-0003-3131-3862}}}}
\date{November 2026}

\begin{document}
\maketitle

\begin{abstract}
For every field $F$ and every $n \ge  3$, we determine the conjugation diameter
$\Delta(\mathrm{SL}_n(F))$ and the abelianised covering number
$\operatorname{cn}_\mathrm{ab}(\mathrm{GL}_n(F))$, proving that they coincide and equal
\[
  \Delta(\mathrm{SL}_n(F)) = \operatorname{cn}_\mathrm{ab}(\mathrm{GL}_n(F)) = n + [\mu_n(F) \neq  \{1\}],
\]
where $\mu_n(F)$ is the group of $n$th roots of unity in $F$ and
$[P] \in \{0,1\}$ is the Iverson bracket. The same value gives the covering number
$\operatorname{cn}_\mathrm{ng}(\mathrm{SL}_n(F))$ for single noncentral $\mathrm{SL}_n(F)$-conjugacy
classes. In every case, transvections attain the worst length: the dimension
sets the scale $n$, while the field enters solely through nonidentity central
scalars in $\mathrm{SL}_n(F)$. Extending existence results of Ellers and Lausch, we
also establish a complete criterion for expressing a matrix in $\mathrm{GL}_n(F)$ as
a product of conjugates of prescribed rank-one transformations at independent
scalar origins, which supplies the common lower bounds.

\end{abstract}

\setcounter{tocdepth}{2}
\tableofcontents
\clearpage

\section{Introduction}
\label{sec-introduction}

This paper studies word length when assembling an arbitrary invertible linear
transformation from matrices of prescribed types, such as shears, reflections,
or conjugates of a fixed matrix and its inverse. Given a finite set $S$ of
\emph{seeds} in a group $G$, allow as letters every conjugate in $G$ of a seed or
its inverse. In a matrix group, conjugation is change of basis: conjugating
individual factors independently lets us place each transformation in whatever
basis we like while keeping its similarity type fixed.

The qualitative question of whether such words eventually reach every element
of $G$, namely \emph{normal generation}, has a well-known answer in the general
linear family. For $n\ge 3$ and any field $F$, a finite set $S$ normally generates
$\mathrm{SL}_n(F)$ if and only if $S$ contains a noncentral matrix
(\Cref{prop-sl-normal-generation}), and it normally generates $\mathrm{GL}_n(F)$ if
and only if $S$ contains a noncentral matrix and the determinants of its
members generate $F^\times$ (\Cref{prop-gl-normal-generation}). \emph{How many factors}
it takes to reach an arbitrary target is a separate question.

We measure a target by its \emph{word length} with respect to $S$, the fewest
letters needed to express it (zero for the identity), counting only the
factors and not the basis changes. The \emph{diameter} of $G$ for $S$ is the
largest word length of any element of $G$, and the \emph{conjugation diameter}
$\Delta(G)$ is the supremum of these diameters over all finite normally
generating sets $S$ (set to $\infty$ if none exists;
\Cref{sec-diameter-definitions}). A bound $\Delta(G)\le k$ is uniform across all
admissible seed sets: it owes nothing to choosing convenient generators.
The invariant $\Delta$ was introduced by K\k{e}dra, Libman and Martin
\cite{KedraLibmanMartin2023} (building on Lawther--Liebeck \cite{LawtherLiebeck1998} and
Burago--Ivanov--Polterovich \cite{BuragoIvanovPolterovich2008}) and named by Libman
and Tarry \cite{LibmanTarry2020}.

Studying $\mathrm{GL}_n(F)$ and $\mathrm{SL}_n(F)$ over an \emph{arbitrary} field $F$ puts
this word-length question in its basic linear-algebraic setting: invertibility
imposes only what makes composition a group, with no preserved metric or
bilinear form, while requiring all factors and conjugating matrices to have
entries in $F$ separates what is forced by dimension from what depends on the
arithmetic of the field.

\subsection{Special linear diameter}
\label{sec-intro-sl}
Consider first $\mathrm{SL}_n(F)$ ($n\ge 3$), where every noncentral seed normally
generates the group and all determinants are fixed at $1$. Our first main
theorem gives the exact conjugation diameter over every field. Write
$\mu_n(F)=\{\zeta \in F^\times : \zeta^n=1\}$ for the $n$th roots of unity in $F$,
and $[P] \in \{0,1\}$ for the Iverson bracket, equal to $1$ when $P$ holds and
$0$ otherwise.

\begin{introresult}{thm-sl-diameter}{Special linear diameter}
Let $F$ be a field and $n\ge 3$. Then
  \[
  \Delta(\mathrm{SL}_n(F)) = n + [\mu_n(F) \neq  \{1\}].
\]
  The supremum is attained by taking a single transvection as seed.
\end{introresult}

Thus the diameter is $n$ if $\mu_n(F)=\{1\}$ and $n+1$ otherwise: about one
factor per dimension. In $\mathrm{SL}_3(\mathbb{R})$, for example, for any matrix $g \neq  I$,
every matrix in $\mathrm{SL}_3(\mathbb{R})$ is a product of at most three conjugates of
$g$ or $g^{-1}$ by real matrices of determinant one.

Why does the dimension $n$ set this scale? To bound word length from below, we
attach to each matrix a conjugation invariant: its \emph{deviation rank}
from a scalar, which measures the number of independent directions in which it
differs from a scalar transformation (\Cref{sec-deviation-geometry}). Because
conjugation preserves this rank, each factor in a product contributes at most
that many independent directions of deviation. A transvection fixes a
hyperplane pointwise and moves vectors along only one line, whereas over every
field $\mathrm{SL}_n(F)$ contains matrices with no nonzero fixed vector; displacing
all $n$ directions therefore requires at least $n$ transvections.

Another choice of seed does not worsen this bound: a transvection, having
minimal deviation rank, already attains the diameter, and no noncentral seed is
slower (though seeds moving more directions can reach all targets in fewer
factors). The proof of the upper bound proceeds in two stages. A matrix $A$ is
\emph{cyclic} (admits a cyclic vector) precisely when $A - c I$ has rank at least
$n-1$ for every scalar $c$ in an algebraic closure of $F$; distance from
cyclicity is therefore again measured by deviation rank, taken over an
algebraic closure and counted down from $n-1$. In the first stage, independent
conjugation allows successive factors to act on new directions: a rank-one seed
adds one direction at each step, whereas any other seed has minimum deviation
rank at least two over an algebraic closure and reaches a cyclic matrix in
roughly half as many factors. In the second stage, two suitably placed cyclic
matrices multiply to produce every nonscalar matrix of the required
determinant. For $|F|\ge 4$ we use Lev's nonscalar coverage theorem
\cite[Theorem 2]{Lev1996}, and for $\mathbb{F}_3$ with $n\neq 4$ we use Kn\"uppel and Nielsen's
extended covering theorem \cite{KnuppelNielsen2006}.

Over $\mathbb{F}_2$ (all $n\ge 3$) and in $\mathrm{SL}_4(\mathbb{F}_3)$, rank-one seeds are
handled directly by \Cref{thm-prescribed-rank-one}, while for every other seed the
same two stages still work with one extra factor in the second stage. Over any
field, two cyclic classes multiply to cover all nonscalar matrices provided one
has a simple eigenvalue in $F$, by B\"unger and Nielsen's cyclic decomposition
theorem \cite[Theorem 1.2]{BungerNielsen1999}. Larger fields supply that eigenvalue
at no extra cost; over $\mathbb{F}_2$ and $\mathbb{F}_3$ we spend one additional factor
to arrange it, so nonscalar coverage takes $2j+1$ factors once $j$ factors
reach a cyclic class. For most seeds this stays within $n$ factors, exceeding
$n$ or leaving no room for the scalars only for seeds of minimum rank two over
$\mathbb{F}_2$ and in dimension four. We handle these residual cases in
\Cref{sec-small-field-proofs}:

\begin{itemize}
  \item \textbf{Even $n\ge 6$ over $\mathbb{F}_2$:} we replace each consecutive pair of
transvections in an $n$-transvection word by two conjugates of the rank-two
seed.
  \item \textbf{Odd $n\ge 3$ over $\mathbb{F}_2$:} when nonscalar coverage uses all $n$ factors,
we construct the identity in $n$ factors using a rotation word.
  \item \textbf{Dimension four over $\mathbb{F}_2$:} we use the exceptional isomorphism
$\mathrm{SL}_4(\mathbb{F}_2) \cong A_8$ and Dvir's class-power covering theorem
\cite{Dvir1985}.
  \item \textbf{Dimension four over $\mathbb{F}_3$:} we analyse the six quadratic seed types in
$\mathrm{SL}_4(\mathbb{F}_3)$ directly, accounting for the distinction between
$\mathrm{SL}_4(\mathbb{F}_3)$- and $\mathrm{GL}_4(\mathbb{F}_3)$-conjugacy.
\end{itemize}

Outside the small-field cases handled directly at their full length, once every
nonscalar matrix is reachable in $n$ factors, only the central scalars remain,
and they turn out to be the hardest targets to build. Nonidentity scalars
$\zeta I$ ($\zeta \in \mu_n(F) \setminus \{1\}$) commute with every matrix and are
therefore invisible to conjugation, yet they are the only targets requiring
$n+1$ transvections. In a product of minimal length $n$, a
target that multiplies every vector in its moving space by a single scalar
$\zeta$ forces every elementary factor to carry the same relative multiplier
(\Cref{sec-two-ended-telescope}), which transvections (whose relative multiplier
is $1$) cannot supply (\Cref{cor-transvection-lengths},
\Cref{cor-sl-diameter-lower-bound}). Conversely, multiplying any scalar target by
the inverse of a suitable conjugate of the seed reduces it to a nonscalar
(\Cref{lem-one-letter-completion}). Scalars therefore require at most one
additional factor, and attain that bound precisely when $\mathrm{SL}_n(F)$ contains
a nonidentity scalar, that is, when $\mu_n(F) \neq  \{1\}$.

Consequently, although all factors and conjugating matrices must have entries
in $F$, the field enters \Cref{thm-sl-diameter} solely through the condition
$\mu_n(F) \neq  \{1\}$. \Cref{tab-intro-fields} summarises the resulting diameter across
standard families of fields: the rationals $\mathbb{Q}$, reals $\mathbb{R}$, complex
numbers $\mathbb{C}$, finite fields $\mathbb{F}_q$, $p$-adic fields $\mathbb{Q}_p$, rational
function fields $k(t_1, \dots, t_r)$ in independent variables, and
algebraically closed fields of positive characteristic. We refer to
\Cref{sec-sl-field-examples} for the derivations and worked interpretations,
including rational families of factorisations over function fields and exact
local factorisations over $\mathbb{Q}_p$.

\begin{table}[htbp]
  \centering
  \begin{tabular}{@{}ll@{}}
    \toprule
    \textbf{Field $F$} & \textbf{When the diameter is $n+1$} \\
    \midrule
    $\mathbb{Q}$ or $\mathbb{R}$ & $n$ is even \\
    $\mathbb{C}$ & Always \\
    $\mathbb{F}_q$ & $\gcd(n, q-1) > 1$ \\
    $\mathbb{Q}_p$, $p$ odd & $\gcd(n, p-1) > 1$ \\
    $\mathbb{Q}_2$ & $n$ is even \\
    $k(t_1, \dots, t_r)$ & $\mu_n(k) \neq \{1\}$ \\
    Algebraically closed, characteristic $p > 0$ & $n$ is not a power of $p$ \\
    \bottomrule
  \end{tabular}
  \caption{Specialisations of \Cref{thm-sl-diameter} for $n \ge 3$. In each case $\Delta(\mathrm{SL}_n(F)) = n+1$ when the stated condition holds and $n$ otherwise.}
  \label{tab-intro-fields}
\end{table}

\subsection{Determinant fibres}
\label{sec-intro-fibres}
We now turn from $\mathrm{SL}_n(F)$ to the full general linear group $\mathrm{GL}_n(F)$,
which decomposes into \emph{determinant fibres}
$\mathcal{D}(a) = \{X \in \mathrm{GL}_n(F) : \det X = a\}$ ($a \in F^\times$), each a coset of
$\mathrm{SL}_n(F)$. Over an infinite field, no finite set of seeds generates
$F^\times$, so $\Delta(\mathrm{GL}_n(F)) = \infty$ (\Cref{cor-gl-infinite-field});
that infinite value reflects only the determinant constraint. To ask whether
any further obstruction remains once the determinant is accounted for, fix a
single noncentral conjugacy class $C$ in $\mathrm{GL}_n(F)$ of determinant $d$ and
consider its exact \emph{positive} powers $C^k$, consisting of products of exactly
$k$ members of $C$ without inverses. Every element of $C^k$ lies in the
permitted fibre $\mathcal{D}(d^k)$. Following Vaserstein and Wheland
\cite[\S~2]{VasersteinWheland1995}, the \emph{abelianised covering number}
$\operatorname{cn}_\mathrm{ab}(\mathrm{GL}_n(F))$ is the least $k$ such that $C^k = \mathcal{D}(d^k)$
for every noncentral full conjugacy class $C$ (\Cref{sec-covering-definitions}).
Our second main theorem shows that this covering number coincides with the
special linear diameter over every field.

\begin{introresult}{thm-gl-fibre-covering}{Determinant-fibre covering}
Let $F$ be a field and $n\ge 3$. Then
  \[
  \operatorname{cn}_\mathrm{ab}(\mathrm{GL}_n(F)) = \Delta(\mathrm{SL}_n(F))
    = n + [\mu_n(F) \neq  \{1\}].
\]
  At this common length every noncentral conjugacy class in $\mathrm{GL}_n(F)$ fills
  its determinant fibre, and the transvection class attains the supremum.
\end{introresult}

\Cref{thm-gl-fibre-covering} concerns full $\mathrm{GL}_n(F)$-conjugacy classes. Within
$\mathrm{SL}_n(F)$ a conjugacy class can be strictly smaller, because only changes of
basis of determinant one are permitted; nevertheless, exact positive powers of
any single noncentral $\mathrm{SL}_n(F)$-conjugacy class still cover $\mathrm{SL}_n(F)$ at
the same length.

\begin{introresult}{cor-sl-class-covering}{Covering by one special linear class}
Let $F$ be a field and $n\ge 3$, and write $k = n + [\mu_n(F) \neq  \{1\}]$. Then
  $C^k = \mathrm{SL}_n(F)$ for every noncentral conjugacy class $C$ of $\mathrm{SL}_n(F)$.
  No smaller exponent suffices for the transvection class.
\end{introresult}

Returning to \Cref{thm-gl-fibre-covering}, the equality
$\operatorname{cn}_\mathrm{ab}(\mathrm{GL}_n(F)) = \Delta(\mathrm{SL}_n(F))$ does not follow automatically
from the fact that each determinant fibre is a copy of $\mathrm{SL}_n(F)$. Two
genuine obstacles separate the two problems.

First, a seed in $\mathrm{GL}_n(F)$ can have any determinant $d \in F^\times$. The
natural way to reduce to $\mathrm{SL}_n(F)$ would be to divide the seed by an $n$th
root of $d$, but over an arbitrary field $F$ no such root need exist. Thus a
covering theorem for $\mathrm{SL}_n(F)$ cannot simply be rescaled to cover the fibre
$\mathcal{D}(d^k)$.

Second, exact positive covering cannot skip the scalar target that plays the
role of the identity. In the diameter problem for $\mathrm{SL}_n(F)$, once every
nonscalar matrix is reachable in $n$ factors and $I$ is the only scalar, the
proof is finished because $I$ is the empty word of length zero. In the fibre
$\mathcal{D}(d^n)$, the scalar matrices are $d \zeta I$ for $\zeta \in \mu_n(F)$, with
$d I$ playing the role of $I$. When $\mu_n(F) = \{1\}$, we are not allowed to
stop at length zero: we must actually write $d I$ as a product of $n$ positive
conjugates of the seed.

The proof of \Cref{thm-gl-fibre-covering} overcomes these two obstacles according
to the size of the field:

\begin{itemize}
  \item \textbf{Fields with $|F|\ge 4$ (\Cref{prop-gl-determinant-scalar}):} Lev's nonscalar
coverage theorem \cite[Theorem 2]{Lev1996} already works at arbitrary determinant
$d$, leaving only the scalar $d I$ when $\mu_n(F) = \{1\}$. Rank-one seeds
are settled by the rank-one criterion of \Cref{sec-intro-rank-one}, and dimensions
$n\ge 5$ by Lev's scalar coverage theorem \cite[Theorem 4(3)]{Lev1996}. In
dimensions $3$ and $4$, we build $d I$ with a \emph{rotation word}
(\Cref{lem-rotation-identity}): conjugating a seed $c$ by successive powers of
a single matrix $p$ of order $n$ causes adjacent powers of $p$ to cancel,
leaving $(c p)^n$. Choosing $p$ so that $(c p)^n = d I$ (via B\"unger and
Nielsen's cyclic decomposition theorem \cite[Theorem 1.2]{BungerNielsen1999} or
Lev's cyclic product theorem \cite[Theorem 5]{Lev1996}) gives $d I$ in $n$
factors. When $\mu_n(F) = \{1\}$, the conjugator $p$ automatically has
determinant one, so all $n$ conjugates lie in a single $\mathrm{SL}_n(F)$-class,
proving \Cref{cor-sl-class-covering} at the same time.
  \item \textbf{The field $\mathbb{F}_3$ (\Cref{prop-gl-f3-fibre}, \Cref{sec-fibre-f3}):}
Kn\"uppel and Nielsen's theorem applies only inside $\mathrm{SL}_n(\mathbb{F}_3)$. In odd
dimensions, $-1 = (-1)^n$ is an $n$th power, so rescaling reduces the
nonscalars to $\mathrm{SL}_n(\mathbb{F}_3)$; we then construct the remaining scalar
$d I$ either by rotation (for nonprimary seeds) or, for primary seeds, by
multiplying commuting conjugates whose scalar parts multiply to the field
norm $d$ while their unipotent parts cancel. In even dimensions with
$\det C = -1$, rescaling into $\mathrm{SL}_n(\mathbb{F}_3)$ is impossible because $-1$
has no even root in $\mathbb{F}_3$. Here, however, the target fibre at length
$n+1$ has determinant $-1$ and contains no scalar matrices at all, so
covering its nonscalars by our arbitrary-determinant cyclic construction
fills the fibre outright.
  \item \textbf{The field $\mathbb{F}_2$ (\Cref{sec-fibre-f2}):} here
$\mathrm{GL}_n(\mathbb{F}_2) = \mathrm{SL}_n(\mathbb{F}_2)$, and the positive constructions of
\Cref{sec-intro-sl} already produce the identity in $n$ factors.
\end{itemize}

\subsection{Prescribed rank-one products}
\label{sec-intro-rank-one}
The preceding theorems determine worst-case lengths across all noncentral
seeds. A complementary linear-algebraic factorisation question prescribes the
elementary factors individually: given a sequence $A_1, \dots, A_k$ of
transformations of \emph{rank-one deviation} and a target $X \in \mathrm{GL}_n(F)$, when
can conjugates of $A_1, \dots, A_k$ multiply to $X$?

A matrix has rank-one deviation at a nonzero scalar \emph{origin} $b$ if it agrees
with $b I$ on a hyperplane and differs from it along a single line;
transvections, reflections, and one-direction dilations are the examples at
origin $1$. For $n\ge 3$, the conjugacy class of such a matrix is determined by
its origin $b$ and its \emph{relative ratio}, the ratio of its determinant to $b^n$
(\Cref{sec-rank-one-criterion}). Our third main theorem answers the factorisation
question completely.

\begin{introresult}{thm-prescribed-rank-one}{Rank-one product criterion}
Let $F$ be a field and $n\ge 3$. Take $A_1, \dots, A_k \in \mathrm{GL}_n(F)$ ($k\ge 1$),
  each of rank-one deviation at an origin $b_i \in F^\times$, and let
  $a = b_1 \dots b_k$. A target $X \in \mathrm{GL}_n(F)$ is a product of conjugates of
  $A_1, \dots, A_k$ if and only if $\det X = \det A_1 \dots \det A_k$ and
  $k \ge  \operatorname{rank}(X - a I)$, subject to two boundary restrictions:
\begin{enumerate}[label=(\roman*)]
  \item when $k = \operatorname{rank}(X - a I) > 0$ and $X$ acts as multiplication by a
      single scalar $\beta$ on the image of $X - a I$, every factor $A_i$ must
      have relative ratio $\beta / a$;
  \item when $X = a I$ (so $\operatorname{rank}(X - a I) = 0$), we must have $k \ge  2$.
\end{enumerate}
\end{introresult}

Specialising all factors to transvections (origin $1$ and relative ratio $1$)
immediately yields the exact transvection length of every matrix in
$\mathrm{SL}_n(F)$ (\Cref{cor-transvection-lengths}, recovering Ellers and Lausch's
transvection-length formula \cite[Proposition 7]{EllersLausch1989}) and the common
lower bounds (\Cref{cor-sl-diameter-lower-bound}) for \Cref{thm-sl-diameter} and
\Cref{thm-gl-fibre-covering}.

\subsection{Neighbouring invariants}
\label{sec-intro-invariants}
Several neighbouring invariants appear in the literature: the worst
single-seed diameter $\Delta_1$, signed distance within an abelianisation
fibre, the classical covering number $\operatorname{cn}_\mathrm{class}$ (over all nonidentity
classes), the normally generating covering number $\operatorname{cn}_\mathrm{ng}$ (repeating
one normally generating class), the extended covering number $\operatorname{ecn}$
(prescribing the classes independently at each position), and $\operatorname{cn}_\mathrm{ab}$.
\Cref{sec-invariants-literature} (\Cref{tab:word-conventions}) gives definitions
for all of these invariants, surveys the known values from the literature, and
uses our results to complete the comparison for $\mathrm{SL}_n(F)$ and $\mathrm{GL}_n(F)$:

\begin{itemize}
  \item \textbf{For $\mathrm{SL}_n(F)$ ($n\ge 3$):} every finite normally generating set contains a
noncentral seed, so
$\Delta_1 = \Delta = \operatorname{cn}_\mathrm{ng} = \operatorname{cn}_\mathrm{ab} = n + [\mu_n(F) \neq  \{1\}]$
over every field, and the signed fibrewise diameter has the same value
(\Cref{prop-comparison-sl}). By
contrast, $\operatorname{cn}_\mathrm{class} = \infty$ whenever $\mu_n(F) \neq  \{1\}$ because
central classes cannot cover, while prescribing classes independently is
stated by Kn\"uppel and Nielsen \cite[p. 635]{KnuppelNielsen2006} to give
$\operatorname{ecn}(\mathrm{SL}_n(F)) = n+1$ for $n\neq 4$ and $F \neq  \mathbb{F}_2$ (see
\Cref{sec-intro-prior} and \Cref{sec-invariants-literature} regarding the printed
proof).
  \item \textbf{For $\mathrm{GL}_n(F)$ ($n\ge 3$, $F \neq  \mathbb{F}_2$):} exact positive whole-group
covering always fails ($\operatorname{cn}_\mathrm{class} = \operatorname{cn}_\mathrm{ng} = \operatorname{ecn} = \infty$),
and $\Delta_1$ and $\Delta$ are finite if and only if $F$ is finite. Over every
field, however, $\operatorname{cn}_\mathrm{ab}(\mathrm{GL}_n(F)) = \Delta(\mathrm{SL}_n(F))$ and the
signed fibrewise diameter is at most $2(n+1)$ (\Cref{prop-comparison-gl}).
\end{itemize}

\subsection{Related work and contributions}
\label{sec-intro-prior}
We summarise what was previously known and what this paper contributes for
each of our three main results:

\begin{itemize}
  \item \textbf{Special linear diameter (\Cref{thm-sl-diameter}):} K\k{e}dra, Libman and Martin
proved the linear bounds $\Delta(\mathrm{SL}_n(F)) \le  12(n-1)$ over every field
\cite[Theorem 6.3]{KedraLibmanMartin2023} and $4(n-1)$ over algebraically closed
fields \cite[Theorem 1.3]{KedraLibmanMartin2022}. Narrower bounds in small
dimensions and over finite fields include Aseeri's bounds
$3 \le  \Delta(\mathrm{SL}_3(\mathbb{R})) \le  4$ and $3 \le  \Delta(\mathrm{SL}_3(\mathbb{C})) \le  6$
(reported in \cite[\S~6]{KedraLibmanMartin2022}), Turfah Alshahaitan's refinement
$3 \le  \Delta(\mathrm{SL}_3(\mathbb{C})) \le  4$ \cite{AlshahaitanPrivComm}, and Rekv\'enyi's exact
projective widths \cite[Theorem 4.5]{Rekvenyi2022}, which give
$\Delta(\mathrm{SL}_3(\mathbb{F}_2)) = \Delta(\mathrm{SL}_3(\mathbb{F}_3)) = 3$ and, together with
the $q\ge 4$ bounds there, $\Delta(\mathrm{SL}_n(\mathbb{F}_q)) = n$ for centreless finite
groups ($\gcd(n,q-1)=1$). In dimension two, Shaun Isherwood
\cite[Theorems 1.1.5--1.1.8]{Isherwood2026HeldDraft} determined $\Delta(\mathrm{SL}_2(F))$
over finite and algebraically closed fields, also locating the faulty step in
Vaserstein and Wheland's $q=5$ argument \cite{VasersteinWheland1995} (whose
counterexample was noted in \cite[p. 635, fn. 1]{KnuppelNielsen2006}). For
$|F|\ge 4$, the signed formula of \Cref{thm-sl-diameter} already follows by
combining Lev's nonscalar coverage theorem \cite[Theorem 2]{Lev1996} with
one-letter completion and classical transvection lengths; our contribution is
the unified deviation-rank framework (\Cref{sec-deviation-geometry}), the
small-field proofs over $\mathbb{F}_2$ and $\mathrm{SL}_4(\mathbb{F}_3)$
(\Cref{sec-small-field-proofs}), and the exact all-field statement covering both
trivial and nontrivial centres.
  \item \textbf{Determinant fibres and single-class covering (\Cref{thm-gl-fibre-covering},
\Cref{cor-sl-class-covering}):} covering numbers of conjugacy classes go back
to Arad and Herzog \cite{AradHerzog1985} (see also \cite{EllersGordeevHerzog1999},
\cite{LiebeckShalev2001}, \cite{GordeevSaxl2002}, and \cite{LiebeckSimion2023}), with the
determinant-fibre convention introduced by Vaserstein and Wheland
\cite[\S~2]{VasersteinWheland1995} and a neighbouring full-similarity length problem
studied by Kn\"uppel \cite{Knuppel1996}. Over broad fields ($|F|\ge 4$), Lev
\cite[Theorems 1, 2, 4]{Lev1996} covered nonscalars for $n\ge 3$, scalars for $n\ge 5$,
and proved $\overline{C}^3 = \mathrm{PSL}_3(F)$ for finite or algebraically closed
$F$ (see also Orevkov \cite[Corollary 1.9]{Orevkov2013} for finite fields), while
Ambrosiewicz \cite[Corollary 6.1]{Ambrosiewicz1996} proved
$\overline{C}^4 = \mathrm{PSL}_3(F)$ for infinite $F$. Kn\"uppel and Nielsen
\cite[p. 635]{KnuppelNielsen2006} state $\operatorname{ecn}(\mathrm{SL}_n(F)) = n+1$ for
$F \neq  \mathbb{F}_2$ and $n\neq 4$, and in \cite[Theorem 1.1]{KnuppelNielsen2013} gave
$n$-letter nonscalar coverage (and, immediately before that theorem,
full-fibre coverage at $n+1$) for $\mathrm{GL}_n(F)$ when $|F|\ge 4$. We prove
\Cref{thm-gl-fibre-covering} and \Cref{cor-sl-class-covering} over every field $F$
and $n\ge 3$, supplying the $n$-factor determinant-scalar rotation in
dimensions $3$ and $4$ (\Cref{prop-gl-determinant-scalar}) and the complete
$\mathbb{F}_2$ and $\mathbb{F}_3$ fibre proofs (\Cref{sec-small-field-proofs}). In dimension
three, where Lev \cite[p. 61]{Lev1996} left open whether the projective identity
lies in $\overline{C}^3$ for irreducible classes over general infinite fields,
\Cref{prop-gl-determinant-scalar} yields the stronger identity $I \in C^3$ in
$\mathrm{SL}_3(F)$ whenever $|F|\ge 4$ (see the remark at the end of
\Cref{sec-determinant-scalar-construction}).
  \item \textbf{Prescribed rank-one products (\Cref{thm-prescribed-rank-one}):} for rank-one
factors at origin $1$, Ellers \cite{Ellers1977} and Ellers and Lausch
\cite[Theorem 5, Propositions 6--7]{EllersLausch1989} proved existence at minimal
length (and with one extra factor in the scalar case) together with the
transvection-length formula. Kn\"uppel and Nielsen
\cite[Observation 7.2, Propositions 7.4--7.5, Corollary 7.6]{KnuppelNielsen2006}
stated the normalised existence assertions for prescribed simple classes,
including at lengths exceeding the deviation rank (omitting only the
one-letter identity exception $k=1$, $X=I$ in Proposition 7.5, also present
in \cite[Proposition 2.3]{KnuppelNielsen2013}). We provide the geometric necessity
framework (\Cref{sec-deviation-geometry}), the explicit one-letter boundary, a
self-contained longer-length proof by factor splitting and ratio grouping,
independent scalar origins, and the resulting complete criterion of
\Cref{thm-prescribed-rank-one}.
\end{itemize}

Three local issues in the published proofs we cite require explicit care. First,
in Kn\"uppel and Nielsen's proof of the extended covering theorem
\cite{KnuppelNielsen2006}, the final case of Proposition 8.2 omits one class from the
product it constructs \cite[p. 654, step (vi)]{KnuppelNielsen2006}, while Propositions
7.4--7.5 and Corollaries 7.6 and 7.10 are deduced via an abbreviated pointer to
Lemmas 7.3 and 7.7 \cite[p. 649]{KnuppelNielsen2006}; we use their theorem only for
$n$ copies of one class together with one further class, where Proposition
8.2(vi) does not arise (\Cref{sec-repeated-class-coverage}). Second, the printed proof
of \cite[Lemma 6.10(c)]{KnuppelNielsen2006} (that two quadratic classes in
$\mathrm{SL}_4(\mathbb{F}_3)$ have a cyclic product) treats only repeated-root quadratics
and omits the distinct-root cases; our proof of \Cref{prop-minimum-rank-product}
avoids that lemma by constructing the quadratic pair directly. Third, the
printed proofs of Theorems 1.2 and 4.2 in B\"unger and Nielsen
\cite[pp. 47--49]{BungerNielsen1999} (cited directly and via
\cite[Propositions 4.2 and 4.8]{KnuppelNielsen2006}) contain three repairable
technical errors in determinant signs, zero-trace completions, and an inverted
constant term; we record the short corrections in the remark following
\Cref{prop-cyclic-power} in \Cref{sec-small-field-cyclic-constructions}.
Finally, for dimension four over $\mathbb{F}_2$, we rely on Dvir's covering theorem
for $A_8$ \cite{Dvir1985} as quoted by Lev \cite[pp. 60--61]{Lev1996}.

\subsection{Organisation of the paper}
\label{sec-intro-organisation}
The body of Part I is organised as follows:

\begin{itemize}
  \item \textbf{\Cref{sec-questions}} fixes the core word-length and covering definitions,
classifies normal generation, and introduces determinant-fibre covering.
  \item \textbf{\Cref{sec-deviation-geometry}} develops deviation rank, the two-ended telescope,
and scalar action at equality.
  \item \textbf{\Cref{sec-rank-one}} proves \Cref{thm-prescribed-rank-one} and deduces transvection
lengths and the common lower bounds.
  \item \textbf{\Cref{sec-class-products}} establishes the upper bounds and proves
\Cref{thm-sl-diameter}, \Cref{thm-gl-fibre-covering}, and
\Cref{cor-sl-class-covering}.
  \item \textbf{\Cref{sec-invariants-and-fields}} compares the neighbouring invariants
(\Cref{tab:word-conventions}) and evaluates the field examples.
  \item \textbf{\Cref{sec-small-field-proofs}} supplies the small-field proofs over $\mathbb{F}_2$ and
$\mathbb{F}_3$.
\end{itemize}

Two problems lie outside Part I and are taken up in Part II: dimension two,
extending Shaun Isherwood's results for $\mathrm{SL}_2(F)$ \cite{Isherwood2026HeldDraft}
from finite and algebraically closed to arbitrary fields (where word lengths
detect quadratic arithmetic), and the whole-group conjugation diameter
$\Delta(\mathrm{GL}_n(\mathbb{F}_q))$ over finite fields, where determinant movement
interacts with fibre filling.

\subsection*{Acknowledgements}
\label{sec-acknowledgements}
I thank Shaun Isherwood, whose PhD thesis inspired this work, for our
conversations and for his assistance with tracking down references, and
Benjamin Martin for his helpful comments particularly on covering numbers.
I am also grateful to Thomas Wright for listening patiently as ideas took
shape, reading drafts, and for his support throughout.

\section{Conjugacy words, diameter and fibre covering}
\label{sec-questions}

\subsection{Signed words and conjugation diameter}
\label{sec-diameter-definitions}
\begin{figure}[htbp]
  \centering
  \includegraphics[width=368pt]{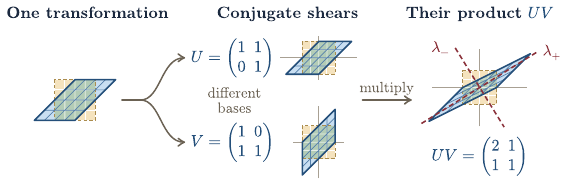}
  \caption{Two conjugates $U$ and $V$ of a single shear multiply (applying
    $V$ first) to a diagonalisable matrix $U V$ with expanding and contracting
    eigenlines. The $2 \times 2$ blocks act in $\mathrm{SL}_3(\mathbb{R})$ by fixing the third
    coordinate.}\label{fig-conjugate-shears}
\end{figure}

Multiplying a matrix by itself keeps all powers in the same basis and
preserves much of its structure: for instance, every power of a shear
(defined in \Cref{sec-shears-and-dilations}) remains unipotent. Conjugation changes
the basis in which a transformation is represented while preserving its
similarity type (its determinant, eigenvalues and Jordan form), whereas
choosing bases independently across successive factors allows multiplication
to change that type while respecting multiplicative invariants such as the
determinant: as \Cref{fig-conjugate-shears} shows for a shear, multiplying two
conjugates already creates expanding and contracting eigendirections.

Starting from a finite collection of matrices, we form products of conjugates
of those matrices and their inverses, and ask how many factors are required to
express an arbitrary target in the group.

\begin{definition}[{Conjugacy words and normal generation}]
\label{def-conjugacy-words}
Let $G$ be a group and $S$ a finite subset of $G$, whose members we call
  \emph{seeds}. The set of \emph{letters} associated with $S$ is
  \[
  L_S = \{h s^\varepsilon h^{-1} : s \in S, h \in G, \varepsilon \in \{-1,1\}\}.
\]
  A \emph{word of length $j$} is an expression $\ell_1 \dots \ell_j$ with each
  $\ell_i \in L_S$, and $L_S^j$ denotes the set of elements represented by words
  of exactly $j$ letters, with $L_S^0=\{1\}$. The set $S$ \emph{normally generates}
  $G$ if every element of $G$ is represented by a word.
\end{definition}

Each conjugate $h s^\varepsilon h^{-1}$ counts as one letter regardless of $h$,
with seed, sign and conjugator chosen independently at each position. The
ambient group $G$ fixes which conjugators are allowed: matrices of determinant
one in $\mathrm{SL}_n(F)$, and arbitrary invertible matrices in $\mathrm{GL}_n(F)$. For a
single seed $g$ with conjugacy class $C$ in $G$, we write
$L_g = C \cup C^{-1}$.

Even when $S$ normally generates $G$, a given element admits many expressions
of different lengths; we measure it by its shortest one.

\begin{definition}[{Word length}]
\label{def-word-length}
The \emph{word length} of $X \in G$ with respect to $S$ is
  \[
  \lVert X \rVert_S = \min \{j \in \mathbb{Z}_{\ge 0} : X \in L_S^j\},
\]
  with $\lVert X \rVert_S=\infty$ if no such word exists. The empty word has length
  zero and represents the identity.
\end{definition}

Thus $\ell \ell^{-1}$ is a two-letter word for the identity, while its word
length is zero.

Normal generation makes every individual target's length finite, yet those
lengths may still be unbounded across $G$. Bounding all targets at once means
asking whether a ball of finite radius contains the whole group.

\begin{definition}[{Balls and diameter for a fixed set}]
\label{def-signed-balls}
For an integer $k\ge 0$, the \emph{ball of radius $k$} is
  \[
  B_S(k) = \{X \in G : \lVert X \rVert_S \le  k\} = \bigcup_{j=0}^k L_S^j,
\]
  and the \emph{diameter with respect to $S$} is
  \[
  \operatorname{diam}_S(G) = \sup_{X \in G} \lVert X \rVert_S,
\]
  allowing the value $\infty$. For $S=\{g\}$ we write $B_g(k)$, or $B_C(k)$
  for the conjugacy class $C$ of $g$; the notation $B_C(k)$ still allows
  letters from both $C$ and $C^{-1}$ and words of at most $k$ letters.
\end{definition}

\begin{figure}[htbp]
  \centering
  \includegraphics[width=368pt]{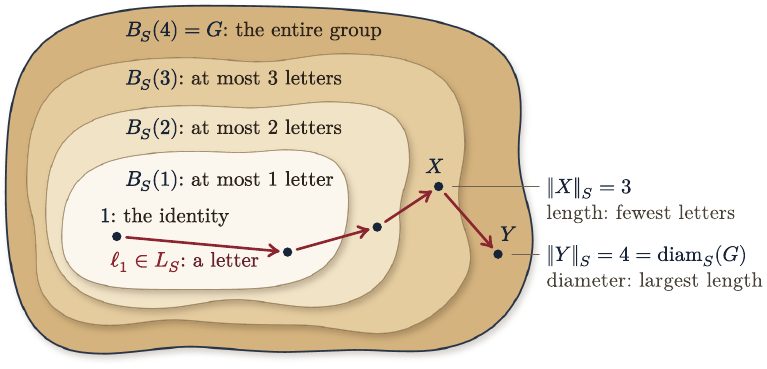}
  \caption{Nested balls $B_S(k)$ and word lengths for a fixed seed set $S$
    (schematic). Each step multiplies by a letter from $L_S$; here $X$ has
    length $3$ and $Y$ attains $\operatorname{diam}_S(G) = 4$.}\label{fig-word-balls}
\end{figure}

Equivalently, $\operatorname{diam}_S(G)\le k$ when $B_S(k)=G$. For a normally
generating set $S$, $\operatorname{diam}_S(G)$ is the diameter of the Cayley graph on
the generating set $L_S$ (\Cref{fig-word-balls}). Because shorter words remain
allowed, the balls $B_S(k)$ are nested, whereas the exact-length sets $L_S^k$
need not be.

To compare word lengths across targets and across different seed sets, we
record three elementary properties of these balls and diameters (compare
\cite[Lemma 2.3]{KedraLibmanMartin2023}): conjugating a target preserves its length
because the letter set $L_S$ is invariant under conjugation, and enlarging the
seed set can only shorten words.

\begin{lemma}[{Basic properties}]
\label{lem-word-length-properties}
Let $S \subseteq T$ be finite subsets of a group $G$. For $X,h \in G$ and
  an integer $k\ge 0$,
\begin{enumerate}
  \item $\lVert h X h^{-1} \rVert_S=\lVert X \rVert_S$;
  \item $B_S(k) \subseteq B_S(k+1)$;
  \item $\lVert X \rVert_T\le \lVert X \rVert_S$, and hence
$B_S(k) \subseteq B_T(k)$ and
$\operatorname{diam}_T(G)\le \operatorname{diam}_S(G)$.
\end{enumerate}
\end{lemma}

\begin{proof}
Conjugating every letter in an expression for $X$ gives an expression of
  the same length for $h X h^{-1}$, and conjugating back by $h^{-1}$ gives
  the reverse inequality. Nesting of balls is immediate from the definition,
  and $L_S \subseteq L_T$ yields the third claim.
\end{proof}

A small diameter for one seed set may reflect a convenient choice of seeds. To
remove this dependence on the seeds, we take the worst diameter over all finite
normally generating sets.

\begin{definition}[{Conjugation diameter}]
\label{def-conjugation-diameter}
If $G$ admits a finite normally generating set, its \emph{conjugation diameter} is
  \[
  \Delta(G) = \sup_S \operatorname{diam}_S(G),
\]
  where $S$ ranges over the finite normally generating subsets of $G$.
  We set $\Delta(G)=\infty$ when these diameters have no finite upper bound,
  and also when $G$ has no finite normally generating set.
\end{definition}

\begin{figure}[htbp]
  \centering
  \includegraphics[width=368pt]{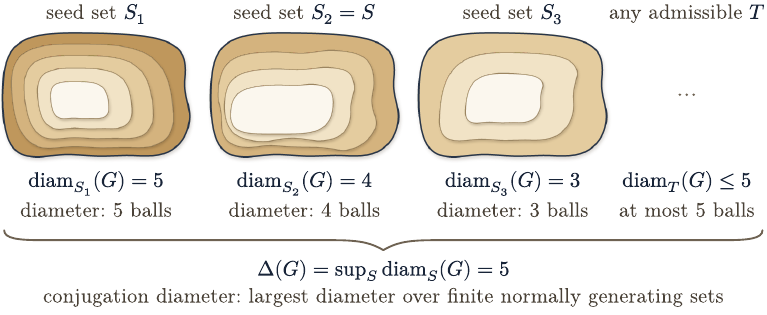}
  \caption{Varying the finite normally generating seed set in the same group
    (schematic), with $S_2 = S$ from \Cref{fig-word-balls}. Different seed sets need
    different numbers of balls to fill $G$; here $\Delta(G) = 5$ is attained by
    $S_1$.}\label{fig-seed-diameters}
\end{figure}

For finitely normally generated groups this is the invariant $\Delta$ of K\k{e}dra,
Libman and Martin \cite[Definition 2.6]{KedraLibmanMartin2023}, termed \emph{conjugation
diameter} by Libman and Tarry \cite{LibmanTarry2020} in their work on symmetric
groups; the convention for groups with no finite normally generating set is our
own. Proving a finite upper bound on $\Delta(G)$ therefore requires showing that
a finite normally generating set exists and bounding $\operatorname{diam}_S(G)$
uniformly across all such sets, whereas a lower bound $\Delta(G)\ge d$ needs only
one admissible seed set and one target of length at least $d$
(\Cref{fig-seed-diameters}).

\subsection{Normal generation in the general linear family}
\label{sec-matrix-normal-generation}
To bound $\Delta(\mathrm{SL}_n(F))$ and $\Delta(\mathrm{GL}_n(F))$, we first determine which
finite subsets normally generate these groups and whether the supremum reduces
to single seeds. Throughout this paper, unless explicitly stated otherwise, $F$
is an arbitrary field and $n\ge 3$; dimension two is treated in Part II.

\subsubsection{Special linear groups}

The centre $Z$ of $\mathrm{SL}_n(F)$ consists of its scalar matrices, which are fixed
by conjugation and generate only scalars. Every other matrix already normally
generates the whole group; because adding further seeds can only shorten words
(\Cref{lem-word-length-properties}), the worst diameter is attained by singletons.

\begin{proposition}[{Normal generation in the special linear group}]
\label{prop-sl-normal-generation}
Let $F$ be a field and $n\ge 3$. A finite subset
  $S \subseteq \mathrm{SL}_n(F)$ normally generates $\mathrm{SL}_n(F)$ if and only
  if it contains a noncentral matrix. In particular, $\mathrm{SL}_n(F)$ admits
  a normally generating singleton over every field, and
  \[
  \Delta(\mathrm{SL}_n(F)) = \sup_{g \in \mathrm{SL}_n(F) \setminus Z}
    \operatorname{diam}_{\{g\}} (\mathrm{SL}_n(F)),
\]
  where $Z$ is the centre of $\mathrm{SL}_n(F)$, consisting of its scalar matrices.
\end{proposition}

\begin{proof}
The group $\mathrm{SL}_n(F)$ is perfect, and its quotient
  $\mathrm{PSL}_n(F)=\mathrm{SL}_n(F) / Z$ is simple; see
  \cite[\S~2, p. 168]{VasersteinWheland1995} for both facts.
  Let $g$ be noncentral and let $N$ be its normal closure in
  $\mathrm{SL}_n(F)$. The image of $N$ in $\mathrm{PSL}_n(F)$ is a nontrivial
  normal subgroup, so simplicity gives $N Z=\mathrm{SL}_n(F)$.
  Thus $\mathrm{SL}_n(F) / N$ is abelian, being a quotient of $Z$,
  and perfect, being a quotient of $\mathrm{SL}_n(F)$. It is therefore
  trivial, proving that $\{g\}$ normally generates $\mathrm{SL}_n(F)$.

  A set of central matrices generates only central matrices, whereas any set
  containing a noncentral $g$ contains the normally generating singleton
  $\{g\}$. Since noncentral matrices exist over every field, this proves the
  criterion and the existence of a normally generating singleton.

  For any finite normally generating $S$, choose a noncentral $g \in S$.
  By \Cref{lem-word-length-properties},
  \[
  \operatorname{diam}_S(\mathrm{SL}_n(F)) \le  \operatorname{diam}_{\{g\}} (\mathrm{SL}_n(F)).
\]
  Taking the supremum over $S$ gives one inequality in the displayed
  formula, and the reverse inequality holds because every noncentral
  singleton is itself a finite normally generating set.
\end{proof}

\subsubsection{General linear groups}

In $\mathrm{GL}_n(F)$, the determinant partitions the group into cosets of
$\mathrm{SL}_n(F)$, the \emph{determinant fibres}
\[
  \mathcal{D}(a) = \{X \in \mathrm{GL}_n(F) : \det X = a\} \quad (a \in F^\times).
\]
Conjugation preserves each fibre, so every word in $L_S$ has its determinant in
the subgroup of $F^\times$ generated by the seed determinants. Once a seed set
contains a noncentral matrix, however, its normal closure already contains
$\mathrm{SL}_n(F)$, so generating $F^\times$ by determinants is the only remaining
condition.

\begin{proposition}[{Normal generation in the general linear group}]
\label{prop-gl-normal-generation}
Let $F$ be a field, $n\ge 3$ and $S$ a finite subset of $\mathrm{GL}_n(F)$.
  Then $S$ normally generates $\mathrm{GL}_n(F)$ if and only if
\begin{enumerate}
  \item $S$ contains a noncentral matrix; and
  \item the determinants $\det s$, for $s \in S$, generate $F^\times$.
\end{enumerate}

  In particular, a singleton $\{g\}$ normally generates $\mathrm{GL}_n(F)$
  if and only if $g$ is noncentral and $\det g$ generates $F^\times$.
\end{proposition}

\begin{proof}
Let $N$ be the normal closure of $S$ in $\mathrm{GL}_n(F)$. If $N=\mathrm{GL}_n(F)$,
  then $S$ cannot consist only of central matrices, and $\det N=F^\times$, where
  $\det N$ is the subgroup generated by $\{\det s : s \in S\}$. Both conditions are
  therefore necessary.

  Conversely, the normal closure in $\mathrm{GL}_n(F)$ of any noncentral matrix
  contains $\mathrm{SL}_n(F)$ \cite[\S~2, p. 168]{VasersteinWheland1995}. A noncentral member
  of $S$ therefore gives $\mathrm{SL}_n(F) \subseteq N$, regardless of its
  determinant. If the seed determinants generate $F^\times$, then for any
  $X \in \mathrm{GL}_n(F)$ there is a $Y \in N$ with $\det Y=\det X$, whence
  $X Y^{-1} \in \mathrm{SL}_n(F) \subseteq N$ and $X \in N$.
\end{proof}

Because $F^\times$ is finitely generated precisely when $F$ is finite, this
determinant condition separates finite from infinite fields.

\begin{corollary}[{The finite-field boundary}]
\label{cor-gl-infinite-field}
Let $F$ be a field and $n\ge 3$. The group $\mathrm{GL}_n(F)$ admits a finite
  normally generating set if and only if $F$ is finite. When $F$ is
  finite, it admits a normally generating singleton. When $F$ is infinite,
  \[
  \Delta(\mathrm{GL}_n(F))=\infty
\]
  under our convention for groups with no finite normally generating set.
\end{corollary}

\begin{proof}
If $F$ is finite with more than two elements, choose a generator $a$
  of $F^\times$; then $\operatorname{diag}(a,1,\dots,1)$ is noncentral and normally
  generates $\mathrm{GL}_n(F)$ by \Cref{prop-gl-normal-generation}. Over $\mathbb{F}_2$,
  any noncentral matrix does so, since $\mathbb{F}_2^\times=\{1\}$.

  Conversely, a finite normally generating set would make $F^\times$
  a finitely generated abelian group, so every subgroup of $F^\times$ would be
  finitely generated. In characteristic zero, however, $F^\times$
  contains $\mathbb{Q}^\times$, which is not finitely generated because a finite
  collection of rational numbers involves only finitely many primes
  in its numerators and denominators.

  In characteristic $p$, if $F$ contains an element $t$ transcendental
  over $\mathbb{F}_p$, then $F^\times$ contains $\mathbb{F}_p(t)^\times$, which is not
  finitely generated because there are infinitely many monic irreducible
  polynomials in $\mathbb{F}_p [t]$. Thus $F$ must be algebraic over $\mathbb{F}_p$.
  Every element of $F^\times$ then lies in a finite subfield and has finite
  order; a finitely generated torsion abelian group is finite, so $F$ is finite.
\end{proof}

\begin{figure}[htbp]
  \centering
  \includegraphics[width=368pt]{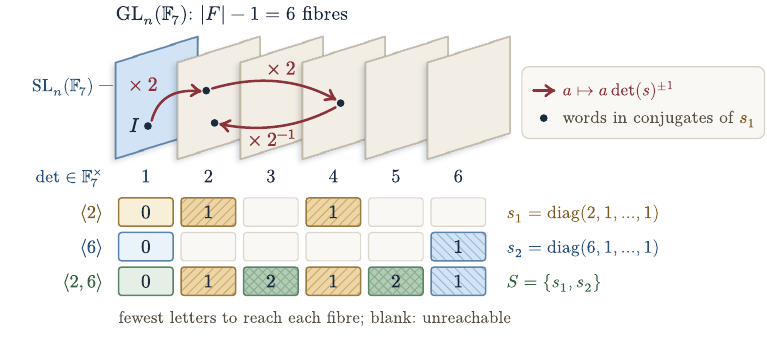}
  \caption{Determinant fibres and collective normal generation in
    $\mathrm{GL}_n(\mathbb{F}_7)$. Neither $\det s_1$ nor $\det s_2$ generates $\mathbb{F}_7^\times$
    alone, yet $\{s_1, s_2\}$ reaches every fibre and normally generates
    $\mathrm{GL}_n(\mathbb{F}_7)$.}\label{fig-determinant-fibres}
\end{figure}

Even over a finite field, where \Cref{cor-gl-infinite-field} guarantees that a
normally generating singleton exists, a prescribed finite normally generating
set need not contain one: several seeds may be needed together to generate
$F^\times$, as \Cref{fig-determinant-fibres} illustrates. Unlike in $\mathrm{SL}_n(F)$, the
singleton reduction in the proof of \Cref{prop-sl-normal-generation} therefore
does not apply to $\mathrm{GL}_n(F)$.

Normal generation in $\mathrm{GL}_n(F)$ thus separates into two problems: reaching
every determinant fibre, and filling each fibre. Setting aside the first leaves
a quantitative question that remains meaningful over every field, even when
$\Delta(\mathrm{GL}_n(F))=\infty$.

\subsection{Abelianised covering number}
\label{sec-covering-definitions}
To study fibre filling on its own, a single noncentral seed is enough: its
normal closure in $\mathrm{GL}_n(F)$ already contains $\mathrm{SL}_n(F)$, and every
determinant fibre $\mathcal{D}(a)$ is a coset of $\mathrm{SL}_n(F)$. To keep all products
of a given length inside a single target fibre, however, we drop the signed
balls of \Cref{sec-diameter-definitions} and use exact positive lengths instead.
Indeed, if $C$ is a conjugacy class of determinant $d$, then allowing inverses
or shorter words would generally mix different determinants, whereas every
product of \emph{exactly} $k$ elements of $C$ has determinant $d^k$ and lands in the
single fibre $\mathcal{D}(d^k)$.

For any conjugacy class $C$ in a group $G$, write
\[
  C^k = \{c_1 \dots c_k : c_1,\dots,c_k \in C\} \quad (k\ge 1),
  \quad C^0 = \{1\}.
\]
In $\mathrm{GL}_n(F)$, we ask for the least positive integer $k$ at which
$C^k = \mathcal{D}(d^k)$ (meaning that determinant is the only remaining constraint
on a $k$-fold product), and how large this length can be as the noncentral
class $C$ varies.

Because the determinant is the abelianisation map of $\mathrm{GL}_n(F)$ with kernel
$\mathrm{SL}_n(F)=[\mathrm{GL}_n(F),\mathrm{GL}_n(F)]$, the same question makes sense in any
group $G$. Conjugation preserves images in $G / [G,G]$, so $C^k$ always
lies in the coset $c^k [G,G]$ for $c \in C$; filling this coset means that the
abelianised image is the sole obstruction to factorisation
\cite[\S~2, pp. 167--168]{VasersteinWheland1995}. Since a class can fill $c^k [G,G]$
only when its normal closure contains $[G,G]$, we restrict the definition to
such classes.

\begin{definition}[{Abelianised covering number}]
\label{def-fibre-covering}
Let $G$ be a group and $C$ a conjugacy class in $G$ whose normal
  closure contains $[G,G]$. For $c \in C$ and $k\ge 0$, put
  \[
  D_k(C) = c^k [G,G].
\]
  The \emph{filling length} of $C$ is the least positive integer $k$ such that
  $C^k=D_k(C)$, or $\infty$ if no such integer exists.

  The \emph{abelianised covering number} $\operatorname{cn}_\mathrm{ab}(G)$ is the
  supremum of these filling lengths over all such classes. Its value
  is $\infty$ if these lengths have no finite upper bound or if there
  is no such class. Conjugacy and normal closure are both taken in $G$.
\end{definition}

In $\mathrm{GL}_n(F)$, every noncentral matrix has $\mathrm{SL}_n(F)$ in its normal closure
(as seen in the proof of \Cref{prop-gl-normal-generation}), so the admissible
classes in \Cref{def-fibre-covering} are precisely the noncentral
$\mathrm{GL}_n(F)$-conjugacy classes, with target fibre $D_k(C) = \mathcal{D}(d^k)$ when
$\det C = d$. In $\mathrm{SL}_n(F)$, the commutator subgroup is the whole group, so
$D_k(C) = \mathrm{SL}_n(F)$ at every length for every noncentral
$\mathrm{SL}_n(F)$-conjugacy class. Both settings therefore cover cosets of
$\mathrm{SL}_n(F)$, yet they are distinct problems: in $\mathrm{GL}_n(F)$ a seed may have
any nonzero determinant, and conjugation uses all invertible matrices rather
than only those of determinant one.

\begin{figure}[htbp]
  \centering
  \includegraphics[width=368pt]{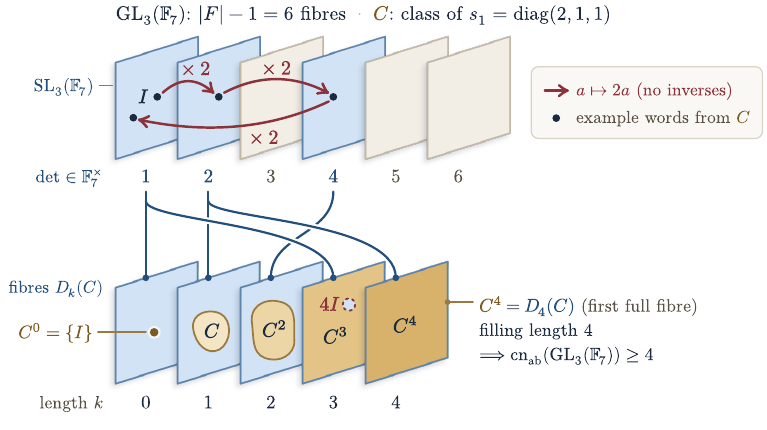}
  \caption{Exact positive powers of the $\mathrm{GL}_3(\mathbb{F}_7)$-class $C$ of
    $s_1 = \operatorname{diag}(2, 1, 1)$, continuing \Cref{fig-determinant-fibres}. Each power
    $C^k$ lies in its target fibre $D_k(C)$, missing only the scalar $4I$ at
    length $3$ before filling $D_4(C)$ at length $4$.}\label{fig-class-powers}
\end{figure}

Unlike the nested balls $B_S(k)$ of \Cref{sec-diameter-definitions}, the powers
$C^k$ move from one fibre to the next as $k$ increases. In \Cref{fig-class-powers},
for instance, products first reach $\mathcal{D}(2)$ at length $1$ and fill it only
on returning at length $4$. Once a power does fill its target fibre, however,
it never loses full coverage at larger lengths.

\begin{lemma}[{Persistence of fibre coverage}]
\label{lem-fibre-persistence}
In the setting of \Cref{def-fibre-covering}, if $C^k=D_k(C)$ for some
  $k\ge 1$, then $C^m=D_m(C)$ for every $m\ge k$.
\end{lemma}

\begin{proof}
Fix $c \in C$. Since $[G,G]$ is normal in $G$,
  \[
  D_k(C)c = c^k [G,G]c = c^{k+1} [G,G] = D_{k+1}(C).
\]
  Hence
  \[
  D_{k+1}(C) = C^k c \subseteq C^{k+1} \subseteq D_{k+1}(C).
\]
  Equality follows, and induction proves the claim.
\end{proof}

\Cref{lem-fibre-persistence} shows that once the filling length is reached, every
subsequent fibre visited by the powers $d^k$ is filled completely. Note,
however, that a class's filling length does not bound the shortest positive
word of every reachable matrix: a matrix missed at an earlier length must wait
until the powers return to its fibre. In \Cref{fig-class-powers}, $C^4$ fills
$\mathcal{D}(2)$, but the scalar $4I \in \mathcal{D}(1)$ missed at length $3$ is only
reached at length $6$, when products next return to $\mathcal{D}(1)$.

Persistence also lets us read $\operatorname{cn}_\mathrm{ab}(G)$ as a single common covering
length rather than just a bound on separate filling lengths. If every
admissible class has filling length at most $k$, then \Cref{lem-fibre-persistence}
guarantees that every such class is full at length $k$ itself; hence, whenever
the admissible family in \Cref{def-fibre-covering} is nonempty and $k\ge 1$,
\[
  \operatorname{cn}_\mathrm{ab}(G)\le k \quad \iff  \quad
  C^k=D_k(C) \;\textnormal{for every admissible class } C.
\]

We now have both of our main invariants in place: the conjugation diameter
$\Delta(\mathrm{SL}_n(F))$, reduced by \Cref{prop-sl-normal-generation} to single
noncentral seeds, and the abelianised covering number
$\operatorname{cn}_\mathrm{ab}(\mathrm{GL}_n(F))$, asking when exact positive powers of a noncentral
class fill their determinant fibre. To establish lower bounds for both, we turn
in the next section to the geometry of elementary factors, beginning with
shears and one-direction dilations.

\section{Deviation rank and products}
\label{sec-deviation-geometry}

\subsection{Transvections, dilations and scalars}
\label{sec-shears-and-dilations}
Shears and one-direction dilations generate all invertible linear
transformations, as Gaussian elimination shows. Yet their simplicity limits
what they can accomplish as seeds: shears, also called \emph{transvections}, will
later give us sharp lower bounds for both the special linear diameter and
determinant-fibre covering. Some of the targets that force long words are
simpler still: scalar transformations, which multiply every vector by the
same number.

\begin{figure}[htbp]
  \centering
  \includegraphics[width=368pt]{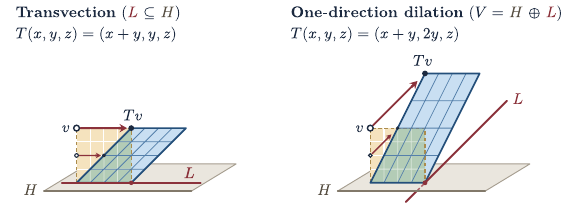}
  \caption{A transvection and a one-direction dilation on a grid in the
    $x y$-plane, with fixed hyperplane $H = \{y = 0\}$ and displacement line $L$.
    On the left \mbox{$L \subseteq H$}; on the right $L$ meets $H$ obliquely (it need
    not be orthogonal to $H$).}\label{fig-shears-dilations}
\end{figure}

Both a transvection and a one-direction dilation fix a hyperplane $H$
pointwise and displace every vector parallel to a single line $L$; what
separates them is whether $L$ lies inside $H$ or outside it
(\Cref{fig-shears-dilations}). In a transvection, $L \subseteq H$, so vectors along
$L$ are fixed and the remaining vectors shear parallel to $H$. In a
one-direction dilation, $L$ is complementary to $H$, so vectors along $L$ are
rescaled while $H$ stays fixed.

For either transformation $T$, choose a nonzero vector $u$ along $L$ and
write the displacement of $v$ as $\varphi(v)u$. Since the displacement map $T-I$
is linear and nonzero, $\varphi$ is a nonzero linear functional, and we can write
\[
  T(v)=v+\varphi(v)u, \quad T=I+u \varphi,
\]
where $u \varphi$ denotes the map $v \mapsto \varphi(v)u$. The functional vanishes
precisely on the fixed hyperplane $H=\ker \varphi$. On the displacement line
$L=F u$, we have $T u=(1+\varphi(u))u$: a transvection fixes this line pointwise
when $\varphi(u)=0$, while a dilation rescales it when $\varphi(u)\neq 0$.

\begin{definition}[{Transvections and one-direction dilations}]
\label{def-elementary-transformations}
Let $V=F^n$ with $n\ge 2$, and let $u \in V$ and $\varphi \in V^*$ be nonzero,
  with $1+\varphi(u)\neq 0$. The map $T=I+u \varphi$ is a \emph{transvection} if
  $\varphi(u)=0$, and a \emph{one-direction dilation} of multiplier
  $\theta=1+\varphi(u)$ if $\varphi(u)\neq 0$.
\end{definition}

\begin{lemma}[{Determinants of the elementary transformations}]
\label{lem-elementary-transformations}
A transvection has determinant $1$. A one-direction dilation of multiplier
  $\theta$ has determinant $\theta$. In both cases the transformation $T$
  satisfies $\operatorname{rank}(T-I)=1$.
\end{lemma}

\begin{proof}
The displacement map $T-I=u \varphi$ has kernel $H$ and image $L$, so its
  rank is one. The map $T$ fixes $H$ pointwise, and $\varphi$ identifies
  $V / H$ with $F$. The identity $\varphi(T v)=(1+\varphi(u))\varphi(v)$ shows that
  $T$ acts on this quotient by multiplication by $1+\varphi(u)$.
  Hence $\det T=1+\varphi(u)$.
\end{proof}

By \Cref{lem-elementary-transformations}, transvections lie in $\mathrm{SL}_n(F)$, and
they generate it by Gaussian elimination (as the elementary row-addition
matrices); adjoining one-direction dilations supplies the determinant changes
needed to generate $\mathrm{GL}_n(F)$. When we use either type as a seed, conjugation by $g$ moves $H$
and $L$ to $g H$ and $g L$, preserving their incidence and the determinant,
while inversion fixes both subspaces and inverts a dilation's multiplier.
Although conjugation lets the displacement line vary from one letter to the
next, every permitted letter still displaces vectors along only a single line.

It is useful to think of each independent displacement direction of a target
as a cost that the letters in a word must pay for. Since each elementary
letter can contribute at most one independent direction, a target whose
displacements span many directions cannot be built from a short word; we make
this rank bound precise in \Cref{sec-two-ended-telescope}.

\begin{figure}[htbp]
  \centering
  \begin{minipage}[c]{172pt}
    \centering
    \includegraphics[width=\linewidth]{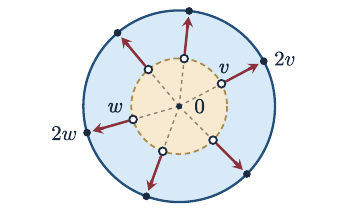}
  \end{minipage}\hfill
  \begin{minipage}[c]{\dimexpr\linewidth-172pt-16pt\relax}
    \captionsetup{justification=raggedright,singlelinecheck=false}
    \caption{Multiplication by $2$ on $\mathbb{R}^2$: every direction is scaled radially,
    doubling the radius of the circle.}\label{fig-uniform-scalar-action}
  \end{minipage}
\end{figure}

Among the most demanding targets are the simplest matrices of all: nonidentity
scalar transformations, which multiply every vector by the same constant
$b\neq 1$ and therefore displace vectors in every direction of $V$
(\Cref{fig-uniform-scalar-action}).

\begin{definition}[{Scalar transformations}]
\label{def-scalar-transformations}
An invertible \emph{scalar transformation} of $V=F^n$ is a map $b I$ with
  $b \in F^\times$, acting by $v \mapsto b v$.
\end{definition}

Indeed, when $b\neq 1$, the displacements $(b-1)v$ span all of $V$, so paying for
all $n$ independent directions takes at least $n$ elementary letters. Moreover,
those letters must combine to scale every direction by the same multiplier $b$.
In a minimal $n$-letter product, this uniformity forces every elementary factor
to act by the same relative multiplier, which transvections (whose multiplier
on their displacement line is $1$) cannot supply; reaching a nonidentity
scalar therefore requires at least $n+1$ transvections, as
\Cref{lem-tight-scalar-action} and \Cref{sec-rank-one} will establish.

Although scalar matrices are excluded as normally generating seeds
(\Cref{sec-matrix-normal-generation}), every scalar matrix in the group is a valid
target. Every $b \in F^\times$ gives a scalar target in $\mathrm{GL}_n(F)$, whereas
membership in $\mathrm{SL}_n(F)$ requires $\det(b I)=b^n=1$. A nonidentity scalar
target therefore exists in $\mathrm{SL}_n(F)$ precisely when $F$ contains a
nontrivial $n$th root of unity, which accounts for the root-of-unity condition
in our diameter and fibre-covering formulas.

\subsection{Deviation from a scalar}
\label{sec-scalar-origins}
Each letter from a transvection or one-direction dilation has only one
displacement direction, which is why targets that displace vectors in many
independent directions require long words in these seeds. For an arbitrary
matrix $A$, the displacements $A v-v$ form $\im(A-I)$, so $\operatorname{rank}(A-I)$
counts their independent directions. More generally, we can compare $A$ with
any nonzero scalar $c I$, over $F$ or an algebraic closure $\overline{F}$.

\begin{definition}[{Deviation from a scalar}]
\label{def-deviation-rank}
Let $A \in \mathrm{GL}_n(F)$ and $c \in \overline{F}^\times$. The \emph{deviation operator}
  of $A$ at $c$ is $A-c I$. Its image $\im(A-c I)$ is the \emph{moving space}
  of $A$ at $c$, and its dimension is the \emph{deviation rank}:
  \[
  \rho_c(A)=\dim \im(A-c I)=\operatorname{rank}(A-c I).
\]
  The \emph{deviation profile} of $A$ is the function $c \mapsto \rho_c(A)$.
  The operator, moving space and rank are taken over a field containing
  $F$ and $c$.
\end{definition}

At $c=1$, the deviation operator is the displacement map from
\Cref{sec-shears-and-dilations}, so $\rho_1(T)=1$ for transvections and
one-direction dilations. For general $c$, the moving space consists of the
differences $A v-c v$, and the kernel consists of the vectors on which $A$
agrees with multiplication by $c$. The kernel is therefore the $c$-eigenspace:
the larger this eigenspace, the smaller the deviation rank.

\begin{lemma}[{Deviation rank and eigenspaces}]
\label{lem-deviation-eigenspaces}
For $A \in \mathrm{GL}_n(F)$ and $c \in \overline{F}^\times$,
  \[
  \rho_c(A)=n-\dim \ker(A-c I).
\]
  In particular, $0\le \rho_c(A)\le n$, and $\rho_c(A)=n$ unless $c$ is an
  eigenvalue of $A$. There are therefore at most $n$ scalar values at which
  the deviation rank is less than $n$.
\end{lemma}

\begin{proof}
The equality follows from rank-nullity for $A-c I$. Its kernel is nonzero
  exactly when $c$ is an eigenvalue, and the characteristic polynomial has
  degree $n$.
\end{proof}

\begin{figure}[htbp]
  \centering
  \includegraphics[width=368pt]{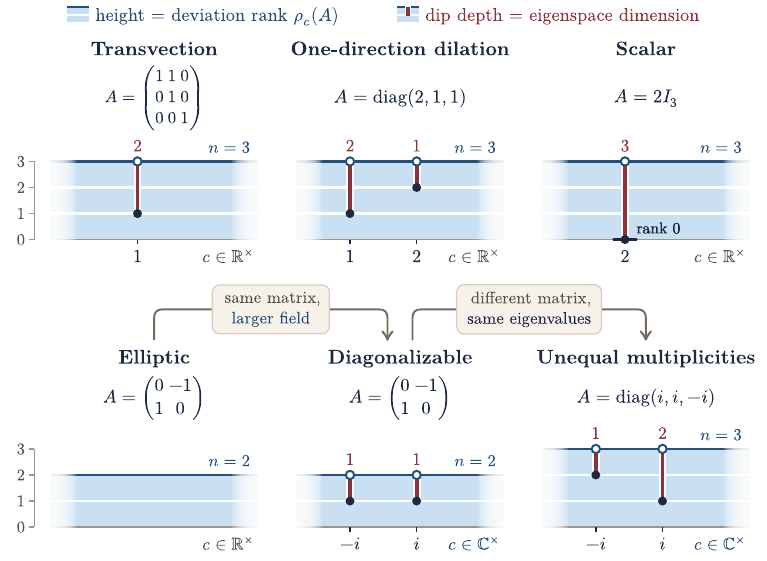}
  \caption{Six deviation profiles $c \mapsto \rho_c(A)$ on a common rank scale
    (horizontal positions schematic). The rank stays at $n$ away from
    eigenvalues, drops by the eigenspace dimension at each eigenvalue, and can
    gain new dips over $\mathbb{C}$.}\label{fig-deviation-profiles}
\end{figure}

\Cref{fig-deviation-profiles} illustrates these deviation profiles across six
examples, showing how each dip depth records an eigenspace dimension and how
passing to an algebraic closure can expose eigenvalues absent over the base
field.

To use deviation rank as a word-length invariant, we need to know how it
behaves across the letters built from a seed $A$. Conjugation moves the image
and kernel of $A-c I$ without changing their dimensions, so deviation rank is
constant on a conjugacy class $C$; we write its common value as $\rho_c(C)$.
Inversion and scalar multiplication also preserve the rank once the comparison
scalar is adjusted accordingly.

\begin{lemma}[{Conjugation, inversion and rescaling}]
\label{lem-deviation-letters}
Let $A,g \in \mathrm{GL}_n(F)$, $b \in F^\times$ and $c \in \overline{F}^\times$.
  Then
  \begin{align*}
  \rho_c(g A g^{-1}) &= \rho_c(A), \\
    \rho_{c^{-1}} (A^{-1}) &= \rho_c(A), \\
    \rho_{b c} (b A) &= \rho_c(A).
\end{align*}
\end{lemma}

\begin{proof}
The corresponding operator identities are
  \begin{align*}
  g A g^{-1}-c I &= g(A-c I)g^{-1}, \\
    A^{-1}-c^{-1}I &= -c^{-1}A^{-1}(A-c I), \\
    b A-b c I &= b(A-c I).
\end{align*}
  Multiplication by an invertible matrix or a nonzero scalar preserves rank.
  Conjugation sends the image and kernel to their images under $g$;
  inversion and rescaling leave the kernel unchanged.
\end{proof}

Thus every letter from $A$ has the same deviation rank $\rho_c(A)$,
measured at $c$ for positive letters and at $c^{-1}$ for inverse letters.
Regardless of which conjugates or inverses are chosen, the seed's deviation
rank therefore fixes how much each letter can contribute towards paying for a
target's deviation.

\begin{figure}[htbp]
  \centering
  \includegraphics[width=\linewidth]{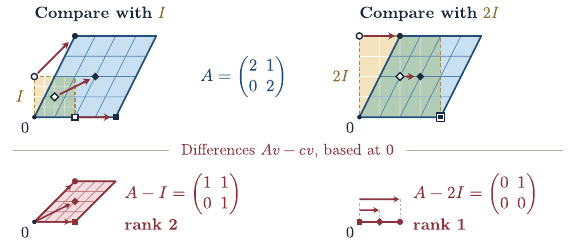}
  \caption{Comparing $A$ with $I$ versus $2I$ on the unit square in $\mathbb{R}^2$.
    Translating the displacements $A v - c v$ to the origin shows that $A - I$
    has rank $2$, whereas $A - 2I$ has rank $1$.}\label{fig-scalar-comparison}
\end{figure}

We consider deviation from a general scalar rather than only the identity,
since comparison with the identity can hide a seed's small deviation from
another scalar, as \Cref{fig-scalar-comparison} illustrates. For such a seed,
each letter pays for little deviation, so reaching a target with large
deviation from the corresponding scalar requires a long word. General scalar
comparisons can therefore give stronger word-length bounds.
The comparison scalars multiply when the matrices multiply, making the
formulas slightly more complex.

\subsection{Deviation of products}
\label{sec-two-ended-telescope}
To turn deviation rank into a word-length bound, we need to see how factor
deviations combine when matrices are multiplied. The case of two matrices
already reveals the mechanism.

\begin{example}[{Deviation of a two-factor product}]
\label{ex-two-factor-deviation}
For $A,B \in \mathrm{GL}_n(F)$, consider their deviations $A-a I$ and $B-b I$,
  with $a,b \in \overline{F}^\times$. The comparison scalars multiply,
  so the product deviation is $A B-a b I$. Adding and subtracting $b A$ gives
  \begin{align*}
  A B-a b I &= (b A-a b I)+(A B-b A) \\
              &= b(A-a I)+A(B-b I).
\end{align*}
  Adding and subtracting $a B$ instead gives
  \[
  A B-a b I=a(B-b I)+(A-a I)B.
\]
  Applied to a vector $v$, these two telescoping identities give the paths
  in \Cref{fig-product-decomposition}.

\begin{figure}[htbp]
  \centering
  \begin{minipage}[c]{224pt}
    \centering
    \includegraphics[width=\linewidth]{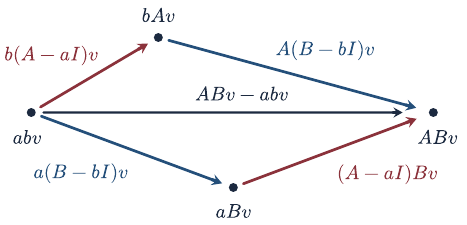}
  \end{minipage}\hfill
  \begin{minipage}[c]{\dimexpr\linewidth-224pt-14pt\relax}
    \captionsetup{justification=raggedright,singlelinecheck=false}
    \caption{Two ways to decompose the total deviation \mbox{$A B v - a b v$}
      into a step from \mbox{$A - a I$} and a step from \mbox{$B - b I$}, passing
      through $b A v$ or $a B v$.}\label{fig-product-decomposition}
  \end{minipage}
\end{figure}

  In the first identity, multiplication by $A$ can move the image of
  $B-b I$, but leaves its kernel unchanged because $A$ is invertible.
  Thus if $A v=a v$ and $B v=b v$, then
  $A B v=a b v$. Equivalently,
  \[
  \ker(A-a I) \cap \ker(B-b I) \subseteq \ker(A B-a b I).
\]
  In the second identity, multiplication on the right by $B$ can move
  the kernel of $A-a I$, but leaves its image unchanged.
  The two summands therefore take their values in the
  images of $B-b I$ and $A-a I$, respectively. Hence
  \[
  \im(A B-a b I) \subseteq \im(A-a I)+\im(B-b I).
\]
  The dimension of this sum is at most $\rho_a(A)+\rho_b(B)$, so
  \[
  \rho_{a b}(A B)\le \rho_a(A)+\rho_b(B).
\]
\end{example}

The two telescoping identities extend to products of any length: expanding from
the left places the prefix products on the left of each factor deviation, while
expanding from the right places the suffix products on the right.

\begin{lemma}[{Two-ended telescope}]
\label{lem-two-ended-deviation}
Let $S_1,\dots,S_k \in M_n(F)$ ($k\ge 1$) and $c_1,\dots,c_k \in \overline{F}$.
  Put $X=S_1 \dots S_k$, $a=c_1 \dots c_k$ and $N_i=S_i-c_i I$. Then
  \begin{align*}
  X-a I &= \sum_{i=1}^k (S_1 \dots S_{i-1}) N_i(c_{i+1} \dots c_k) \\
          &= \sum_{i=1}^k (c_1 \dots c_{i-1}) N_i(S_{i+1} \dots S_k),
\end{align*}
  with empty matrix products equal to $I$ and empty scalar products equal to $1$.
\end{lemma}

\begin{proof}
Put $\Pi_i=S_1 \dots S_i$ and $a_i=c_1 \dots c_i$, with $\Pi_0=I$ and $a_0=1$.
  The deviation of the next partial product satisfies
  \[
  \Pi_i-a_i I
    =c_i(\Pi_{i-1}-a_{i-1}I)+\Pi_{i-1}N_i.
\]
  The first term scales the preceding deviation by $c_i$; the second
  is the new factor deviation multiplied by the preceding partial product.
  To give the partial products a common comparison scalar $a$, multiply
  each by the remaining scalars, writing $Q_i=\Pi_i(c_{i+1} \dots c_k)$.
  The endpoints are $Q_0=a I$ and $Q_k=X$. Multiplying the recurrence
  by $c_{i+1} \dots c_k$ gives
  \begin{align*}
  Q_i-Q_{i-1}
    &= \Pi_{i-1}(S_i-c_i I)(c_{i+1} \dots c_k) \\
    &= \Pi_{i-1}N_i(c_{i+1} \dots c_k).
\end{align*}
  Summing these differences cancels every intermediate $Q_i$, leaving
  $X-a I$ and proving the first identity.

  From the other end, use the suffix products $S_{i+1} \dots S_k$,
  multiplied by the product $a_i$ of the preceding comparison scalars.
  The difference between successive stages is
  \[
  a_{i-1}S_i \dots S_k-a_i S_{i+1} \dots S_k
    =a_{i-1}N_i S_{i+1} \dots S_k.
\]
  Summing for $i=1,\dots,k$ cancels the intermediate terms, leaving
  $X-a I$ and proving the second identity.
\end{proof}

Just as in \Cref{ex-two-factor-deviation}, the left-ended expansion gives a kernel
inclusion and the right-ended expansion gives an image sum, yielding a general
upper bound on the deviation rank of a product.

\begin{corollary}[{Deviation rank, image and kernel}]
\label{cor-product-deviation}
With the notation of \Cref{lem-two-ended-deviation} (extending
  $\rho_c(S)=\operatorname{rank}(S-c I)$ over $\overline{F}$ to all $S \in M_n(F)$ and
  $c \in \overline{F}$),
  \begin{align*}
  \bigcap_{i=1}^k \ker N_i &\subseteq \ker(X-a I), \\
    \im(X-a I) &\subseteq \sum_{i=1}^k \im N_i,
\end{align*}
  and consequently
  \[
  \rho_a(X) \le  \sum_{i=1}^k \rho_{c_i}(S_i), \quad
    \dim \ker(X-a I) \ge  n - \sum_{i=1}^k \rho_{c_i}(S_i).
\]
\end{corollary}

\begin{proof}
A vector in every $\ker N_i$ is killed by every summand in the first
  expansion, giving the kernel inclusion.
  In the second expansion, each summand has image contained in $\im N_i$,
  so their sum has image contained in $\sum_i \im N_i$.
  Taking dimensions gives $\rho_a(X)\le \sum_i \operatorname{rank}N_i=\sum_i \rho_{c_i}(S_i)$,
  and rank-nullity gives the kernel dimension bound.
\end{proof}

Applying \Cref{cor-product-deviation} to a word of $k$ conjugates of $A$ and
$A^{-1}$ makes the payment rule of \Cref{sec-shears-and-dilations} precise. By
\Cref{lem-deviation-letters}, every letter contributes deviation rank $\rho_c(A)$
(at comparison scalar $c$ for positive letters and $c^{-1}$ for inverse
letters), giving a total \emph{rank budget} of $k \rho_c(A)$ that must cover the
target's deviation at the product scalar $c^\sigma$, where $\sigma$ is the number
of positive letters minus the number of inverse letters.

\begin{corollary}[{Budget inequality}]
\label{cor-deviation-word-length}
Let $A \in \mathrm{GL}_n(F)$ and $c \in \overline{F}^\times$. Suppose $X \in \mathrm{GL}_n(F)$
  is a product of $k\ge 0$ conjugates of $A$ and $A^{-1}$ in $\mathrm{GL}_n(F)$, and
  let $\sigma$ be the number of positive factors minus the number of inverse
  factors. Then
  \[
  \rho_{c^\sigma}(X)\le k \rho_c(A).
\]
  In particular, if $C$ is a conjugacy class in $\mathrm{GL}_n(F)$ and $X \in C^k$,
  then $\rho_{c^k}(X)\le k \rho_c(C)$.
\end{corollary}

\begin{proof}
Write $X=\prod_{i=1}^k h_i A^{\varepsilon_i} h_i^{-1}$ with
  $h_i \in \mathrm{GL}_n(F)$ and $\varepsilon_i \in \{-1,1\}$, so
  $\sigma=\sum_{i=1}^k \varepsilon_i$. Give the $i$th letter comparison scalar
  $c_i=c^{\varepsilon_i}$. By \Cref{lem-deviation-letters}, its deviation rank is
  $\rho_c(A)$, and the comparison scalars multiply to $c^\sigma$, so
  \Cref{cor-product-deviation} gives the bound. For a positive class power all
  signs are $+1$, so $\sigma=k$; the empty word has $X=I$ and $\sigma=0$, and
  satisfies the bound as well.
\end{proof}

\begin{remark}[{Relation to the rank-one case}]
For invertible factors whose deviation from $I$ has rank one, the
  inclusions in \Cref{cor-product-deviation} are the moved-space and fixed-space
  containments of delMas and Lewis \cite[Proposition 6]{DelMasLewis2022}.
  The telescope here also allows arbitrary deviation ranks, singular factors
  and zero comparison scalars. Its proof uses neither matrix inverses nor
  division by scalars.
\end{remark}

To extract sharp lower bounds from \Cref{cor-deviation-word-length}, we want seeds
whose letters contribute as little deviation rank as possible. Transvections
have the smallest positive deviation rank, $\rho_1(T)=1$, and because
$1^{-1}=1$ their comparison scalar stays fixed at $1$ for both positive and
inverse letters. Every product of $k$ transvections therefore satisfies
$\rho_1(X)\le k$. In particular, any matrix $X$ with no eigenvalue $1$ has no
nonzero fixed vectors, hence $\rho_1(X)=n$, and requires at least $n$
transvections. Such targets exist in $\mathrm{SL}_n(F)$ (and therefore in
$\mathrm{GL}_n(F)$) over every field.

\begin{example}[{Full deviation rank over every field}]
\label{ex-full-deviation-target}
For $n\ge 2$, let $X$ be the companion matrix with characteristic polynomial
  \[
  p(t)=t^n+b t+(-1)^n, \quad b \in F.
\]
  The constant term gives $\det X=(-1)^n p(0)=1$. To exclude the eigenvalue $1$,
  choose $b$ so that $p(1)=1+b+(-1)^n\neq 0$. Only one value of $b$ is
  forbidden, so every field admits such a choice.
\end{example}

The budget inequality alone gives only a lower bound of $n$ transvections,
because a nonidentity scalar $\zeta I$ ($\zeta \in \mu_n(F) \setminus \{1\}$) has
deviation rank $n$ at $1$. Why, then, does $\zeta I$ require at least $n+1$
transvections? The obstruction appears at \emph{equality} in
\Cref{cor-product-deviation}, where the factors must supply the target's
deviation rank with no overlap or cancellation.

We call a product $X=S_1 \dots S_k$ \emph{tight} at comparison scalars
$b_1,\dots,b_k$ (with product $a=b_1 \dots b_k$) if
$\rho_a(X)=\sum_{i=1}^k \rho_{b_i}(S_i)$, so no rank is lost in the telescope
sum. Because $X$ commutes with $X-a I$, it preserves the moving space
$W=\im(X-a I)$, and we call the restriction $X|_W$ the \emph{moving action} of $X$
at $a$. For a nonidentity scalar target $X=\zeta I$ at $a=1$, for instance, we
have $W=F^n$ and the moving action is simply multiplication by $\zeta$.

When a product is tight, the telescope terms in \Cref{lem-two-ended-deviation}
have linearly independent images whose direct sum is $W$. If $X$ acts as a
single scalar $\beta$ on $W$, those independent terms cannot mix at all,
forcing each factor deviation to satisfy the identities below.

\begin{lemma}[{Scalar action at equality}]
\label{lem-tight-scalar-action}
Let $S_1,\dots,S_k \in \mathrm{GL}_n(F)$ and $b_1,\dots,b_k \in F^\times$, and put
  $X=S_1 \dots S_k$, $a=b_1 \dots b_k$ and $N_i=S_i-b_i I$. Suppose the product
  is tight ($\rho_a(X)=\sum_{i=1}^k \operatorname{rank} N_i$) and $X$ acts on its moving
  space $W=\im(X-a I)$ as multiplication by a scalar $\beta \in F$. Then
  \[
  N_i N_j=0 \quad (i\neq j), \quad
    N_i^2=b_i \left(\frac{\beta}{a}-1\right)N_i \quad (1\le i\le k).
\]
\end{lemma}

\begin{proof}
Write $\Pi_i=S_1 \dots S_i$ and $a_i=b_1 \dots b_i$, with $\Pi_0=I$ and
  $a_0=1$. The telescope gives
  \[
  D=X-a I=\sum_{i=1}^k D_i, \quad
    D_i=\left(\frac{a}{a_i}\right)\Pi_{i-1}N_i.
\]
  Each $D_i$ has the same rank as $N_i$. Tightness forces equality in
  \[
  \operatorname{rank} D\le \dim \sum_{i=1}^k \im D_i
    \le \sum_{i=1}^k \operatorname{rank} D_i.
\]
  Thus the images form a direct sum equal to $\im D$, as in the
  rank-additivity argument of Marsaglia and Styan
  \cite[p. 451]{MarsagliaStyan1972}:
  \[
  W=\bigoplus _{i=1}^k W_i, \quad W_i=\im D_i.
\]
  For $v \in W_j$, the scalar action of $X$ gives
  \[
  \sum_{i=1}^k D_i v=(\beta-a)v.
\]
  Each term on the left belongs to its own summand $W_i$, whereas the
  right-hand side belongs to $W_j$. Comparing components therefore gives
  $D_i v=0$ for $i\neq j$ and $D_j v=(\beta-a)v$. Hence
  \[
  D_i D_j=0 \quad (i\neq j), \quad D_i^2=(\beta-a)D_i.
\]

  To recover the factor deviations, apply the telescope to the preceding
  prefix:
  \[
  \Pi_{i-1}=a_{i-1}\left(I+\frac{1}{a} \sum_{h<i}D_h\right).
\]
  Every earlier $D_h$ vanishes on $W_i$, so this prefix acts on $W_i$
  by $a_{i-1}$. Since $D_i$ takes its values in $W_i$,
  \[
  N_i=\left(\frac{a_i}{a}\right)\Pi_{i-1}^{-1}D_i=\left(\frac{b_i}{a}\right)D_i.
\]
  The identities for the $N_i$ now follow from those for the $D_i$.
\end{proof}

Applying \Cref{lem-tight-scalar-action} to a nonidentity scalar target $\zeta I$
($\zeta\neq 1$) resolves the $n+1$ obstruction previewed in
\Cref{sec-shears-and-dilations}. Any $n$-transvection product for $\zeta I$ is tight
at $b_i=1$ (so $a=1$), with moving action $\beta=\zeta$. For each transvection we
have $N_i\neq 0$ and $N_i^2=0$, so the second identity of
\Cref{lem-tight-scalar-action} forces $\beta=1$, contradicting $\zeta\neq 1$. Hence no
$n$ transvections can multiply to $\zeta I$.

In \Cref{sec-rank-one} we specialise these necessary conditions to arbitrary
prescribed sequences of rank-one factors, show that they are also sufficient,
and deduce the exact transvection length of every matrix in $\mathrm{SL}_n(F)$.

\section{Rank-one products and lower bounds}
\label{sec-rank-one}
\Cref{sec-deviation-geometry} gave two lower-bound mechanisms for transvection words:
full deviation rank forces at least $n$ transvections
(\Cref{ex-full-deviation-target}), and \Cref{lem-tight-scalar-action} rules out
$n$-transvection expressions for any nonidentity scalar. To turn these
obstructions into exact transvection lengths, and hence sharp lower bounds for
both $\Delta(\mathrm{SL}_n(F))$ and $\operatorname{cn}_\mathrm{ab}(\mathrm{GL}_n(F))$, we need to know when
these necessary conditions are also sufficient. Because we will also encounter
rank-one factors at arbitrary scalar origins later in the paper, it pays to
settle a broader factorisation question first: given any prescribed sequence of
rank-one transformations in $\mathrm{GL}_n(F)$, which targets can be written as a
product of their conjugates? We establish the complete criterion in
\Cref{sec-rank-one-criterion} and deduce exact transvection lengths and both
lower bounds in \Cref{sec-rank-one-lower-bounds}.

\subsection{The rank-one product criterion}
\label{sec-rank-one-criterion}
In \Cref{sec-shears-and-dilations} we met the simplest rank-one transformations:
transvections and one-direction dilations, both of which have deviation rank $1$
at the identity. General rank-one transformations in $\mathrm{GL}_n(F)$ can live at
any scalar origin $b \in F^\times$, like the scaled shear at $2 I$ in
\Cref{fig-scalar-comparison}. Since our factorisation question allows replacing each
prescribed factor by any $\mathrm{GL}_n(F)$-conjugate, we first classify the conjugacy
classes of rank-one transformations.

\begin{definition}[{Origin and relative ratio}]
\label{def-rank-one-class}
Let $s \in \mathrm{GL}_n(F)$, with $n\ge 2$. If $\rho_b(s)=1$ for
  $b \in F^\times$, we call $s$ a \emph{rank-one factor at origin $b$} and
  define its \emph{relative ratio} to be $\theta=\det(s)/b^n$.
  For $b,\theta \in F^\times$, write
  \[
  C(b,\theta)=\{s \in \mathrm{GL}_n(F) : \rho_b(s)=1, \det s=b^n \theta\}.
\]
\end{definition}

\begin{lemma}[{Rank-one deviation}]
\label{lem-rank-one-classes}
Let $F$ be a field, $n\ge 2$ and $b \in F^\times$.
  A matrix $s \in \mathrm{GL}_n(F)$ has $\rho_b(s)=1$ if and only if
  $b^{-1}s$ is a transvection or a one-direction dilation.
  Its relative ratio $\theta$ is $1$ in the transvection case and is
  the dilation multiplier otherwise. Each $C(b,\theta)$ is a single
  $\mathrm{GL}_n(F)$ conjugacy class.
\end{lemma}

\begin{proof}
By \Cref{lem-deviation-letters}, $\rho_b(s)=1$ if and only if
  $\rho_1(b^{-1}s)=1$, that is, $b^{-1}s=I+u \varphi$ for nonzero $u$ and $\varphi$.
  By \Cref{lem-elementary-transformations}, $\theta=\det(b^{-1}s)=1+\varphi(u)$, so
  $\theta=1$ gives a transvection and $\theta\neq 1$ gives a one-direction dilation.

  For $\theta\neq 1$, the hyperplane $\ker \varphi$ and line $F u$ are
  complementary eigenspaces of $s$, giving the form
  $\operatorname{diag}(b \theta,b,\dots,b)$ in a suitable basis.
  For $\theta=1$, the deviation $s-b I$ has rank one and square zero,
  so $s$ has one Jordan block of size two at $b$ and $n-2$ blocks of
  size one. These forms show that each $C(b,\theta)$ is one conjugacy class.
\end{proof}

Geometrically, a matrix $s \in C(b,\theta)$ acts as multiplication by $b$ on a
hyperplane and as multiplication by $b \theta$ on its deviation line, so the
relative ratio $\theta$ compares those two multipliers. In particular, $C(1,1)$
is the class of transvections and $C(1,\theta)$ ($\theta\neq 1$) is the class of
one-direction dilations with multiplier $\theta$. Because each $C(b,\theta)$ is a single conjugacy class, we only need to match
determinants: once we construct a product of origin-$1$ factors $T_i$ with
determinants $\theta_i$, rescaling each by $b_i$ automatically places $b_i T_i$
in the prescribed class $C(b_i,\theta_i)$.

Now fix a list of rank-one factors $A_1,\dots,A_k$ with $A_i \in C(b_i,\theta_i)$,
and ask which targets $X \in \mathrm{GL}_n(F)$ can be written as a product of their
conjugates. Because scalars are central, pulling $b_i$ out of each factor
replaces $X$ by $Y=a^{-1}X$, where $a=b_1 \dots b_k$, and reduces the problem to
factoring $Y$ into origin-$1$ rank-one transformations $T_1,\dots,T_k$ of
determinants $\theta_1,\dots,\theta_k$.

From \Cref{sec-deviation-geometry} we already know what can block such a
factorisation. First, determinants must match, so
$\det X = \prod_{i=1}^k \det A_i$ (equivalently
$\det Y = \prod_{i=1}^k \theta_i$). Second, each origin-$1$ factor contributes
deviation rank $1$, so \Cref{cor-product-deviation} forces $k \ge  \rho_a(X)$.
Third, at minimal length $k = \rho_a(X) > 0$, the product is tight, so if $X$
acts as a scalar $\beta$ on its moving space $W=\im(X-a I)$, then
\Cref{lem-tight-scalar-action} forces $\theta_i = \beta/a$ for every $i$. Finally,
if $X = a I$ (so $\rho_a(X) = 0$), then $k$ cannot be $1$, because a single
rank-one factor is never scalar.

Once $n \ge  3$, these obstructions are the only ones. In particular, the
scalar-action restriction disappears as soon as $k > \rho_a(X)$: whereas
tightness forces the factor displacement lines to be independent, even a single
spare factor lets them overlap so that the factors can trade ratios freely while
preserving $\det Y$.

The proof of the following theorem builds the required factorisations in two
stages. First, for minimal length $k = \rho_a(X)$ (and for length
$\rho_a(X)+1$ in the scalar-action case), existence theorems of Ellers and
Lausch supply the base factorisations over any field. Second, to reach all
longer lengths, we use a \emph{factor-splitting} move: because $n \ge  3$, the affine
space of functionals taking a prescribed value on a nonzero vector has
dimension $n-1 \ge  2$, hence at least four points even over $\mathbb{F}_2$. That
room lets us split any rank-one transformation along the same displacement line
into two rank-one factors with any prescribed ratios of the same product.
Iterating this split extends the base products to every length $k > \rho_a(X)$,
and also builds $a I$ whenever $k \ge  2$ by pairing one factor with a split copy
of its inverse.

\begin{theorem}[{Rank-one product criterion}]
\label{thm-prescribed-rank-one}
Let $F$ be a field and $n\ge 3$. For $i=1,\dots,k$ ($k\ge 1$), let
  $A_i \in C(b_i,\theta_i)$ be a rank-one factor at origin $b_i \in F^\times$ with
  relative ratio $\theta_i \in F^\times$. Given $X \in \mathrm{GL}_n(F)$, put
  $a=b_1 \dots b_k$ and $r=\rho_a(X)$.

  Then $X$ is a product of conjugates of $A_1,\dots,A_k$ in $\mathrm{GL}_n(F)$ if and
  only if
  \[
  \det X=\prod_{i=1}^k \det A_i, \quad k\ge r,
\]
  subject to two boundary conditions:
\begin{enumerate}[label=(\roman*)]
  \item if $k=r>0$ and $X$ acts on its moving space $W=\im(X-a I)$ as
      multiplication by a scalar $\beta \in F$, then $\theta_i=\beta/a$ for every
      $i$;
  \item if $X=a I$ (so $r=0$), then $k\ge 2$.
\end{enumerate}
\end{theorem}

\begin{proof}
First divide each factor by its origin, so that all deviations are
  measured from $I$. Put $Y=a^{-1}X$. By \Cref{lem-rank-one-classes}, the
  required factorisation is equivalent to writing $Y=T_1 \dots T_k$, where each
  $T_i$ has rank-one deviation at $1$ and determinant $\theta_i$. We recover the
  original factors by setting $S_i=b_i T_i$.

  \Cref{cor-product-deviation} and multiplicativity of the determinant give
  $k\ge r$ and $\det Y=\prod_i \theta_i$. Suppose $k=r>0$ and the moving
  action of $Y$ is $\nu I_W$, where $\nu=\beta/a$. By
  \Cref{lem-tight-scalar-action}, each $N_i=T_i-I$ satisfies
  $N_i^2=(\nu-1)N_i$. Thus $T_i$ acts by $\nu$ on its deviation line,
  so its relative ratio is $\theta_i=\nu$. A single rank-one factor cannot
  equal $I$, which also gives the necessary restriction when $r=0$.

  We now prove that these necessary conditions suffice.

  \textbf{Tight products.} Suppose $k=r>0$. If $r=1$, take $T_1=Y$.
  For $r\ge 2$, we apply Ellers and Lausch's results over the field $F$.
  Their simple mappings are our origin-one rank-one transformations,
  and their types are the determinants $\theta_i$. The map $Y-I$ induces
  an isomorphism from $F^n / \ker(Y-I)$ to $W$ and intertwines the
  induced action of $Y$ with $Y|_W$. Their condition on the action modulo
  the fixed space is therefore exactly our condition on the moving action.
  When this action is nonscalar, their Theorem 5 gives $r$ factors with
  the first $r-1$ ratios chosen as required; the determinant fixes the
  last ratio. When it is $\nu I_W$ and all ratios equal $\nu$, their
  Proposition 6 gives the required $r$ factors
  \cite[Theorem 5 and Proposition 6]{EllersLausch1989}.

  \textbf{Splitting a factor.} To extend these products to longer lists, we
  replace one rank-one factor by two with any nonzero ratios whose product
  is its ratio. Both new deviations will have the same image as the old
  one. Write
  $T=I+w \varphi$ with $w,\varphi\neq 0$ and $\det T=p q\neq 0$.
  Choose a functional $\psi$ satisfying $\psi(w)=p-1$, with
  $\psi\neq 0,\varphi$. The possible functionals form an affine space of dimension
  $n-1\ge 2$, with at least four elements even over $\mathbb{F}_2$, so these two
  exclusions leave a choice. Set $\chi=(\varphi-\psi)/p$. Then
  \[
  (I+w \psi)(I+w \chi)=I+w(\psi+p \chi)=T,
\]
  and the two factors have determinants $p$ and $q$. Both deviations
  are nonzero and have rank one. Iteration splits $T$ into any ordered
  list of nonzero ratios with product $\det T$.

  \textbf{Longer products.} Now suppose $k>r>0$. If $r=1$, split $Y$ itself.
  If the moving action is nonscalar, use the tight case with ratios
  \[
  \theta_1,\dots,\theta_{r-1}, \prod_{i=r}^k \theta_i,
\]
  then split the last factor into the remaining ratios.

  For scalar moving action $\nu I_W$ with $r\ge 2$, a tight product would
  require every ratio to equal $\nu$. Instead we start with $r+1$ factors,
  using the scalar clause of Ellers and Lausch's Theorem 5. This lets us
  choose the first $r$ ratios freely if $\nu=1$ or one of those ratios
  differs from $\nu$; the determinant fixes the last. Group the required
  ratios into the $r+1$ entries
  \[
  \theta_1,\dots,\theta_r,\gamma, \quad \gamma=\prod_{i=r+1}^k \theta_i.
\]
  Because $(Y-I)(F^n) \subseteq W$, the map $Y$ induces the identity on
  $F^n / W$, so their product is $\det Y = \nu^r$.
  If $\nu\neq 1$, the grouped ratios cannot all equal $\nu$, as that would
  give product $\nu^{r+1}\neq \nu^r$. Place a differing ratio among the first
  $r$ if necessary and apply the theorem. The identity
  \[
  A B=B(B^{-1}A B)
\]
  exchanges adjacent factor classes without changing their product, so we can
  restore the grouped order and, if $k>r+1$, split the factor of ratio $\gamma$.
  This gives the original ordered list.

  \textbf{Zero deviation.} If $Y=I$ and $k\ge 2$, the determinant condition is
  $\prod_i \theta_i=1$. Choose a rank-one factor $T$ of ratio
  $\theta_1$. Its inverse has rank-one deviation and ratio
  $\theta_1^{-1}=\prod_{i=2}^k \theta_i$. Split $T^{-1}$ into these
  remaining ratios and prefix $T$. For $k=2$, use $T^{-1}$ itself.
  This gives product $I$. Restoring the origins by $S_i=b_i T_i$
  gives $X$ in the required conjugacy classes in every case.
\end{proof}

\begin{remark}[{Relation to earlier length theorems}]
\label{rem-rank-one-attribution}
Over fields, prescribing the types of simple-mapping factors goes back to
  Ellers \cite{Ellers1977} and Ellers and Lausch \cite[Theorem 5,
  Proposition 6]{EllersLausch1989}, whose scalar-action clause also supplies
  determinant-compatible lists one above the deviation rank. After dividing $X$
  by $a = \prod b_i$ so that $Y = a^{-1} X$, the path of Kn\"uppel and Nielsen
  \cite[Observation 7.2, Propositions 7.4--7.5, Corollary 7.6]{KnuppelNielsen2006} is
  precisely $\operatorname{im}(Y-I)$, their simple mappings are the normalised rank-one
  factors, and their path homotheties are the matrices with scalar moving
  action; their Proposition 7.5 already states strict-rank existence for
  prescribed simple classes, apart from the omitted one-letter identity
  exception $k=1$, $Y=I$ (also present in \cite[Proposition 2.3]{KnuppelNielsen2013}).
  The proof above gives a self-contained derivation for all lengths via
  splitting and grouping, pairs it with the geometric necessity conditions of
  \Cref{sec-deviation-geometry}, and recovers independent origins by rescaling.

  Ellers and Lausch \cite[Propositions 7--8]{EllersLausch1989} also give the
  transvection minimum and treat reflections (see also
  \cite[\S~7]{KnuppelNielsen2006}); their commutative-local-ring predecessor is Ellers
  and Ishibashi \cite{EllersIshibashi1987}, and Phadke \cite{Phadke1974} discusses
  Dieudonn\'e's transvection bound and the additional complications over
  noncommutative division rings.
\end{remark}

\subsection{Transvection lengths and lower bounds}
\label{sec-rank-one-lower-bounds}
We now specialise \Cref{thm-prescribed-rank-one} to the transvection class $C(1,1)$
to determine exact transvection lengths in $\mathrm{SL}_n(F)$, recovering the length
formula of Ellers and Lausch \cite[Proposition 7]{EllersLausch1989}.

\begin{corollary}[{Transvection lengths}]
\label{cor-transvection-lengths}
Let $F$ be a field, $n\ge 3$, and let $T$ be the $\mathrm{SL}_n(F)$ conjugacy
  class of a transvection. For $X \in \mathrm{SL}_n(F)$ with $X\neq I$, put
  $W=\im(X-I)$ and $r=\rho_1(X)$. Then
  \[
  \lVert X \rVert_T=\begin{cases}r+1 & \textnormal{if } X|_W=\beta I_W \textnormal{ for some } \beta\neq 1 \\ r & \textnormal{otherwise.}\end{cases}
\]
  The identity has length zero. These lengths are unchanged if only
  positive factors from $T$ are allowed.
\end{corollary}

\begin{proof}[{Proof of \Cref{cor-transvection-lengths}}]
First, all transvections are conjugate in $\mathrm{SL}_n(F)$.
  By \Cref{lem-rank-one-classes} they form the single $\mathrm{GL}_n(F)$ class
  $C(1,1)$. Fix a transvection $t$. In a Jordan basis for $t$, a matrix
  that is the identity on its nontrivial two-dimensional block and scales
  one complementary basis vector commutes with $t$. Such a vector exists
  because $n\ge 3$, and its scaling can be chosen to give any nonzero
  determinant. Given a $\mathrm{GL}_n(F)$ conjugator $g$, choose such a matrix
  $h$ with $\det h=(\det g)^{-1}$. Then $g h$ has determinant one and
  gives the same conjugate of $t$. Inverses are themselves transvections,
  so $T^{-1}=T$.

  Apply \Cref{thm-prescribed-rank-one} with each $A_i$ a transvection at
  origin $1$.
  For $X \in \mathrm{SL}_n(F)$ the determinant condition holds, and the
  rank condition excludes lengths below $r$. Length $r$ is allowed
  unless $X|_W$ is scalar with multiplier different from $1$.
  In that case length $r+1$ is allowed. This proves the length formula
  for positive and signed words alike.
\end{proof}

\begin{figure}[htbp]
  \centering
  \includegraphics[width=\linewidth]{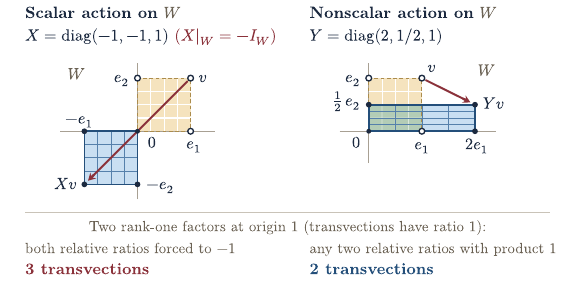}
  \caption{Two targets in $\mathrm{SL}_3(\mathbb{Q})$ fixing $e_3$ with the same moving
    plane $W = \operatorname{span}\{e_1, e_2\}$ and deviation rank $2$. Because $X$ acts
    as a scalar ($-I_W$) on $W$, \Cref{cor-transvection-lengths} forces three
    transvections, whereas $Y$ needs only two.}\label{fig-moving-action-comparison}
\end{figure}

The scalar-action penalty is not confined to scalar matrices: a nonscalar target
$X \in \mathrm{SL}_n(F)$ can still act as a nonidentity scalar on its moving space $W$,
and therefore need $r+1$ transvections even when its deviation rank is $r < n$.
\Cref{fig-moving-action-comparison} illustrates this in $\mathrm{SL}_3(\mathbb{Q})$: the matrices
$X=\operatorname{diag}(-1,-1,1)$ and $Y=\operatorname{diag}(2,1 / 2,1)$ share the same
determinant, the same deviation rank $2$ and the same moving plane
$W=\operatorname{span}\{e_1,e_2\}$, yet $X$ reverses every vector in $W$ and forces three
transvections, whereas $Y$ acts nonscalarly on $W$ and needs only two.

A nonscalar target can never require $n+1$ transvections: if a matrix has full
deviation rank $r=n$, its moving space is all of $F^n$, so scalar action on $W$
forces the matrix itself to be a scalar. In fact, the rank-one product
criterion (\Cref{thm-prescribed-rank-one}) shows something stronger: not only does
every nonscalar have transvection length at most $n$, but every nonscalar (together
with the identity) can be written as a product of \emph{exactly} $n$ positive
transvections, leaving only the nonidentity scalars outside $T^n$.

\begin{corollary}[{The $n$th power of the transvection class}]
\label{cor-transvection-nth-power}
Let $F$ be a field, $n\ge 3$, and let $T$ be the transvection class in
  $\mathrm{SL}_n(F)$. Then
  \[
  T^n = \mathrm{SL}_n(F) \setminus \{\zeta I : \zeta \in \mu_n(F), \zeta \neq  1\}.
\]
\end{corollary}

\begin{proof}
Apply the rank-one product criterion (\Cref{thm-prescribed-rank-one}) to $n$
  transvections, all with origin $1$ and relative ratio $1$. For a nonscalar
  $X \in \mathrm{SL}_n(F)$, either $\rho_1(X) < n$, or $\rho_1(X) = n$ and the moving
  action on $W = F^n$ is $X$ itself, hence nonscalar; in either case
  \Cref{thm-prescribed-rank-one} gives $X \in T^n$. It also gives $I \in T^n$ by
  the zero-deviation condition, since $n \ge  2$. By contrast, a nonidentity
  scalar $\zeta I$ has $\rho_1(\zeta I) = n$ and scalar moving action $\zeta I$,
  which would require every factor to have relative ratio $\zeta \neq  1$. This
  excludes precisely the scalars shown.
\end{proof}

\Cref{cor-transvection-nth-power} pins down the largest transvection length in
$\mathrm{SL}_n(F)$. Every nonscalar lies in $T^n$ and therefore has length at most
$n$, with length $n$ always attained over every field by
\Cref{ex-full-deviation-target}. Length $n+1$ occurs precisely at the nonidentity
scalars $\zeta I$ in $\mathrm{SL}_n(F)$, that is, for nontrivial $n$th roots of unity
$\zeta \in \mu_n(F) \setminus \{1\}$. Because a single transvection normally
generates $\mathrm{SL}_n(F)$, this maximum gives a lower bound for
$\Delta(\mathrm{SL}_n(F))$; and because the transvection class has determinant $1$ (so
its target fibre at every length is $\mathrm{SL}_n(F)$), the same targets yield the
matching lower bound for $\operatorname{cn}_\mathrm{ab}(\mathrm{GL}_n(F))$.

\begin{corollary}[{Lower bounds for diameter and fibre covering}]
\label{cor-sl-diameter-lower-bound}
Let $F$ be a field, $n\ge 3$, and let $T$ be the transvection class in
  $\mathrm{SL}_n(F)$. Then
  \[
  \operatorname{diam}_T(\mathrm{SL}_n(F))=n+[\mu_n(F)\neq \{1\}],
\]
  where $\mu_n(F)=\{\zeta \in F^\times : \zeta^n=1\}$ is the set of $n$th roots
  of unity in $F$. Consequently,
  \begin{align*}
  \Delta(\mathrm{SL}_n(F)) &\ge n+[\mu_n(F)\neq \{1\}], \\
    \operatorname{cn}_\mathrm{ab}(\mathrm{GL}_n(F)) &\ge n+[\mu_n(F)\neq \{1\}].
\end{align*}
\end{corollary}

\begin{proof}[{Proof of \Cref{cor-sl-diameter-lower-bound}}]
By \Cref{cor-transvection-nth-power}, every nonscalar $X \in \mathrm{SL}_n(F)$ lies in
  $T^n$ and hence has transvection length at most $n$. The companion matrix in
  \Cref{ex-full-deviation-target} is nonscalar and has no eigenvalue $1$, so its
  moving space is all of $F^n$ and \Cref{cor-transvection-lengths} gives length
  exactly $n$.

  A nonidentity scalar $\zeta I$ acts by $\zeta$ on its moving space $F^n$,
  giving length $n+1$ by \Cref{cor-transvection-lengths}. Such a scalar belongs to
  $\mathrm{SL}_n(F)$ exactly when $\zeta^n=1$. The diameter is therefore $n+1$
  if $\mu_n(F)\neq \{1\}$, and $n$ otherwise.

  A single transvection $t \in T$ normally generates $\mathrm{SL}_n(F)$
  and gives the same word lengths as $T$. Its diameter is therefore a
  lower bound for $\Delta(\mathrm{SL}_n(F))$.

  By the proof of \Cref{cor-transvection-lengths}, $T$ is also the full
  $\mathrm{GL}_n(F)$ transvection class. Its determinant is one, so its target fibre
  at every length is $\mathrm{SL}_n(F)$. If $T^k=\mathrm{SL}_n(F)$, every target has signed
  word length at most $k$. The filling length of $T$ is therefore at least its
  signed diameter. Taking the supremum over the eligible full GL classes proves
  the lower bound for $\operatorname{cn}_\mathrm{ab}(\mathrm{GL}_n(F))$.
\end{proof}

Thus the transvection class alone establishes the lower bounds for both
$\Delta(\mathrm{SL}_n(F))$ and $\operatorname{cn}_\mathrm{ab}(\mathrm{GL}_n(F))$; \Cref{sec-class-products} turns
to the matching upper bounds for arbitrary noncentral classes.

\section{Class products and upper bounds}
\label{sec-class-products}
In \Cref{sec-rank-one}, the transvection class gave the common lower bound
$n + [\mu_n(F) \neq  \{1\}]$ for both $\Delta(\mathrm{SL}_n(F))$ and
$\operatorname{cn}_\mathrm{ab}(\mathrm{GL}_n(F))$ (\Cref{cor-sl-diameter-lower-bound}). We now prove
that these bounds are sharp: no noncentral conjugacy class requires longer
words than a transvection. After reducing both upper bounds to covering
nonscalar targets in $n$ letters (\Cref{sec-class-product-obligations}), we
establish the special linear diameter in \Cref{sec-repeated-class-coverage}
and determinant-fibre covering in
\Cref{sec-determinant-scalar-construction}.

\subsection{Reduction to nonscalar targets}
\label{sec-class-product-obligations}
We saw in \Cref{cor-transvection-nth-power} that $n$ positive transvections
already cover every nonscalar in $\mathrm{SL}_n(F)$ and miss only the nonidentity
scalars. Each transvection has the smallest possible positive deviation rank,
whereas a general noncentral conjugacy class $C$ can displace vectors in
several independent directions per factor, so we may expect $C^n$ to cover all
nonscalars of its determinant fibre as well. Both upper bounds therefore aim
for the nonscalars in $n$ letters first, using one extra letter when length
$n+1$ is allowed. For signed words in $\mathrm{SL}_n(F)$, that extra step is
immediate: any scalar $\zeta I$ can be written as $(\zeta c^{-1}) c$ for any
$c \in C$, because $\zeta c^{-1}$ is nonscalar.

For exact positive fibre covering in $\mathrm{GL}_n(F)$, the extra letter has to do
more: multiplying by a class $C$ of determinant $d$ shifts the determinant
fibre from $\mathcal{D}(a)$ to $\mathcal{D}(a d)$ (in the notation of
\Cref{sec-matrix-normal-generation}), and we must reach every target
$X \in \mathcal{D}(a d)$, nonscalar as well as scalar, at the exact same length.
Even when $X$ is nonscalar, $X c^{-1}$ is nonscalar provided $c$ is not a
scalar multiple of $X$; as the following lemma shows, every noncentral
conjugacy class contains two nonproportional matrices, so such a choice of $c$
is always available.

\begin{lemma}[{One-letter completion of a determinant fibre}]
\label{lem-one-letter-completion}
Let $F$ be a field, $n\ge 2$, and $a,d \in F^\times$. Suppose
  $P \subseteq \mathcal{D}(a)$ contains every nonscalar in $\mathcal{D}(a)$.
  Let $C$ be the conjugacy class of a nonscalar matrix of determinant
  $d$ under either $\mathrm{SL}_n(F)$ or $\mathrm{GL}_n(F)$. Then
  \[
  P C = \mathcal{D}(a d).
\]
\end{lemma}

\begin{proof}
Fix $X \in \mathcal{D}(a d)$. A choice of $c \in C$ gives a scalar
  $X c^{-1}$ only if $c$ is proportional to $X$. Thus it suffices to
  find two nonproportional members of $C$.

  Take $c_0 \in C$. Some elementary transvection $h=I+E_{i j}$ fails
  to commute with $c_0$, since a matrix commuting with every off-diagonal
  matrix unit is scalar. Put $c_1=h c_0 h^{-1}$, which belongs to $C$ since
  $\det h=1$. If $c_1=\lambda c_0$ for some $\lambda \in F^\times$, then
  $c_0^{-1} h c_0=\lambda h$; comparing the sole eigenvalue $1$ on the left with
  $\lambda$ on the right gives $\lambda=1$, contrary to the choice of $h$.

  Choose $c \in \{c_0,c_1\}$ so that $X c^{-1}$ is nonscalar. It has determinant
  $a$, hence lies in $P$, so $X=(X c^{-1})c \in P C$. Every element of $P C$ has
  determinant $a d$, giving the reverse inclusion.
\end{proof}

Once $C^n$ covers the nonscalars of $\mathcal{D}(d^n)$, \Cref{lem-one-letter-completion}
settles both upper bounds immediately whenever $\mu_n(F) \neq  \{1\}$, giving
$C^{n+1} = \mathcal{D}(d^{n+1})$. When $\mu_n(F) = \{1\}$, the only scalar in
$\mathcal{D}(d^n)$ is $d I$: for $\Delta(\mathrm{SL}_n(F))$ ($d=1$) the identity already
has length zero, whereas exact positive fibre covering still requires
$d I \in C^n$, which we prove in \Cref{sec-determinant-scalar-construction}.
We now apply this reduction to determine $\Delta(\mathrm{SL}_n(F))$ in
\Cref{sec-repeated-class-coverage}.

\subsection{Special linear diameter}
\label{sec-repeated-class-coverage}
To determine $\Delta(\mathrm{SL}_n(F))$, we first establish $n$-letter coverage for
every noncentral conjugacy class in $\mathrm{SL}_n(F)$. How we obtain that coverage
depends on the size of the field: published covering theorems handle all fields
with $|F| \ge  4$ as well as $\mathbb{F}_3$ outside dimension four, leaving only
$\mathbb{F}_2$ and $\mathrm{SL}_4(\mathbb{F}_3)$ for direct construction in
\Cref{sec-small-field-proofs}.

Over any field with at least four elements, the following theorem of Lev shows
that $n$ letters from any noncentral $\mathrm{SL}_n(F)$-conjugacy class cover every
nonscalar of the matching determinant. We state it in its full generality for a
seed of arbitrary determinant $d \in F^\times$: setting $d=1$ gives the special
linear coverage needed here, while keeping general $d$ lets us apply the same
statement directly to determinant fibres of $\mathrm{GL}_n(F)$ in
\Cref{sec-determinant-scalar-construction}.

\begin{theorem}[{Nonscalar coverage \cite[Theorem 2]{Lev1996}}]
\label{prop-sl-nonscalar-coverage}
Let $n\ge 3$, let $F$ be a field with $\lvert F \rvert\ge 4$, and let $C$ be the
  $\mathrm{SL}_n(F)$-conjugacy class of a nonscalar element of $\mathrm{GL}_n(F)$
  of determinant $d$. Then $C^n$ contains every nonscalar
  $X \in \mathrm{GL}_n(F)$ with $\det X=d^n$.
\end{theorem}

Over $\mathbb{F}_3$, we obtain the same $n$-letter nonscalar coverage in every
dimension $n\ge 3$ except $n=4$ from Kn\"uppel and Nielsen's extended covering
theorem, which lets us prescribe $n+1$ noncentral classes whose product
contains the identity. Reserving one slot for the class of $X^{-1}$ leaves $n$
slots for copies of our seed class $C$.

\begin{proposition}[{Nonscalar targets over $\mathbb{F}_3$}]
\label{prop-sl-f3-nonscalar-coverage}
Let $n\ge 3$ with $n\neq 4$, and let $C$ be a noncentral conjugacy class
  in $\mathrm{SL}_n(\mathbb{F}_3)$. Then $C^n$ contains every nonscalar in
  $\mathrm{SL}_n(\mathbb{F}_3)$.
\end{proposition}

\begin{proof}
Fix a nonscalar target $X \in \mathrm{SL}_n(\mathbb{F}_3)$. In these dimensions
  the identity belongs to the product of any $n+1$ prescribed noncentral
  classes \cite[p. 635]{KnuppelNielsen2006}. Prescribe $C$ in the first $n$ places and the class of
  $X^{-1}$ in the last place. The product of the first $n$ letters is
  then conjugate to $X$ in $\mathrm{SL}_n(\mathbb{F}_3)$. Since $C^n$ is closed
  under conjugation by this group, it contains $X$.
\end{proof}

\begin{remark}[{The part of the extended covering theorem used}]
\label{rem-kn-repeated-class}
Our repeated-class specialisation avoids a gap in the printed proof of
  Kn\"uppel and Nielsen's general theorem: the final case of their
  Proposition 8.2 omits one class from the product it constructs
  \cite[p. 654, step (vi)]{KnuppelNielsen2006}. Because that case requires at least
  three and at most $n-1$ classes of rank-one deviation, it cannot occur when
  $n$ of the $n+1$ prescribed classes coincide.
\end{remark}

In the two remaining cases, the transvection lower bound
(\Cref{cor-sl-diameter-lower-bound}) tells us how many letters we may use: over
$\mathbb{F}_2$ there is no nonidentity scalar, so we must stay within $n$ letters,
whereas in $\mathrm{SL}_4(\mathbb{F}_3)$ the central element $-I$ already requires five
transvections, allowing five letters. Unlike \Cref{prop-sl-nonscalar-coverage} and
\Cref{prop-sl-f3-nonscalar-coverage}, our small-field constructions cover the
entire group (including the scalar matrices) directly.

\begin{proposition}[{Positive coverage over the small fields}]
\label{prop-sl-small-field-coverage}
\begin{enumerate}[label=(\roman*)]
  \item For $n\ge 3$, every noncentral conjugacy class $C$ in
     $\mathrm{SL}_n(\mathbb{F}_2)$ satisfies $C^n=\mathrm{SL}_n(\mathbb{F}_2)$.
  \item Every noncentral conjugacy class $C$ in $\mathrm{SL}_4(\mathbb{F}_3)$
     satisfies $C^5=\mathrm{SL}_4(\mathbb{F}_3)$.
\end{enumerate}
\end{proposition}

\begin{proof}
The constructions are given in \Cref{sec-fibre-f2} (Positive products over
  $\mathbb{F}_2$) and \Cref{sec-sl-f2-comparison} (Positive products in
  $\mathrm{SL}_4(\mathbb{F}_3)$), respectively.
\end{proof}

Most cases of \Cref{prop-sl-small-field-coverage} are proved in
\Cref{sec-small-field-proofs} by multiplying conjugates of $C$ into a cyclic
intermediate (\Cref{sec-small-field-cyclic-constructions}) and applying B\"unger and
Nielsen's cyclic decomposition theorem \cite[Theorem 1.2]{BungerNielsen1999}. The
one exception where that cyclic route takes a letter too many is a rank-two
seed over $\mathbb{F}_2$ in even dimension $n\ge 6$: there
\Cref{sec-f2-rank-two-comparison} instead starts from an $n$-transvection word
supplied by \Cref{cor-transvection-nth-power} and replaces each consecutive pair
$t_1 t_2 = (t_1 p)(p t_2)$ (using a suitable transvection $p$ with $p^2 = I$)
by two letters of $C$ (\Cref{lem-f2-pairs}), illustrating how the rank-one
criterion also serves as an upper-bound tool.

Combining \Cref{prop-sl-nonscalar-coverage}, \Cref{prop-sl-f3-nonscalar-coverage}
and \Cref{prop-sl-small-field-coverage} with the one-letter completion lemma
(\Cref{lem-one-letter-completion}) shows that no noncentral seed in $\mathrm{SL}_n(F)$
requires more letters than a transvection. Over every field, the hardest
nonscalar targets take $n$ transvections and a nonidentity scalar takes one
more, so the field enters the diameter only through whether it contains a
nontrivial $n$th root of unity.

\begin{theorem}[{Special linear diameter}]
\label{thm-sl-diameter}
Let $F$ be a field and $n\ge 3$. Then
  \[
  \Delta(\mathrm{SL}_n(F)) = n + [\mu_n(F) \neq  \{1\}],
\]
  where $\mu_n(F) = \{\zeta \in F^\times : \zeta^n=1\}$. The supremum is
  attained by taking a single transvection as seed.
\end{theorem}

\begin{proof}
The normal generation proposition for the special linear group
  (\Cref{prop-sl-normal-generation}) reduces the diameter to the supremum over
  noncentral single seeds, and the transvection lower bound
  (\Cref{cor-sl-diameter-lower-bound}) gives the stated value for one such seed.

  For the upper bound, fix any noncentral seed and let $C$ be its conjugacy
  class in $\mathrm{SL}_n(F)$. There are four cases.

\begin{itemize}
  \item \textbf{$\lvert F \rvert\ge 4$.} Lev's nonscalar coverage theorem
(\Cref{prop-sl-nonscalar-coverage}, with $d=1$) puts every nonscalar in $C^n$.
  \item \textbf{$F=\mathbb{F}_3$, $n\neq 4$.} \Cref{prop-sl-f3-nonscalar-coverage} does the same.
  \item \textbf{$F=\mathbb{F}_2$.} Small-field positive coverage
(\Cref{prop-sl-small-field-coverage}(i)) gives $C^n=\mathrm{SL}_n(\mathbb{F}_2)$, and
$1$ is the only root of unity.
  \item \textbf{$F=\mathbb{F}_3$, $n=4$.} Small-field positive coverage
(\Cref{prop-sl-small-field-coverage}(ii)) gives $C^5=\mathrm{SL}_4(\mathbb{F}_3)$, as
required since $-I$ belongs to the centre.
\end{itemize}

  In the first two cases every nonscalar has signed length at most $n$ and the
  identity has length zero; any nonidentity scalar $\zeta I$
  ($\zeta \in \mu_n(F) \setminus \{1\}$) lies in $C^{n+1}$ by the one-letter
  completion lemma (\Cref{lem-one-letter-completion}). In the last two cases the
  whole group is covered at once. Positive words are also signed words, so the
  upper bound holds for every seed in the supremum.
\end{proof}

\subsection{Determinant-fibre covering}
\label{sec-determinant-scalar-construction}
For the special linear diameter we could count the identity as reached at
length zero, without writing down a single letter. Exact positive
determinant-fibre covering in $\mathrm{GL}_n(F)$ removes that freedom: as noted in
\Cref{sec-class-product-obligations}, once $C^n$ covers the nonscalars of
$\mathcal{D}(d^n)$, the case $\mu_n(F) \neq  \{1\}$ is finished immediately by the
one-letter completion lemma (\Cref{lem-one-letter-completion}), whereas the case
$\mu_n(F) = \{1\}$ requires writing the unique scalar $d I \in \mathcal{D}(d^n)$ as a
product of exactly $n$ conjugates from $C$. We treat broad fields $|F| \ge  4$
first and then collect the small-field results over $\mathbb{F}_2$ and $\mathbb{F}_3$.

Over fields with $|F| \ge  4$, Lev's nonscalar coverage theorem
(\Cref{prop-sl-nonscalar-coverage}) already puts every nonscalar of determinant
$d^n$ in $C_0^n$ for any noncentral $\mathrm{SL}_n(F)$-conjugacy class $C_0$ of
determinant $d$, and hence also in $C^n$ for the full $\mathrm{GL}_n(F)$-conjugacy
class $C$ containing $C_0$.\footnote{Kn\"uppel and Nielsen give the nonscalar
bound for full general linear classes in \cite[Theorem 1.1]{KnuppelNielsen2013};
their inequality immediately preceding that theorem gives full coverage of the
next determinant fibre.} In fact, over every such field, $d I$ lies in $C^n$
regardless of $\mu_n(F)$ (and in $C_0^n$ when $\mu_n(F) = \{1\}$). Two regimes
are already covered: the rank-one product criterion
(\Cref{thm-prescribed-rank-one}) handles all rank-one classes over every field,
and Lev's scalar coverage theorem handles non-rank-one classes for $n \ge  5$.
This leaves non-rank-one classes in dimensions $3$ and $4$.

In dimensions $3$ and $4$, we construct $d I \in C^n$ by taking a
representative $c \in C$ and multiplying $n$ conjugates of $c$ formed by
successive powers of a single matrix $p$: factoring $c = s p^{-1}$ makes the
interior powers of $p$ telescope. Because \Cref{sec-small-field-proofs} reuses this
identity at lengths other than $n$, we state it for an arbitrary group and
length $m$.

\begin{lemma}[{Rotation identity}]
\label{lem-rotation-identity}
Let $c,p$ be elements of a group and let $m\ge 1$ be an integer. In
  increasing order of $i$, the $m$ conjugates have product
  \[
  \prod_{i=0}^{m-1} p^i c p^{-i} = (c p)^m p^{-m}.
\]
  In particular, if $c=s p^{-1}$ and $p^m=1$, then their product is $s^m$.
\end{lemma}

\begin{proof}
Writing $p^i c p^{-i} = p^i (c p) p^{-(i+1)}$ exposes the cancellation
  between consecutive factors:
  \[
  \prod_{i=0}^{m-1} p^i c p^{-i}
    = (c p) p^{-1} \cdot p (c p) p^{-2} \, \cdots \, p^{m-1} (c p) p^{-m}
    = (c p)^m p^{-m}.
\]
  Substituting $c p = s$ and $p^m = 1$ gives $s^m$.
\end{proof}

To apply the rotation identity (\Cref{lem-rotation-identity}) with $m=n$, we seek
a factorisation $c = s p^{-1}$ with $s^n = d I$ and $p^n = I$. Both powers are
realised naturally by cyclic shifts of a basis: cycling
$e_1 \mapsto e_2 \mapsto \cdots \mapsto e_n \mapsto d e_1$ multiplies every vector
by $d$ after $n$ steps, while an unscaled cyclic shift returns to the identity.
Choosing their bases independently via cyclic-class decomposition theorems
allows the product $s p^{-1}$ to represent any nonscalar $c$ in dimensions $3$
and $4$.

\begin{proposition}[{Reaching the determinant scalar}]
\label{prop-gl-determinant-scalar}
Let $n\ge 3$ and let $F$ be a field with at least four elements. If $C$ is a
  noncentral conjugacy class in $\mathrm{GL}_n(F)$ of determinant $d$, then
  $d I \in C^n$. If moreover $\mu_n(F)=\{1\}$, then $d I \in C_0^n$ for
  every $\mathrm{SL}_n(F)$-conjugacy class $C_0$ contained in $C$.
\end{proposition}

\begin{proof}
Fix $c \in C$, and let $C_0 \subseteq C$ be the $\mathrm{SL}_n(F)$-conjugacy class
  of $c$. Suppose first that $c$ has rank-one deviation at an origin $b$. For
  $n$ copies of $c$ the product origin is $a=b^n$, and every factor has
  relative ratio $\theta=d/a$. The target $d I$ has the required determinant
  $d^n$. If $\theta=1$, then $d I=a I$, so the zero-deviation clause of the
  rank-one product criterion (\Cref{thm-prescribed-rank-one}) applies. If
  $\theta\neq 1$, then $d I-a I$ has rank $n$ and $d I$ acts by $d$ on its moving
  space $F^n$, so the required ratio is $d/a$ in every slot and
  \Cref{thm-prescribed-rank-one} again applies. This proves $d I \in C^n$ over any
  field in the rank-one case. Moreover, $C_0 = C$ here: writing $c=b I+v \varphi$
  with a nonzero vector $v$ and linear functional $\varphi$, choose an
  $(n-2)$-dimensional subspace $U \subseteq \ker \varphi$ not containing $v$, and a
  complement $W$ of $U$ containing $v$. Any invertible map that is arbitrary
  on $U$ and the identity on $W$ fixes $v$ and preserves $\varphi$, hence commutes
  with $c$; since $n\ge 3$, such maps realise every nonzero determinant, so any
  $\mathrm{GL}_n(F)$-conjugator can be adjusted to determinant one.

  Now suppose $c$ has no rank-one deviation at any origin. If $n\ge 5$, Lev's
  scalar coverage theorem \cite[Theorem 4(3)]{Lev1996} puts every scalar of
  determinant $d^n$ in $C_0^n$ (and hence in $C^n$); its excluded classes,
  generalised transvections and generalised dilatations, are precisely the
  rank-one classes just treated.

  It remains to treat $n=3,4$. For a monic polynomial $f \in F[t]$ of degree
  $r\ge 1$, write $\operatorname{Comp}(f)$ for its companion matrix (multiplication by $t$
  on $F[t] / (f)$ in the basis $1,t,\dots,t^{r-1}$, with ones on the
  subdiagonal), which satisfies $f(\operatorname{Comp}(f))=0$. Set
  \[
  R=\operatorname{Comp}(t^n-d), \quad P=\operatorname{Comp}(t^n-1),
\]
  so $R^n=d I$ and $P^n=I$. Both $R$ and $P^{-1}$ are cyclic, and their
  determinants multiply to $d$, since $\det R=(-1)^{n-1}d$ and
  $\det P=(-1)^{n-1}$. If the characteristic of $F$ does not divide $n$, then
  $1$ is a simple eigenvalue of $P$, hence also of $P^{-1}$, and the
  prescribed cyclic class decomposition of B\"unger and Nielsen
  \cite[Theorem 1.2]{BungerNielsen1999}, in the form given by Kn\"uppel and Nielsen
  \cite[Proposition 4.2]{KnuppelNielsen2006} (see the remark following
  \Cref{prop-cyclic-power} in \Cref{sec-small-field-cyclic-constructions} for
  corrections to three technical slips in the printed proof), puts the
  nonscalar $c$ in the product of the full classes of $P^{-1}$ and $R$. If the
  characteristic divides $n$, then in dimensions $3$ and $4$ we have
  $t^n-1=(t-1)^n$, so all eigenvalues of $P^{-1}$ lie in $F$, and Lev's cyclic
  product theorem \cite[Theorem 5]{Lev1996} (originally \cite{Lev1994}) gives the same
  inclusion.

  Full conjugacy classes commute as sets, since $u v=v (v^{-1} u v)$, so we
  can write $c=s p^{-1}$ with $s$ conjugate to $R$ and $p$ conjugate to $P$.
  The rotation identity (\Cref{lem-rotation-identity}) then gives $d I=s^n$ as a
  product of $n$ conjugates of $c$ by powers of $p$. Here
  $\det p=\det P=(-1)^{n-1}$: in dimension three this gives $\det p=1$
  unconditionally, so all three letters lie in $C_0$; in dimension four,
  $\det p=-1$ is a fourth root of unity, so when $\mu_4(F)=\{1\}$ we must have
  $-1=1$ and again $\det p=1$, placing the product in $C_0^4$.
\end{proof}

For cyclic matrices in characteristic zero, $d I \in C^n$ under full general
linear conjugacy also follows by evaluating the matrix-polynomial factorisation
of Laffey and Meehan \cite{LaffeyMeehan1992} at zero; for the special linear clause
of \Cref{prop-gl-determinant-scalar}, however, the determinant-one control on the
conjugator $p$ is essential.

Over the two smallest fields we obtain full-fibre coverage directly in
\Cref{sec-small-field-proofs}. Over $\mathbb{F}_2$ the general and special linear groups
coincide, so \Cref{prop-sl-small-field-coverage}(i) already gives
$C^n = \mathrm{GL}_n(\mathbb{F}_2)$. Over $\mathbb{F}_3$, $-1$ is the only nontrivial root of
unity, so the transvection lower bound (\Cref{cor-sl-diameter-lower-bound})
allows $n$ letters in odd dimension and $n+1$ in even dimension.

\begin{proposition}[{Determinant fibres over $\mathbb{F}_3$}]
\label{prop-gl-f3-fibre}
Let $n\ge 3$, and let $C$ be a noncentral conjugacy class in
  $\mathrm{GL}_n(\mathbb{F}_3)$ of determinant $d$. Put $k=n$ for odd $n$ and
  $k=n+1$ for even $n$. Then $C^k=\mathcal{D}(d^k)$.
\end{proposition}

\begin{proof}
The construction is given in \Cref{sec-fibre-f3} (Determinant fibres over
  $\mathbb{F}_3$).
\end{proof}

Combining the broad-field inputs (\Cref{prop-sl-nonscalar-coverage},
\Cref{prop-gl-determinant-scalar} and \Cref{lem-one-letter-completion}) with the
small-field coverage results (\Cref{prop-sl-small-field-coverage}(i) and
\Cref{prop-gl-f3-fibre}) now determines $\operatorname{cn}_\mathrm{ab}(\mathrm{GL}_n(F))$ over every
field. The two invariants agree because the scalars in $\mathcal{D}(d^n)$ are
precisely $d$ times the centre of $\mathrm{SL}_n(F)$: the identity translates to
$d I$ (reached in $n$ letters when $\lvert F \rvert\ge 4$ by \Cref{prop-gl-determinant-scalar}), whereas any
other scalar corresponds to a nonidentity central element that already requires
$n+1$ transvections.

\begin{theorem}[{Determinant-fibre covering}]
\label{thm-gl-fibre-covering}
Let $F$ be a field and $n\ge 3$. Then
  \[
  \operatorname{cn}_\mathrm{ab}(\mathrm{GL}_n(F)) = \Delta(\mathrm{SL}_n(F))
    = n + [\mu_n(F) \neq  \{1\}].
\]
  At this common length every noncentral conjugacy class in $\mathrm{GL}_n(F)$
  fills its determinant fibre. The transvection class attains the
  supremum of the filling lengths.\footnote{For the abelianised covering
  number see Vaserstein and Wheland \cite[\S~2]{VasersteinWheland1995}.
  In dimension three over finite fields, Orevkov classifies the products
  of prescribed conjugacy classes that contain the identity, including
  those in the general linear group \cite{Orevkov2013}.}
\end{theorem}

\begin{proof}
Write $k=n+[\mu_n(F)\neq \{1\}]$. The transvection lower bound
  (\Cref{cor-sl-diameter-lower-bound}) gives a filling length of at least $k$ for
  the transvection class. For the upper bound, fix any noncentral full
  conjugacy class $C$ of determinant $d$, whose target at length $k$ is
  $D_k(C)=\mathcal{D}(d^k)$. There are four cases.

\begin{itemize}
  \item \textbf{$\lvert F \rvert\ge 4$, $\mu_n(F)=\{1\}$.} By Lev's nonscalar coverage theorem
(\Cref{prop-sl-nonscalar-coverage}), $C^n$ contains every nonscalar in
$\mathcal{D}(d^n)$. The only scalar there is $d I$, and the determinant-scalar
proposition (\Cref{prop-gl-determinant-scalar}) gives $d I \in C^n$, so
$C^n=\mathcal{D}(d^n)$.
  \item \textbf{$\lvert F \rvert\ge 4$, $\mu_n(F)\neq \{1\}$.} Again $C^n$ contains every nonscalar in
$\mathcal{D}(d^n)$, so the one-letter completion lemma
(\Cref{lem-one-letter-completion}) gives $C^{n+1}=\mathcal{D}(d^{n+1})$.
  \item \textbf{$F=\mathbb{F}_2$.} Small-field positive coverage
(\Cref{prop-sl-small-field-coverage}(i)) gives $C^n=\mathrm{GL}_n(\mathbb{F}_2)$, the
only fibre.
  \item \textbf{$F=\mathbb{F}_3$.} \Cref{prop-gl-f3-fibre} gives $n$ letters in odd dimension and
$n+1$ in even dimension; since $-1$ is the only nontrivial root of unity in
$\mathbb{F}_3$, these lengths are precisely $k$.
\end{itemize}

  Thus $C^k=D_k(C)$ over every field. Together with the lower bound this gives
  equality, and the special linear diameter theorem (\Cref{thm-sl-diameter})
  identifies the same value for $\Delta(\mathrm{SL}_n(F))$.
\end{proof}

\begin{figure}[htbp]
  \centering
  \includegraphics[width=304pt]{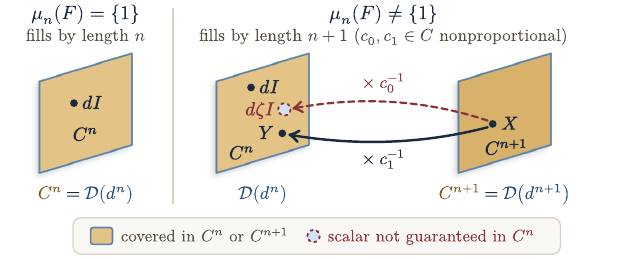}
  \caption{Completing the determinant fibre for $n \ge  3$ and $\lvert F \rvert \ge  4$.
    When $\mu_n(F) = \{1\}$ (left), $C^n$ already fills $\mathcal{D}(d^n)$; when
    $\mu_n(F) \neq  \{1\}$ (right), at least one of two nonproportional choices
    $c_0, c_1 \in C$ sends $X \in \mathcal{D}(d^{n+1})$ to a covered nonscalar
    $Y \in C^n$.}\label{fig-fibre-completion}
\end{figure}

\Cref{fig-fibre-completion} summarises the two broad-field branches of
\Cref{thm-gl-fibre-covering}.

\Cref{thm-gl-fibre-covering} concerns full general linear classes. Within the
special linear group, a conjugacy class can be smaller because only
determinant-one conjugators are allowed, yet it still fills $\mathrm{SL}_n(F)$ at the
same exact positive length.

\begin{corollary}[{Covering by one special linear class}]
\label{cor-sl-class-covering}
Let $F$ be a field and $n\ge 3$, and write $k=n+[\mu_n(F)\neq \{1\}]$. Then
  $C^k=\mathrm{SL}_n(F)$ for every noncentral conjugacy class $C$ of
  $\mathrm{SL}_n(F)$. No smaller exponent has this property for the
  transvection class.
\end{corollary}

\begin{proof}
Over $\mathbb{F}_2$, and for $n=4$ over $\mathbb{F}_3$, this is small-field positive
  coverage (\Cref{prop-sl-small-field-coverage}).

  In every other case $C^n$ contains every nonscalar of $\mathrm{SL}_n(F)$: for
  $\lvert F \rvert\ge 4$ by Lev's nonscalar coverage theorem
  (\Cref{prop-sl-nonscalar-coverage}, with $d=1$), and over $\mathbb{F}_3$ with $n\neq 4$
  by \Cref{prop-sl-f3-nonscalar-coverage}. Both results apply to $\mathrm{SL}_n(F)$-classes.
  If $\mu_n(F)\neq \{1\}$, the one-letter completion lemma
  (\Cref{lem-one-letter-completion}) gives $C^{n+1}=\mathrm{SL}_n(F)$.

  Suppose $\mu_n(F)=\{1\}$. The identity is then the only scalar in
  $\mathrm{SL}_n(F)$, so it remains to show $I \in C^n$. For $\lvert F \rvert\ge 4$ this is the
  $\mathrm{SL}_n(F)$-class clause of the determinant-scalar proposition
  (\Cref{prop-gl-determinant-scalar}) with $d=1$. If $F=\mathbb{F}_3$, then $n$ is
  odd, so $-I$ has determinant $-1$ and
  $\mathrm{GL}_n(\mathbb{F}_3)=\{\pm I\} \times \mathrm{SL}_n(\mathbb{F}_3)$; since conjugation
  by $-I$ is trivial, $C$ is a full conjugacy class of $\mathrm{GL}_n(\mathbb{F}_3)$, and
  \Cref{prop-gl-f3-fibre} gives $C^n=\mathrm{SL}_n(\mathbb{F}_3)$.

  Every positive word is a signed word, and the special linear diameter theorem
  (\Cref{thm-sl-diameter}) shows that some target has signed length $k$ from a
  transvection; that target lies in no shorter power of the transvection class.
\end{proof}

\begin{remark}[{Dimension three}]
\label{rem-sl3-identity}
For a noncentral class $C$ in $\mathrm{SL}_3(F)$, write $\overline{C}$ for its image
  in $\mathrm{PSL}_3(F)$. Over any infinite field $F$, Ambrosiewicz
  \cite[Corollary 6.1]{Ambrosiewicz1996} proved $\overline{C}^4 = \mathrm{PSL}_3(F)$. Under
  his standing hypothesis $|F|\ge 4$, Lev \cite[Theorem 1]{Lev1996} proved
  $\overline{C}^3 = \mathrm{PSL}_3(F)$ for finite and algebraically closed fields (see
  also Orevkov \cite[Corollary 1.9]{Orevkov2013} for finite fields) and left open
  whether the projective identity belongs to $\overline{C}^3$ when $C$ consists
  of irreducible matrices over a general infinite field \cite[p. 61]{Lev1996}. For
  every non-rank-one class $C$ in $\mathrm{SL}_3(F)$ over $|F|\ge 4$, the rotation word
  in \Cref{prop-gl-determinant-scalar} has conjugators of determinant one
  unconditionally, supplying the stronger conclusion $I \in C^3$ (not merely
  that $C^3$ meets the centre), while rank-one classes are already covered by
  \Cref{thm-prescribed-rank-one}.

  Passing to the quotient yields $\operatorname{cn}_\mathrm{class}(\mathrm{PSL}_3(F)) = 3$ over
  every field $F$: nonscalar coverage (\Cref{prop-sl-nonscalar-coverage}) and
  \Cref{prop-gl-determinant-scalar} give $\overline{C}^3 = \mathrm{PSL}_3(F)$ for
  $|F|\ge 4$, while \Cref{prop-sl-small-field-coverage}(i) and \Cref{prop-gl-f3-fibre}
  cover $\mathbb{F}_2$ and $\mathbb{F}_3$. Sharpness survives quotienting: for any
  nonzero $b \in F$, the matrix $\operatorname{Comp}(t^3 + b t - 1) \in \mathrm{SL}_3(F)$ has
  characteristic polynomial taking the nonzero value $b \zeta$ at every
  $\zeta \in \mu_3(F)$, so every central rescaling has deviation rank three at
  $1$ and no lift is a product of two transvections.
\end{remark}

\section{Neighbouring invariants and dependence on the field}
\label{sec-invariants-and-fields}
With the special linear diameter (\Cref{thm-sl-diameter}), determinant-fibre
covering (\Cref{thm-gl-fibre-covering}) and single-class special linear covering
(\Cref{cor-sl-class-covering}) established, this section interprets our common
formula $n + [\mu_n(F) \neq  \{1\}]$ in two directions: first by placing
$\Delta(G)$ and $\operatorname{cn}_\mathrm{ab}(G)$ alongside the neighbouring covering
invariants of $\mathrm{SL}_n(F)$ and $\mathrm{GL}_n(F)$
(\Cref{sec-invariants-literature}), and then by examining how the value
varies across families of fields (\Cref{sec-sl-field-examples}).

\subsection{Conjugation diameter and covering numbers}
\label{sec-invariants-literature}
The broader family of covering invariants grows out of the classical covering
problem: the \emph{classical covering number} $\operatorname{cn}_\mathrm{class}(G)$ is the least
$k$ such that $C^k=G$ for every nonidentity conjugacy class $C$
\cite{AradHerzog1985} \cite[\S~1]{VasersteinWheland1995}, and the corresponding signed
question asks for $B_C(k)=G$, allowing inverses and shorter words. Both
questions are natural for finite nonabelian simple groups, where every
nonidentity class normally generates the group.

Outside that setting, two obstructions can force the classical answer to be
$\infty$: a nonidentity class may generate only a proper normal subgroup,
and nontrivial abelianisation confines $C^k$ to a single coset of $[G,G]$. To avoid these degeneracies without losing the underlying covering question, one
can either restrict the seeds (keeping $G$ as the target) or restrict the target
to a single coset of $[G,G]$. This gives the remaining three columns of
\Cref{tab:word-conventions}:

\begin{itemize}
  \item \textbf{Requiring a single normally generating seed (whole-group target).} Keeping
$G$ as the target and restricting to single seeds or classes that normally
generate $G$ gives the \emph{single-seed conjugation diameter} $\Delta_1(G)$ for
signed words and the \emph{normally generating covering number}
$\operatorname{cn}_\mathrm{ng}(G)$ for positive powers.
  \item \textbf{Allowing multiple seeds or classes (whole-group target).} In general, normal
generation is a property of sets rather than single elements: a group may
require several seeds together to normally generate it, and even when a single
normal generator exists, a finite normally generating set need not contain one
(\Cref{fig-determinant-fibres}). Taking all finite normally generating sets $S$ for
signed words gives the conjugation diameter $\Delta(G)$. On the positive side,
allowing each factor to come from a different normally generating class rather
than repeating one class gives the \emph{extended covering number} $\operatorname{ecn}(G)$,
the least $k$ such that $C_1 \dots C_k=G$ for all normally generating classes
$C_1, \dots, C_k$.
  \item \textbf{Restricting the target (abelianisation fibre).} Alternatively, keeping a
single class $C$ whose normal closure contains $[G,G]$ and targeting only the
permitted coset of $[G,G]$ gives $\operatorname{cn}_\mathrm{ab}(G)$ from
\Cref{sec-covering-definitions} for positive powers. Its signed counterpart is the
worst \emph{signed fibrewise diameter}, measuring distances within a coset of
$[G,G]$ by the signed word length of quotients in $[G,G]$.
\end{itemize}

\Cref{tab:word-conventions} summarises these eight invariants.\footnote{Across all eight questions, we assign $\infty$ when no common finite bound
  exists or the admissible family is empty.}

\begin{table}[htbp]
  \centering
  \small
  \renewcommand{\arraystretch}{1.2}
  \setlength{\tabcolsep}{4pt}
  \begin{tabular}{@{}
    >{\centering\arraybackslash}p{2.15cm}
    >{\centering\arraybackslash}p{2.55cm}
    >{\centering\arraybackslash}p{2.4cm}
    >{\centering\arraybackslash}p{2.45cm}
    >{\centering\arraybackslash}p{2.4cm}@{}}
    \toprule
    \rowcolor{tabband}
    & {\bfseries Single\par seed}\par\vspace{3pt}\includegraphics[width=1.65cm]{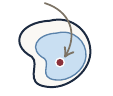}
    & {\bfseries Must normally\par generate}\par\vspace{3pt}\includegraphics[width=1.65cm]{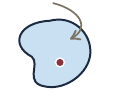}
    & {\bfseries Multiple\par seeds}\par\vspace{3pt}\includegraphics[width=1.65cm]{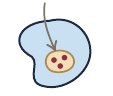}
    & {\bfseries Only correct\par fibre}\par\vspace{3pt}\includegraphics[width=1.65cm]{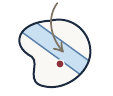} \\
    \midrule
    {\bfseries Nested\par lengths}\par{\footnotesize\color{tabmuted}(signed)}\par\vspace{4pt}\includegraphics[width=1.65cm]{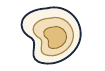}\par\vspace{2pt}$B_S(k)$
    & Worst single-seed diameter\par\vspace{4pt}$\sup\limits_{s \neq 1} \operatorname{diam}_{\{s\}}(G)$\par\vspace{4pt}{\footnotesize\color{tabmuted}\cite[\S 2]{KedraLibmanMartin2023}}
    & Single-seed conjugation diameter\par\vspace{4pt}$\Delta_1(G)$\par\vspace{4pt}{\footnotesize\color{tabmuted}\cite{KedraLibmanMartin2023}}
    & \multicolumn{1}{!{\color[HTML]{C89B49}\vrule width 1.5pt}>{\centering\arraybackslash}p{2.45cm}}{\cellcolor{tabfocal}Conjugation diameter\par\vspace{4pt}$\Delta(G)$\par\vspace{4pt}{\footnotesize\color{tabmuted}\cite{KedraLibmanMartin2023,LibmanTarry2020}}}
    & Fibrewise diameter\textsuperscript{$\dagger$}\par\vspace{4pt}{\footnotesize\color{tabmuted}worst over\par admissible seeds}\par\vspace{4pt}{\footnotesize\color{tabmuted}\cite[\S 2]{KedraLibmanMartin2023}} \\
    \midrule
    {\bfseries Strict\par lengths}\par{\footnotesize\color{tabmuted}(unsigned)}\par\vspace{4pt}\includegraphics[width=1.75cm]{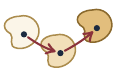}\par\vspace{2pt}$C^k$
    & Classical covering number\par\vspace{4pt}$\operatorname{cn}_{\mathrm{class}}(G)$\par\vspace{4pt}{\footnotesize\color{tabmuted}\cite[\S 1]{VasersteinWheland1995}}
    & Normally generating covering number\par\vspace{4pt}$\operatorname{cn}_{\mathrm{ng}}(G)$\par\vspace{4pt}{\footnotesize\color{tabmuted}\cite{Preusser2022}}
    & Extended covering number\par{\footnotesize\itshape\color{tabmuted}prescribed classes}\par\vspace{2pt}$\operatorname{ecn}(G)$\par\vspace{4pt}{\footnotesize\color{tabmuted}\cite{Preusser2022,KnuppelNielsen2006}}
    & \multicolumn{1}{!{\color[HTML]{C89B49}\vrule width 1.5pt}>{\centering\arraybackslash}p{2.4cm}@{}}{\cellcolor{tabfocal}Abelianised covering number\par\vspace{4pt}$\operatorname{cn}_{\mathrm{ab}}(G)$\par\vspace{4pt}{\footnotesize\color{tabmuted}\cite[\S 2]{VasersteinWheland1995}}} \\
    \bottomrule
  \end{tabular}
  \par\medskip
  {\footnotesize\color{tabmuted}\raggedright
  \textsuperscript{$\dagger$} Supremum over admissible seeds of the restricted word-metric diameter \cite[\S 2]{KedraLibmanMartin2023}.\newline
  \begingroup\setlength{\fboxsep}{0pt}\colorbox{tabfocal}{\textcolor[HTML]{C89B49}{\rule{1.5pt}{6.5pt}}\rule{6pt}{0pt}}\endgroup\hspace{2.5pt}Main Part~I invariants: conjugation diameter $\Delta(G)$ and abelianised covering number $\operatorname{cn}_{\mathrm{ab}}(G)$.\par}
  \caption{Eight conjugation-diameter and covering-number invariants, arranged by word convention (rows) and seed or target restriction (columns).}
  \label{tab:word-conventions}
\end{table}

Combining \Cref{sec-matrix-normal-generation} with \Cref{thm-sl-diameter},
\Cref{thm-gl-fibre-covering} and \Cref{cor-sl-class-covering} gives exact values,
finiteness criteria and upper bounds for the entries of \Cref{tab:word-conventions}
for $\mathrm{SL}_n(F)$ and $\mathrm{GL}_n(F)$ (the finite-field whole-group diameters of
$\mathrm{GL}_n(F)$ are determined in \mbox{Part II}).

In $\mathrm{SL}_n(F)$, perfectness identifies each abelianisation fibre with the
whole group, and every finite normally generating set contains a noncentral
element. Consequently, all of the normally generating signed and
single-class positive invariants collapse to the same value
$n+[\mu_n(F)\neq \{1\}]$, even though exact positive covering must construct the
identity explicitly rather than assigning it length zero. The only exception in
the single-class entries occurs when central seeds are permitted, which forces
$\infty$ whenever the centre is nontrivial.

\begin{proposition}[{Comparison of invariants for $\mathrm{SL}_n(F)$}]
\label{prop-comparison-sl}
Let $F$ be a field and $n\ge 3$.

\begin{enumerate}[label=(\roman*)]
  \item \textbf{(Signed diameters.)} We have
      \[
  \Delta_1(\mathrm{SL}_n(F))=\Delta(\mathrm{SL}_n(F))=n+[\mu_n(F)\neq \{1\}],
\]
      and the worst signed fibrewise diameter has the same value. The supremum
      of single-seed diameters over all nonidentity seeds is $\infty$ if
      $\mu_n(F)\neq \{1\}$, and equals $\Delta_1(\mathrm{SL}_n(F))$ otherwise.
  \item \textbf{(Exact positive covering.)} We have
      \[
  \operatorname{cn}_\mathrm{ab}(\mathrm{SL}_n(F)) = \operatorname{cn}_\mathrm{ng}(\mathrm{SL}_n(F))
        = \Delta(\mathrm{SL}_n(F)) = n+[\mu_n(F)\neq \{1\}],
\]
      and
      \[
  \operatorname{cn}_\mathrm{class}(\mathrm{SL}_n(F))=
        \begin{cases}\infty & \quad \textnormal{if } \mu_n(F)\neq \{1\} \\ \operatorname{cn}_\mathrm{ng}(\mathrm{SL}_n(F)) & \quad \textnormal{if } \mu_n(F)=\{1\}.\end{cases}
\]
\end{enumerate}
\end{proposition}

\begin{proof}
For (i), the normal generation proposition (\Cref{prop-sl-normal-generation})
  gives $\Delta_1=\Delta$, perfectness identifies the signed fibrewise and
  whole-group diameters, central seeds reach only central elements, and the
  special linear diameter theorem (\Cref{thm-sl-diameter}) gives the common value.

  For (ii), a class in $\mathrm{SL}_n(F)$ normally generates if and only if it is
  noncentral, and perfectness gives $\operatorname{cn}_\mathrm{ab} = \operatorname{cn}_\mathrm{ng}$, yielding
  the formula for $\operatorname{cn}_\mathrm{class}$. Since every positive word is a signed
  word, (i) gives the lower bound for $\operatorname{cn}_\mathrm{ng}$, while single-class
  special linear covering (\Cref{cor-sl-class-covering}) supplies the matching
  upper bound.
\end{proof}

\begin{remark}[{The extended covering number}]
\label{rem-ecn-sl}
Prescribing the classes independently changes the behaviour even when the
  centre is trivial, because some sequences of $n$ classes have no product
  equal to the identity: for instance, in the centreless group
  $\mathrm{SL}_3(\mathbb{F}_3) = \mathrm{PSL}_3(\mathbb{F}_3)$, Orevkov \cite[Corollary 1.9]{Orevkov2013}
  gives $\operatorname{ecn} = 4$, whereas $\operatorname{cn}_\mathrm{ng} = 3$. For $n\neq 4$ and
  $F\neq \mathbb{F}_2$, Kn\"uppel and Nielsen state that
  $\operatorname{ecn}(\mathrm{SL}_n(F))=n+1$ \cite[p. 635]{KnuppelNielsen2006}. In the printed
  proof, the last case of Proposition 8.2 leaves one class out of the product
  it constructs \cite[p. 654, step (vi)]{KnuppelNielsen2006}, and Propositions
  7.4--7.5 and Corollaries 7.6 and 7.10 are deduced via an abbreviated pointer
  to Lemmas 7.3 and 7.7 \cite[p. 649]{KnuppelNielsen2006}. We use their theorem
  only for $n$ copies of one class together with one further class, where
  Proposition 8.2(vi) does not arise, as explained in the remark following
  \Cref{prop-sl-f3-nonscalar-coverage}.
\end{remark}

In $\mathrm{GL}_n(F)$, by contrast, the choice of invariant changes the answer
sharply whenever $F\neq \mathbb{F}_2$ (over $\mathbb{F}_2$, $\mathrm{GL}_n = \mathrm{SL}_n$ and
\Cref{prop-comparison-sl} applies). Every fixed positive product has a single
determinant, so all three whole-group positive covering numbers are infinite;
whole-group signed diameters are finite only over finite fields, where
\mbox{Part II} determines $\Delta(\mathrm{GL}_n(F))$. Restricting the
target to a determinant fibre, however, yields finite bounds over every field,
with $\operatorname{cn}_\mathrm{ab}(\mathrm{GL}_n(F))$ matching $\Delta(\mathrm{SL}_n(F))$ exactly.

\begin{proposition}[{Comparison of invariants for $\mathrm{GL}_n(F)$}]
\label{prop-comparison-gl}
Let $F$ be a field and $n\ge 3$.

\begin{enumerate}[label=(\roman*)]
  \item \textbf{(Signed diameters.)} If \mbox{$F\neq \mathbb{F}_2$}, the unrestricted supremum
      of single-seed diameters of $\mathrm{GL}_n(F)$ is $\infty$. Both
      $\Delta_1(\mathrm{GL}_n(F))$ and $\Delta(\mathrm{GL}_n(F))$ are finite if and only if
      $F$ is finite. Over every field, every noncentral seed gives signed
      diameter at most $2(n+1)$ on each determinant fibre.
  \item \textbf{(Exact positive covering.)} If $F\neq \mathbb{F}_2$, then
      \[
  \operatorname{cn}_\mathrm{class}(\mathrm{GL}_n(F))
        =\operatorname{cn}_\mathrm{ng}(\mathrm{GL}_n(F))=\operatorname{ecn}(\mathrm{GL}_n(F))=\infty.
\]
      Over every field,
      \[
  \operatorname{cn}_\mathrm{ab}(\mathrm{GL}_n(F))=\Delta(\mathrm{SL}_n(F))=n+[\mu_n(F)\neq \{1\}].
\]
\end{enumerate}
\end{proposition}

\begin{proof}
For (i), if $F\neq \mathbb{F}_2$ then a determinant-one seed cannot reach other
  determinant fibres. By the finite-field boundary
  (\Cref{cor-gl-infinite-field}), normally generating singletons exist over finite
  fields, where finiteness of the group bounds the diameters, whereas over
  infinite fields no finite set normally generates. For any noncentral class
  $C$, the determinant-fibre covering theorem (\Cref{thm-gl-fibre-covering}) gives
  $C^k=D_k(C)$ for some $k\le n+1$; taking quotients yields
  $C^k (C^k)^{-1}=\mathrm{SL}_n(F)$, bounding signed distances within each fibre by
  $2k$.

  For (ii), every prescribed positive product has a single determinant and
  therefore cannot cover $\mathrm{GL}_n(F)$ when $F^\times\neq \{1\}$ (and empty admissible
  families give $\infty$ by convention). The formula for
  $\operatorname{cn}_\mathrm{ab}(\mathrm{GL}_n(F))$ is \Cref{thm-gl-fibre-covering}, which agrees with
  \Cref{prop-comparison-sl}.
\end{proof}

\subsection{Dependence on the field}
\label{sec-sl-field-examples}
We now evaluate the common formula $n + [\mu_n(F) \neq  \{1\}]$ governing the signed
conjugation diameter $\Delta(\mathrm{SL}_n(F))$, the single-class positive covering
number $\operatorname{cn}_\mathrm{ng}(\mathrm{SL}_n(F))$, and the determinant-fibre covering number
$\operatorname{cn}_\mathrm{ab}(\mathrm{GL}_n(F))$ across the field families summarised in
\Cref{tab-intro-fields}, and spell out what the bound means over local and function
fields. In every example below, the same number gives all three invariants
(whereas finite-field whole-group diameters of $\mathrm{GL}_n(F)$ can be larger, as
studied in \mbox{Part II}).

\begin{example}[{Rational, real and complex scalars}]
\label{ex-field-real-complex}
Over $\mathbb{Q}$ or $\mathbb{R}$ the only candidate for a nonidentity scalar of
  determinant one is $-I$, whose determinant is $(-1)^n$; the common value is
  therefore $n$ in odd dimension and $n+1$ in even dimension. Over $\mathbb{C}$ a
  primitive $n$th root of unity supplies such a scalar in every dimension
  $n\ge 3$, so the common value is always $n+1$.
\end{example}

\begin{example}[{Finite fields}]
\label{ex-field-finite}
For any prime power $q$, the multiplicative group $\mathbb{F}_q^\times$ is cyclic
  of order $q-1$ \cite[Exercise 1-3(e)]{Milne2022}, so exactly \mbox{$\gcd(n,q-1)$}
  elements satisfy $\zeta^n=1$. The extra letter is therefore required precisely
  when $\gcd(n,q-1)>1$. Enlarging the field at fixed dimension and
  characteristic does not increase the $n$-letter scale, but it can supply new
  roots of unity: in dimension three, the common value is $3$ over $\mathbb{F}_2$
  and $4$ over $\mathbb{F}_4$. Over $\mathbb{F}_2$, it equals $n$ in every dimension $n\ge 3$.
\end{example}

Over a $p$-adic field, a transvection can be arbitrarily close to the identity;
the bound nevertheless stays small because conjugation in $\mathrm{SL}_n(\mathbb{Q}_p)$
allows changes of scale.

\begin{example}[{A transvection tending to the identity}]
\label{ex-field-padic}
In $\mathrm{SL}_n(\mathbb{Q}_p)$ take $g_a=I+p^a E_{1 2}$ for positive integers $a$,
  where $E_{1 2}$ has its single nonzero entry $1$ in position $(1,2)$. As $a$
  increases, $g_a$ tends to $I$, yet the determinant-one change of basis
  \[
  P_a=\operatorname{diag}(p^{-a},1,p^a) \oplus  I_{n-3},
    \quad P_a g_a P_a^{-1}=I+E_{1 2},
\]
  (with the last block absent for $n=3$) shows that every $g_a$ lies in the
  same transvection class and produces the same word lengths. As the shear
  shrinks, the denominator $p^a$ in the conjugator grows.

  For an odd prime $p$, the roots of unity in $\mathbb{Q}_p$ form a cyclic group of
  order \mbox{$p-1$}, whereas in $\mathbb{Q}_2$ they are $1$ and $-1$
  \cite[\S~7 and solution to Exercise 7-5]{Milne2020}. In dimension three, every
  determinant-one matrix over $\mathbb{Q}_2$ is therefore a product of at most three
  signed conjugates of any prescribed noncentral seed, while over $\mathbb{Q}_7$ the
  bound is four, attained at a scalar of order three for a transvection seed.
  These are exact algebraic identities over $\mathbb{Q}_p$; in $\mathrm{SL}_n(\mathbb{Z}_p)$ or
  its congruence subgroups, conjugators must be integral and the rescaling
  above is unavailable.
\end{example}

When matrix entries depend rationally on parameters, applying \Cref{thm-sl-diameter}
over the function field produces a single rational family of factorisations
rather than unrelated choices at individual parameter values.

\begin{example}[{A rational family of factorisations}]
\label{ex-field-rational-family}
Fix a noncentral seed $g(t)$ and a target $X(t)$ in $\mathrm{SL}_n(K)$, where
  $K=k(t)$. \Cref{thm-sl-diameter} gives an identity
  \[
  X(t)=\prod_{i=1}^m P_i(t) g(t)^{\varepsilon_i} P_i(t)^{-1},
    \quad m\le n+[\mu_n(K)\neq \{1\}],
\]
  with $\varepsilon_i \in \{1,-1\}$ and each $P_i(t) \in \mathrm{SL}_n(K)$, so the
  conjugating matrices themselves have entries rational in $t$ over $k$.
  Choosing a common nonzero denominator for the entries of $g(t)$, $X(t)$ and
  the $P_i(t)$, the determinant-one conditions and the product identity
  specialise wherever that denominator does not vanish. The bound controls the
  number of factors, not the degrees of the rational entries or the exceptional
  values of $t$.

  Moreover, adjoining independent indeterminates introduces no new roots of
  unity: if $a / b \in k(t)$ with $a,b \in k[t]$ coprime satisfies
  $a^n=b^n$, then neither polynomial can have an irreducible factor, so both
  are constant. By induction in several variables,
  \[
  \mu_n(k(t_1,\dots,t_r)) = \mu_n(k).
\]
  The invariants over $k(t_1,\dots,t_r)$ therefore agree with those over $k$:
  for instance, in dimension three over $\mathbb{Q}(t)$ or $\mathbb{F}_2(t)$ at most
  three factors suffice, whereas over $\mathbb{F}_4(t)$ or $\mathbb{C}(t)$ the bound is
  four.
\end{example}

\begin{example}[{A cube root from a relation}]
\label{ex-field-related-parameters}
An algebraic relation between parameters can, however, enlarge the field of
  constants. The equation \mbox{$u^2+u v+v^2=0$} defines an integral affine
  curve over $\mathbb{Q}$ in whose function field the ratio $u / v$ satisfies
  $x^2+x+1=0$, giving a primitive cube root of unity; in dimension three the
  common value rises from $3$ (over $\mathbb{Q}(t)$) to $4$. More generally, for an
  integral parameter space $S$ over $k$, every root of unity in $k(S)$ lies in
  the full constant field (the algebraic closure of $k$ in $k(S)$), and
  clearing denominators on an affine open of $S$ yields the factorisation on a
  nonempty open subset (which need not contain a $k$-rational point).
\end{example}

\begin{example}[{Algebraically closed fields of characteristic $p$}]
\label{ex-field-closed}
Over $\mathbb{Q}$ a primitive cube root arose by adjoining a root of
  \mbox{$x^2+x+1=0$}. In characteristic three the same polynomial is
  $(x-1)^2$, so even the algebraic closure contains no nontrivial cube root.
  More generally, let $F$ be algebraically closed of characteristic $p>0$, and
  write $n=p^a m$ with $p$ not dividing $m$. Then
  \[
  x^n-1=(x^m-1)^{p^a}
\]
  has exactly $m$ distinct roots in $F$. Thus $\mu_n(F)=\{1\}$ precisely when $n$
  is a power of $p$, in which case the common value equals $n$; in every other
  dimension $n\ge 3$ it equals $n+1$.
\end{example}

\begin{figure}[htbp]
  \centering
  \includegraphics[width=368pt]{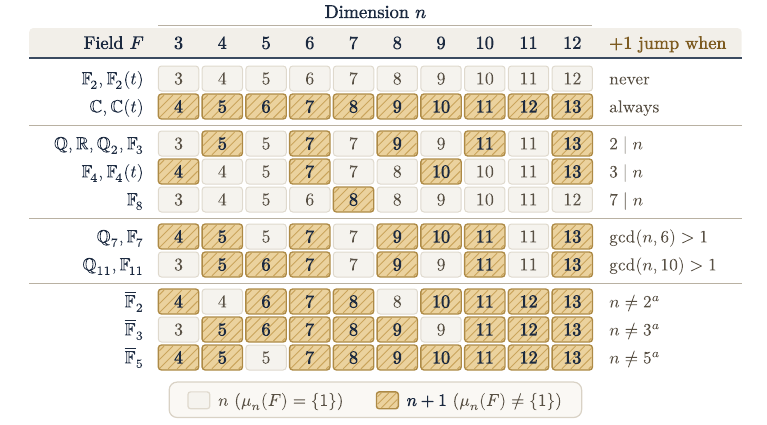}
  \caption{The common value \mbox{$n + [\mu_n(F) \neq  \{1\}]$} for
    $\Delta(\mathrm{SL}_n(F)) = \operatorname{cn}_\mathrm{ab}(\mathrm{GL}_n(F))$ for different example
    fields when $3 \le  n \le  12$.}\label{fig-field-dependence}
\end{figure}

\Cref{fig-field-dependence} summarises these patterns across dimensions
$3 \le  n \le  12$, grouping fields that share the same roots of unity.

Over a finite field one can go on to ask how many letters it takes to reach all
determinants as well as fill each fibre. Part II takes up this whole-group
diameter problem for $\mathrm{GL}_n(\mathbb{F}_q)$ and returns to dimension two, where
quadratic arithmetic enters even within a single fibre.

\appendix
\section{Small-field proofs}
\label{sec-small-field-proofs}
This appendix supplies the small-field positive coverage results used in
\Cref{sec-class-products}: exact $n$-letter coverage of $\mathrm{GL}_n(\mathbb{F}_2)=\mathrm{SL}_n(\mathbb{F}_2)$
(\Cref{prop-sl-small-field-coverage}(i)), five-letter coverage of
$\mathrm{SL}_4(\mathbb{F}_3)$ from a prescribed special linear class
(\Cref{prop-sl-small-field-coverage}(ii)), and determinant-fibre coverage of
$\mathrm{GL}_n(\mathbb{F}_3)$ (\Cref{prop-gl-f3-fibre}). All three proofs draw on a shared
cyclic construction developed in
\Cref{sec-small-field-cyclic-constructions}: once a short power $C^j$
contains a full cyclic class, two such blocks and one extra letter reach every
nonscalar of the matching determinant. The length of those blocks then decides
whether a further letter remains available for scalar targets, whereas on a
scalar-free determinant fibre the cyclic construction already finishes the
proof.

\subsection{Shared cyclic constructions}
\label{sec-small-field-cyclic-constructions}
Because the similarity class of a cyclic matrix is determined by a single
polynomial, products of two full cyclic classes are particularly effective: when
one class has a simple eigenvalue in the ground field, their product reaches
every nonscalar matrix of the prescribed determinant. We apply this cyclic
decomposition twice, first to obtain a partner for a cyclic class already
reached and then to reach the target itself.

\begin{definition}[{Cyclic and nice cyclic classes}]
\label{def-nice-cyclic}
A matrix $g$ is \emph{cyclic} if some vector $v$ has
  $v,g v,\dots,g^{n-1}v$ as a basis. A cyclic matrix is \emph{nice cyclic}
  if it has a simple eigenvalue in the ground field. We use the same
  terms for conjugacy classes of such matrices. A \emph{full} class is taken
  under general linear conjugation.
\end{definition}

\begin{proposition}[{Nonscalar coverage from a cyclic power}]
\label{prop-cyclic-power}
Let $n\ge 3$ and let $C$ be the conjugacy class of a nonscalar matrix
  of determinant $d$ under $\mathrm{SL}_n(F)$ or $\mathrm{GL}_n(F)$. If $C^j$
  contains a full cyclic class $\mathcal{K}$ for some $j\ge 1$, then
  \[
  \{X \in \mathrm{GL}_n(F) : \det X=d^{2j+1}, X \textnormal{ nonscalar}\} \subseteq C^{2j+1}.
\]
\end{proposition}

\begin{figure}[htbp]
  \centering
  \includegraphics[width=274pt]{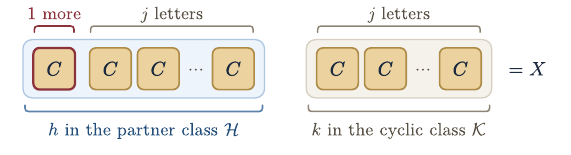}
  \caption{Splitting a nonscalar target $X$ of determinant $d^{2j+1}$ into a
    nice cyclic partner $h \in \mathcal{H}$ ($j + 1$ letters) and a cyclic factor
    $k \in \mathcal{K}$ ($j$ letters), giving $X = h k \in C^{2j+1}$.}\label{fig-appendix-cyclic-count}
\end{figure}

One block of $j$ letters supplies a full cyclic class, and spending one further
letter on a second block lets us pick a nice cyclic partner of the matching
determinant (\Cref{fig-appendix-cyclic-count}). Full-class containment is essential
here: conjugating a product of prescribed special linear classes yields its
special linear orbit, and the simple eigenvalue guarantees that this orbit
already equals the partner's full general linear class.

\begin{lemma}[{Choosing a nice cyclic partner}]
\label{lem-cyclic-partner}
With $C$ as above, let $\mathcal{K}$ be any full cyclic class and
  $\mathcal{H}$ any nice full cyclic class with $\det \mathcal{H}=d \det \mathcal{K}$.
  Then $\mathcal{H} \subseteq C \mathcal{K}$.
\end{lemma}

\begin{proof}
Fix $c \in C$. The cyclic decomposition theorem of B\"unger and Nielsen
  \cite[Theorem 1.2]{BungerNielsen1999}, in the prescribed-class form of
  \cite[Proposition 4.2]{KnuppelNielsen2006} (see the remark below for the required
  corrections to the printed proof), applies to $\mathcal{H}$ and
  $\mathcal{K}^{-1}$: the first is nice, both are cyclic, and their determinants
  multiply to $\det c$. Hence $c=h k^{-1}$ for some $h \in \mathcal{H}$ and
  $k \in \mathcal{K}$, giving $h=c k$ (the factor reversal of
  \cite[Lemma 4.6]{KnuppelNielsen2006}). Because the simple eigenline of $h$ admits
  an $h$-invariant complement, scaling that line while fixing the complement
  realises every nonzero determinant in the centraliser of $h$. Its special
  linear orbit therefore coincides with its full class $\mathcal{H}$, and since
  $C \mathcal{K}$ is $\mathrm{SL}_n(F)$-invariant and contains $h$, it contains $\mathcal{H}$.
\end{proof}

\begin{proof}[{Proof of \Cref{prop-cyclic-power}}]
Recall from the proof of \Cref{prop-gl-determinant-scalar} that $\operatorname{Comp}(f)$
  denotes the companion matrix of a monic polynomial $f$ with ones on the
  subdiagonal. There is a nice cyclic class of every determinant $e \in F^\times$:
  take $\operatorname{Comp}((x-1)(x^{n-1}+u x+a))$ with $a=(-1)^{n+1}e$. Its determinant is
  $e$, and $1$ is a simple root provided $1+u+a\neq 0$; since $n\ge 3$, the linear
  coefficient $u$ is free and only one choice is excluded. Choosing such a
  class $\mathcal{H}$ of determinant $e=d^{j+1}$, the nice cyclic partner lemma
  (\Cref{lem-cyclic-partner}) gives $\mathcal{H} \subseteq C \mathcal{K} \subseteq C^{j+1}$.
  Applying the same cyclic decomposition theorem to $\mathcal{H}$ and $\mathcal{K}$ then
  puts every nonscalar of determinant $d^{2j+1}$ in
  $\mathcal{H} \mathcal{K} \subseteq C^{2j+1}$.
\end{proof}

\begin{remark}[{Corrections to the cited proof of B\"unger and Nielsen}]
\label{rem-bunger-nielsen-corrections}
Both \Cref{prop-gl-determinant-scalar} and \Cref{lem-cyclic-partner} use the
  following consequence of \cite[Theorem 1.2]{BungerNielsen1999} (cited via
  \cite[Proposition 4.2]{KnuppelNielsen2006}): for $n\ge 3$ and full cyclic classes
  $\mathcal{H}, \mathcal{K}$ in $\mathrm{GL}_n(F)$ where $\mathcal{H}$ has a simple eigenvalue
  $\lambda \in F$, the product $\mathcal{H} \mathcal{K}$ contains every nonscalar $M$ of
  determinant $\det \mathcal{H} \det \mathcal{K}$. Indeed, setting
  $q(x) = \chi_{\mathcal{H}} (x) / (x - \lambda)$ and prescribing $\mathcal{K}$ as the
  second class, \cite[Theorem 1.2]{BungerNielsen1999} yields $M = A B$ with
  $B \in \mathcal{K}$ and $q \mid m_A$. Since $\deg q = n-1$ and
  $\det A = \det \mathcal{H}$, the remaining linear factor of $\chi_A$ must be
  $x - \lambda$; as $q(\lambda) \neq  0$, both $q$ and $x - \lambda$ divide $m_A$, so
  $m_A = \chi_{\mathcal{H}}$ and $A \in \mathcal{H}$.

  Three technical slips in the printed proofs of Theorems 4.2 and 1.2 in
  \cite[pp. 47--49]{BungerNielsen1999} (which also underlie
  \cite[Propositions 4.2 and 4.8]{KnuppelNielsen2006} and \Cref{lem-f3-det-scalar} below)
  are repaired as follows:
\begin{enumerate}[label=(\roman*)]
  \item \textbf{Theorem 4.2 (pp. 47--48).} Both coprime factors of $\chi_M$ must have
      positive degree (a condition already stated in
      \cite[Proposition 4.8]{KnuppelNielsen2006} and satisfied in
      \Cref{lem-f3-det-scalar}). In addition, left multiplication by the cycle
      matrix $R$ multiplies both determinants by $\varepsilon = \det R = (-1)^{n-1}$.
      Allocating $\det B \det B' = \varepsilon \det \Delta$ before multiplying by $R$
      (setting the last diagonal entry of $B'$ to $\varepsilon \det \Delta$ and
      adjusting the last entry of $C'$ by the reciprocal sign) preserves the
      block products and gives the prescribed determinants in both classes.
  \item \textbf{Theorem 1.2, Case A (eq. (5.1)) and Case B2.} The bordered matrix $E$
      requires $\beta = \operatorname{tr} \Omega$, which can vanish even when $\Omega$ is
      invertible. Allowing $\beta \in F$ causes no obstruction: since $D$ is a
      single upper nilpotent Jordan block with $d_{n-1} \neq  0$, the vectors
      $d, D d, \dots, D^{n-2} d$ form a basis, so choosing the bordering row $e$
      by its values on this basis prescribes all remaining coefficients of
      $\chi_E(x) = (x - \beta) x^{n-1} - \sum_{j=0}^{n-2} (e D^j d) x^{n-2-j}$.
      The prescribed nonzero constant term guarantees invertibility and forces
      $e_1 \neq  0$, as required in the subsequent step.
  \item \textbf{Theorem 1.2, Case B (p. 49).} Comparing constant terms in
      $(x - \alpha) q(x)$ gives
      $\alpha = (-1)^{n-1} \det M / (q(0) \det \Omega)$ (dividing by $q(0)$,
      rather than multiplying as printed). In Case B1, taking the companion
      matrix of $(x - \alpha) q(x)$ fixes both its characteristic polynomial and
      its determinant; for our nice-cyclic consequence above, $\alpha = \lambda$
      with $q(\lambda) \neq  0$, so the factors are already coprime and corrected
      Theorem 4.2 applies directly.
\end{enumerate}
\end{remark}

To control the block length $j$, we measure progress towards cyclicity (having
a single nonconstant invariant factor) even when the seed has no eigenvalue in
the ground field.

\begin{definition}[{Minimum rank}]
\label{def-appendix-min-rank}
For a class $C$ in dimension $n$, let $\delta(C)$ be $n$ minus the
  number of nonconstant invariant factors. Equivalently, for $c \in C$,
  it is the minimum of $\operatorname{rank}(c-\lambda I)$ over the algebraic
  closure, so $\delta(C)=n-1$ precisely for a cyclic class. Following
  Kn\"uppel and Nielsen, a class is \emph{merely quadratic} if all its
  nonconstant invariant factors have degree two (divisibility then
  forces them all to be equal).
\end{definition}

\begin{proposition}[{Minimum-rank products of full classes}]
\label{prop-minimum-rank-product}
Let $n\ge 3$, and let $F$ be a field with $(n,F)\neq (4,\mathbb{F}_2)$. For
  noncentral full classes $C_1,\dots,C_r$ in $\mathrm{GL}_n(F)$, $r\ge 1$,
  the product contains a full class of minimum rank
  \[
  \min(n-1, \sum_{i=1}^r \delta(C_i)).
\]
  In particular, $j=\lceil (n-1)/\delta(C) \rceil$ copies of a noncentral full
  class $C$ suffice to obtain a full cyclic class in $C^j$.
\end{proposition}

Kn\"uppel and Nielsen obtain cyclicity by combining the cycles attached to the
invariant factors; keeping the constant-term coefficients on those cycles lets
the construction accommodate arbitrary determinants. Two merely quadratic
factors require a separate step, joining two cyclic blocks whose characteristic
polynomials are coprime.

\begin{proof}
Over $\mathbb{F}_2$ the full group is the special linear group, so
  \cite[Corollary 6.9]{KnuppelNielsen2006} applies at the stated dimensions.
  Suppose $F\neq \mathbb{F}_2$. We adapt the construction of Kn\"uppel and Nielsen
  \cite[\S\S~3, 5--6]{KnuppelNielsen2006}. Their later rank-addition result
  \cite[Proposition 5.2 and Corollary 5.3]{KnuppelNielsen2013} also covers
  arbitrary determinants when the sum is at most $n-1$. Here the sum
  can cross that threshold.

  \textbf{Permutation representatives.} Let $U$ be the upper unitriangular
  group. Choose a permutation with cycle lengths the degrees of the
  nonconstant invariant factors of a class. For each cycle
  $(i_1,\dots,i_l)$ take $i_l$ largest and let
  $x^l+a_{l-1}x^{l-1}+\dots+a_0$ be the corresponding invariant factor.
  Let a monomial matrix $M$ send $e_{i_j}$ to $e_{i_{j+1}}$ for
  $j<l$, and $e_{i_l}$ to $-a_0 e_{i_1}$. Let $u$ send $e_{i_l}$ to
  $e_{i_l}-\sum_{j=1}^{l-1} a_j e_{i_j}$ and fix the earlier vectors
  of the cycle. These choices on the disjoint cycles give $u \in U$.
  In the ordered cycle bases $M u$ has the prescribed companion blocks,
  including their constant terms.

  Two such representatives can be multiplied with the monomial parts
  together, since
  \[
  (M_1 u_1)(u_1^{-1} M_2 u_2 u_1)=M_1 M_2 u_2 u_1.
\]
  The second letter stays in its full class. Iteration gives the
  representatives used in \cite[3.3 and 3.7]{KnuppelNielsen2006}.
  The rearrangement in their Lemma 3.2 and the cyclicity criterion
  in 3.4 allow arbitrary determinants. Thus the permutation choices
  of 3.5--3.6 and 6.2--6.3 apply to these representatives: when a
  prefix first reaches rank sum $n-1$, they give a cyclic product,
  except for a pair of merely quadratic classes. A prefix of one
  factor at the threshold is already cyclic.

  \textbf{The quadratic pair.} For $n=2t$, write the minimum polynomials as
  $x^2-\alpha x-\beta$ and $x^2-\alpha' x-\beta'$ and put $\gamma=\beta \beta'$.
  The quadratic-pair construction \cite[Proposition 4.6]{KnuppelNielsen2013}
  puts $Z \oplus  \gamma Z^{-1}$ in the product of the full classes whenever
  $Z \in \mathrm{GL}_t(F)$ and the minimum polynomials of $Z$ and $\gamma Z^{-1}$ are
  coprime. If $Z$ is cyclic, their direct sum is then cyclic. For $\lvert F \rvert\ge 4$, take
  $Z=J_t(\lambda)$, a single Jordan block with $\lambda\neq 0$ and
  $\lambda^2\neq \gamma$. At most two nonzero choices are excluded. Over
  $\mathbb{F}_3$ with $\gamma=-1$, the same choice works with $\lambda=1$.

  Over $\mathbb{F}_3$ with $\gamma=1$, let $Z=\operatorname{Comp}(f)$ with $f=x^t-x-1$.
  The characteristic polynomial of $Z^{-1}$ is the monic reciprocal
  $g=x^t+x^{t-1}-1$. Since $x g-(x+1)f=(x-1)^2$ and $f(1)=-1$,
  these polynomials are coprime. In dimension four this gives a cyclic
  matrix with polynomial $x^4+1$.

  \textbf{Finishing the product.} Once a prefix contains a cyclic element it
  contains its full class. In dimension at least three that class is
  never merely quadratic. The same permutation construction combines it
  with each remaining noncentral factor to keep a cyclic product.
  Finally, if the original sum is $R<n-1$, append the full class of
  $J_{n-R}(1) \oplus  I_R$, of minimum rank $n-1-R$. The enlarged product
  contains a cyclic element $S T$ with $S$ in the original product
  and $T$ in the appended class. Subadditivity of minimum rank
  \cite[Lemma 5.1]{KnuppelNielsen2006} gives
  \[
  n-1=\delta(S T)\le \delta(S)+n-1-R\le n-1.
\]
  Thus $\delta(S)=R$, and simultaneous conjugation supplies its full class.
\end{proof}

\begin{remark}[{The published corollary in dimension four}]
\label{rem-kn-corollary-dimension-four}
Kn\"uppel and Nielsen state the corresponding result for special linear
  classes whenever $n\ge 3$ and $(n,F)\neq (4,\mathbb{F}_2)$
  \cite[Corollary 6.9]{KnuppelNielsen2006}. Over $\mathbb{F}_3$ in dimension four
  their proof rests on Lemma 6.10(c), whose printed proof treats only
  two classes with minimum polynomials $(x \pm 1)^2$. Pairs
  involving $x^2+1$, $x^2-1$ or $x^2 \pm x-1$ are not addressed.
  The quadratic-pair step above covers every pair of full classes at
  $(n,F)=(4,\mathbb{F}_3)$, so the proposition does not use that lemma.
\end{remark}

\subsection{Positive products over \texorpdfstring{$\mathbb{F}_2$}{F2}}
\label{sec-fibre-f2}
Over $\mathbb{F}_2$ there is only one determinant fibre, so the sole target left out
by nonscalar coverage is the identity itself. Whenever the cyclic construction
leaves even one letter to spare, the one-letter completion lemma
(\Cref{lem-one-letter-completion}) absorbs the identity automatically; only when
the length budget is tight do we need an explicit $n$-letter word for $I$.

\textbf{Restatement of \Cref{prop-sl-small-field-coverage}(i) (Positive coverage over
  $\mathbb{F}_2$).} For $n\ge 3$ and every noncentral conjugacy class $C$ in
  $\mathrm{GL}_n(\mathbb{F}_2)=\mathrm{SL}_n(\mathbb{F}_2)$,
  \[
  C^n=\mathrm{GL}_n(\mathbb{F}_2).
\]

\begin{proof}
Put $G=\mathrm{GL}_n(\mathbb{F}_2)$. In dimension four the exceptional isomorphism
  $G \cong A_8$ and Dvir's class-power covering theorem \cite{Dvir1985},
  as quoted by Lev \cite[pp. 60--61]{Lev1996}, give $C^4=G$ for every
  nontrivial class. We henceforth assume $n\neq 4$.

  Let $\delta(C)$ be the minimum rank of \Cref{def-appendix-min-rank}. If
  $\delta(C)=1$, then $C$ is the transvection class, and the rank-one product
  criterion (\Cref{thm-prescribed-rank-one}) gives $C^n=G$, including $I$ by its
  cancellation clause.

  Suppose $\delta(C)\ge 2$ and set
  \[
  j=\left\lceil \frac{n-1}{\delta(C)} \right\rceil, \quad L=2j+1.
\]
  The minimum-rank product proposition (\Cref{prop-minimum-rank-product}) supplies
  a full cyclic class $\mathcal{K} \subseteq C^j$ (when $C$ is already cyclic, this
  simply takes $j=1$ and $\mathcal{K}=C$; the excluded case $(n,F)=(4,\mathbb{F}_2)$ has
  already been treated). Nonscalar coverage from a cyclic power
  (\Cref{prop-cyclic-power}) now gives
  \[
  G \setminus \{I\} \subseteq C^L,
\]
  whose full-class hypothesis holds automatically because all conjugation here
  is in $G$. If $L<n$, the one-letter completion lemma
  (\Cref{lem-one-letter-completion}) yields $C^{L+1}=G$, and full-fibre
  persistence (\Cref{lem-fibre-persistence}) gives $C^n=G$.

  For odd $n$, the bound $\delta(C)\ge 2$ implies $j\le (n-1)/2$ and hence $L\le n$,
  so only $L=n$ remains. Let $P=\operatorname{Comp}(x^n-1)$. Since $n$ is odd, $1$ is a
  simple root of $x^n-1$ over $\mathbb{F}_2$, so $P^{-1}$ is nice cyclic. For any
  $c \in C$, the cyclic-product theorem \cite[Proposition 4.2]{KnuppelNielsen2006}
  gives a factorisation $c=s p^{-1}$ with both $s$ and $p$ conjugate to $P$;
  in particular, $s^n=p^n=I$. The rotation identity (\Cref{lem-rotation-identity})
  then yields
  \[
  \prod_{i=0}^{n-1} p^i c p^{-i}=(c p)^n p^{-n}=I,
\]
  expressing $I$ as $n$ positive letters from $C$ and completing the
  odd-dimensional case (including $n=3$).

  If $n=2m\ge 6$ and $\delta(C)\ge 3$, then
  \[
  j\le \left\lceil \frac{2m-1}{3} \right\rceil\le m-1
\]
  because $2m-1\le 3m-3$; hence $L\le n-1$, and completion applies as above. For
  the remaining case, $n\ge 6$ even and $\delta(C)=2$, we use the rank-two
  comparison (\Cref{lem-f2-pairs}) proved below. The rank-one product criterion
  (\Cref{thm-prescribed-rank-one}) writes every $X \in G$, including $I$, as a
  product of exactly $n$ transvections. Each adjacent pair has deviation rank
  at most two and therefore lies in $C^2$ by \Cref{lem-f2-pairs}; replacing all
  $n/2$ pairs gives $X \in C^n$.
\end{proof}

\subsubsection{Rank-two comparison}
\label{sec-f2-rank-two-comparison}
\begin{figure}[htbp]
  \centering
  \includegraphics[width=368pt]{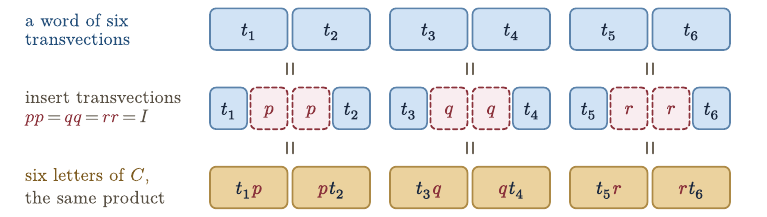}
  \caption{Replacing pairs in a six-transvection word over $\mathbb{F}_2$ by
    inserting cancelling transvections $p^2 = q^2 = r^2 = I$ so each regrouped
    factor lies in $C$. In even dimension $n \ge  6$ the same replacement covers
    all $n/2$ pairs.}\label{fig:sl-transvection-pairs}
\end{figure}

For a seed of minimum rank two in even dimension $n\ge 6$, building a cyclic
class can take $n/2$ letters, pushing the cyclic-power bound $2j+1$ one letter
past $n$. Because the transvection theorem already supplies an $n$-letter word
for every target, it suffices to replace consecutive pairs of transvections by
two letters from $C$ without changing their product (\Cref{fig:sl-transvection-pairs}).

\begin{lemma}[{Rank-two comparison}]
\label{lem-f2-pairs}
Let $n\ge 6$ be even and let
  $C$ be a conjugacy class in $\mathrm{GL}_n(\mathbb{F}_2)$ with $\delta(C)=2$.
  Then
  \[
  \{Y \in \mathrm{GL}_n(\mathbb{F}_2): \operatorname{rank}(Y-I)\le 2\} \subseteq C^2.
\]
\end{lemma}

\begin{proof}
A member $g \in C$ has $n-2$ nonconstant invariant factors. They cannot
  all have degree at least two, since $2(n-2)>n$. The smallest therefore
  has degree one and must be $x+1$. It divides all the others, so
  $\operatorname{rank}(g-I)=2$.

  The nontrivial part of $g$ consequently has one of the three forms
  \[
  J_2(1) \oplus  J_2(1), \quad J_3(1), \quad
    Q=\operatorname{Comp}(x^2+x+1)=\begin{pmatrix}0 & 1 \\ 1 & 1\end{pmatrix},
\]
  with identity on the remaining space. Indeed, any primary component
  without eigenvalue $1$ contributes its full dimension to the deviation
  rank and so can only be the irreducible quadratic component. Otherwise
  the unipotent Jordan blocks contribute a total rank of two, giving the
  first two forms.

  Fix $Y$ of deviation rank at most two. The rank-one criterion expresses
  it as $Y=t_1 t_2$ with both $t_i$ transvections, also when $Y=I$.
  Write
  \[
  t_i=I+v_i \varphi_i, \quad \varphi_i(v_i)=0,
\]
  where $v_i$ and $\varphi_i$ are nonzero. Choose $(a,b)$ to be $(0,0)$,
  $(0,1)$ or $(1,1)$ according to the three seed types above. We seek
  a vector $w$ outside $\operatorname{span}\{v_1,v_2\}$ and a functional $\psi$
  satisfying
  \[
  \varphi_1(w)=\varphi_2(w)=a, \quad
    \psi(v_1)=\psi(v_2)=b, \quad \psi(w)=0,
\]
  with $\psi$ independent of each $\varphi_i$. First we check what these
  conditions give; their simultaneous existence is proved below.

  Set $p=I+w \psi$. This is a transvection and $p^2=I$. Multiplication gives
  \[
  N_i\coloneqq t_i p-I=v_i \varphi_i+(w+a v_i)\psi.
\]
  The two image vectors are independent, as are the two functionals,
  so $N_i$ has rank two. The parameters determine its type:

\begin{itemize}
  \item If $(a,b)=(0,0)$, then $N_i^2=0$, giving $J_2(1) \oplus  J_2(1)$.
  \item If $(a,b)=(0,1)$, then $N_i^2=w \varphi_i\neq 0$ and $N_i^3=0$, giving
$J_3(1)$.
  \item If $(a,b)=(1,1)$, the matrix of $t_i p$ on
$\operatorname{span}\{v_i,w\}$ in that ordered basis is $Q$. The space
$\ker \varphi_i \cap \ker \psi$ is a fixed complement: the restrictions
of $\varphi_i,\psi$ to this two-space have matrix $\begin{pmatrix}0 & 1 \\ 1 & 0\end{pmatrix}$.
This gives the irreducible quadratic type.
\end{itemize}

  Thus $t_1 p,t_2 p \in C$. Since $p t_2=p(t_2 p)p$, also $p t_2 \in C$,
  and $Y=(t_1 p)(p t_2) \in C^2$.

  It remains to make the choices. Over $\mathbb{F}_2$, two nonzero functionals
  are either equal or independent. The equations
  $\varphi_1(w)=\varphi_2(w)=a$ are therefore consistent and their affine solution
  space has dimension at least $n-2$. It has at least eight points,
  whereas $\operatorname{span}\{v_1,v_2\}$ has at most four. Choose $w$ outside
  that span.

  The prescriptions for $\psi$ on $v_1,v_2,w$ are also consistent:
  $v_1,v_2$ are equal or independent, and $w$ lies outside their span.
  If $b=1$, any extension is nonzero and differs from each $\varphi_i$,
  since $\psi(v_i)=1$ and $\varphi_i(v_i)=0$. If $b=0$, the annihilator of
  $\operatorname{span}\{v_1,v_2,w\}$ has dimension at least $n-3\ge 2$.
  Choose $\psi$ in it outside the at most three elements $0,\varphi_1,\varphi_2$.
  Over $\mathbb{F}_2$ this is precisely the required independence.
  All the conditions used in the type checks can therefore be met.
\end{proof}

Because the inserted transvection can be chosen afresh for each pair, the
replacements assemble cleanly even when neighbouring pairs have unrelated
products, and the same argument produces an $n$-letter positive word when the
target is the identity.

\subsection{Positive products in \texorpdfstring{$\mathrm{SL}_4(\mathbb{F}_3)$}{SL4(F3)}}
\label{sec-sl-f2-comparison}
Over $\mathbb{F}_3$ the choice of conjugating group begins to matter: in dimension
four we must build every target from five letters of a single prescribed
special linear class, where a general linear change of basis can move a letter
outside that class. We first dispose of the seeds handled by the rank-one
criterion, cyclic coverage, or a cyclic three-space summand; the distinction
between $\mathrm{SL}_4(\mathbb{F}_3)$ and $\mathrm{GL}_4(\mathbb{F}_3)$ becomes decisive for the
merely quadratic seeds that remain.

\textbf{Restatement of \Cref{prop-sl-small-field-coverage}(ii) (Positive coverage in
  $\mathrm{SL}_4(\mathbb{F}_3)$).} For every noncentral conjugacy class $C$ in
  $\mathrm{SL}_4(\mathbb{F}_3)$,
  \[
  C^5=\mathrm{SL}_4(\mathbb{F}_3).
\]
  In particular, a noncentral full class in $\mathrm{GL}_4(\mathbb{F}_3)$ of determinant
  one has the same fifth power, since it contains such an $\mathrm{SL}_4(\mathbb{F}_3)$
  class.

\begin{lemma}[{Initial reductions}]
\label{lem-sl43-nonquadratic}
Let $C$ be a noncentral conjugacy class in $\mathrm{SL}_4(\mathbb{F}_3)$.
  If its nonconstant invariant factors do not have degrees $(2,2)$,
  then $C^5=\mathrm{SL}_4(\mathbb{F}_3)$.
  More precisely, if $C$ is rank one, then
  $C^4=\mathrm{SL}_4(\mathbb{F}_3) \setminus \{-I\}$; if $C$ is cyclic, then
  $C^4=\mathrm{SL}_4(\mathbb{F}_3)$.
\end{lemma}

\begin{proof}
The possible degree lists, apart from the scalar case, are
  $(1,1,2)$, $(1,3)$, $(2,2)$ and $(4)$.

  In the first case the class has rank one at an origin $b \in \{1,-1\}$. Since
  the determinant is one, its relative ratio is $1$, so $C=b T$, where
  $T$ is the transvection class. This full GL class is a single SL class:
  arbitrary determinants on a fixed two-dimensional complement centralise a
  transvection and adjust any GL conjugator to determinant one. Now
  $C^4=T^4$, and the rank-one product criterion
  (\Cref{thm-prescribed-rank-one}) gives
  $T^4=\mathrm{SL}_4(\mathbb{F}_3) \setminus \{-I\}$; one-letter completion
  (\Cref{lem-one-letter-completion}) then yields $C^5=\mathrm{SL}_4(\mathbb{F}_3)$.

  For a cyclic class, \cite[Theorem 3]{Lev1993}, with three copies of the prescribed
  SL class, puts every nonscalar in $C^3$. One-letter completion
  (\Cref{lem-one-letter-completion}) gives $C^4=\mathrm{SL}_4(\mathbb{F}_3)$, and
  full-fibre persistence (\Cref{lem-fibre-persistence}) gives the fifth power as
  well.

  For degree list $(1,3)$, write a representative as $b \oplus  B$, where
  $b=\pm 1$ and $B$ is cyclic on a three-space. Scaling the first line
  supplies arbitrary centraliser determinants, so again the full GL class is a
  single SL class. Let $\mathcal{B}$ be the full class of $B$ in
  $\mathrm{GL}_3(\mathbb{F}_3)$ and let $\mathcal{Q}$ be the nice cyclic class with
  characteristic polynomial $(x+1)(x^2-x-1)$. Its determinant is
  $1=(\det B)^2$, so the nice cyclic partner lemma (\Cref{lem-cyclic-partner})
  gives $\mathcal{Q} \subseteq \mathcal{B}^2$. For $B_1,B_2 \in \mathcal{B}$, both lifted
  factors $b \oplus  B_i$ lie in $C$, and since $b^2=1$ their product is
  \[
  (b \oplus  B_1)(b \oplus  B_2)=1 \oplus  (B_1 B_2).
\]
  Lifting $\mathcal{Q} \subseteq \mathcal{B}^2$ in this way and using GL invariance of
  $C^2$ gives $\mathcal{H}_0 \subseteq C^2$, where $\mathcal{H}_0$ is the full
  $\mathrm{GL}_4(\mathbb{F}_3)$-class of $1 \oplus  \operatorname{Comp}((x+1)(x^2-x-1))$. Its
  characteristic polynomial $(x-1)(x+1)(x^2-x-1)$ has pairwise coprime factors,
  so $\mathcal{H}_0$ is nice cyclic of determinant one. The cyclic decomposition
  theorem \cite[Theorem 1.2]{BungerNielsen1999} puts every nonscalar in
  $\mathcal{H}_0^2 \subseteq C^4$, and one-letter completion
  (\Cref{lem-one-letter-completion}) gives the fifth power.
\end{proof}

\begin{lemma}[{The quadratic types}]
\label{lem-sl43-quadratic-types}
Let $C$ be a conjugacy class in $\mathrm{SL}_4(\mathbb{F}_3)$ whose invariant
  factors have degrees $(2,2)$. Both factors equal one polynomial $p$,
  where
  \[
  p \in \{(x-1)^2,(x+1)^2,x^2+1,x^2-1,x^2-x-1,x^2+x-1\}.
\]
  For the last four choices, $C$ is also a full GL class.
\end{lemma}

\begin{proof}
The first invariant factor divides the second, and they have the
  same degree, so they coincide. The list consists of all six monic
  quadratic polynomials with nonzero constant term over $\mathbb{F}_3$.
  For an irreducible quadratic, the centraliser is
  $\mathrm{GL}_2(\mathbb{F}_9)$, whose determinant over $\mathbb{F}_3$ is the
  surjective norm of its $\mathbb{F}_9$-determinant. For $p=x^2-1$, the
  centraliser contains arbitrary invertible maps on either two-dimensional
  eigenspace, again giving every determinant. In these four cases a
  GL conjugator can therefore be adjusted to determinant one.
\end{proof}

Under general linear conjugation, a single cyclic element in $C^2$ would
suffice; with a prescribed special linear seed class $C$, however, only its
$\mathrm{SL}_4(\mathbb{F}_3)$-orbit is guaranteed to stay inside $C^2$. For four of the
six quadratic polynomials the seed class is already a full
$\mathrm{GL}_4(\mathbb{F}_3)$-class, whereas for the two repeated-root polynomials we
construct a cyclic product whose special linear orbit is itself a full class.
Once a full cyclic class lies in $C^2$, only the two scalars $\pm I$
still need five-letter words.

\begin{lemma}[{Two ingredients for five factors}]
\label{lem-sl43-two-ingredients}
Let $C$ be a noncentral conjugacy class in $\mathrm{SL}_4(\mathbb{F}_3)$.
  Suppose a full cyclic class $\mathcal{K}$ in $\mathrm{GL}_4(\mathbb{F}_3)$ lies in $C^2$.
  Then every nonscalar in $\mathrm{SL}_4(\mathbb{F}_3)$ lies in $C^5$.
  Moreover, the nice full class $\mathcal{H}$ with characteristic polynomial
  $(x-1)(x+1)(x^2-x-1)$ satisfies
  \[
  \mathcal{H} \subseteq C \mathcal{K} \subseteq C^3.
\]
  Consequently, if $I,-I \in C^5$, then $C^5=\mathrm{SL}_4(\mathbb{F}_3)$.
\end{lemma}

\begin{proof}
Both $\mathcal{K}$ and $\mathcal{H}$ have determinant one. Applying the nice cyclic
  partner lemma (\Cref{lem-cyclic-partner}) and nonscalar coverage from a cyclic
  power (\Cref{prop-cyclic-power}) with $j=2$ and $d=1$ gives the displayed
  inclusion and puts every nonscalar in $\mathcal{H}\mathcal{K} \subseteq C^5$, leaving
  only the scalars $I$ and $-I$.
\end{proof}

For the repeated-root polynomials we start from a block shear and conjugate it
into the opposite block position so that two copies act cyclically on the whole
four-space; choosing their product to have irreducible characteristic
polynomial ensures that its special linear orbit already fills its full
$\mathrm{GL}_4(\mathbb{F}_3)$-class.

\begin{lemma}[{Quadratic unipotent classes}]
\label{lem-sl43-unipotent}
Let $C$ be a prescribed conjugacy class in $\mathrm{SL}_4(\mathbb{F}_3)$ whose
  two invariant factors are both $(x-\varepsilon)^2$, where $\varepsilon=\pm 1$.
  Let $\mathcal{K}$ be the full GL
  class with characteristic polynomial
  \[
  \Phi_{10}(x)=x^4-x^3+x^2-x+1.
\]
  Then $\mathcal{K} \subseteq C^2$, the class $C$ is inverse-stable, and
  $I,-I \in C^3 \cap C^5$. In particular, $C^5=\mathrm{SL}_4(\mathbb{F}_3)$.
\end{lemma}

\begin{proof}
Begin with $\varepsilon=1$ and the SL class $C_0$ of the matrix $A$ below.
  All blocks have size two. Put
  \[
  J=\begin{pmatrix}0 & 2 \\ 1 & 0\end{pmatrix}, \quad
    A=\begin{pmatrix}I & I \\ 0 & I\end{pmatrix}, \quad B=\begin{pmatrix}I & 0 \\ J & I\end{pmatrix}, \quad h=\begin{pmatrix}0 & I \\ J & 0\end{pmatrix}.
\]
  Since $\det J=1$, we have $\det h=1$, and direct multiplication gives
  $B=h A h^{-1}$. Thus $A,B \in C_0$. Their product satisfies
  \[
  \det(x I-A B)
    =\det((x-1)^2 I-x J)
    =(x-1)^4+x^2
    =\Phi_{10}(x).
\]
  The first equality can be computed over $\mathbb{F}_3(x)$ by a block
  determinant and hence holds as a polynomial identity.

  Every root of $\Phi_{10}$ has order ten: use
  $(x+1)\Phi_{10}=x^5+1$ and $\Phi_{10}(-1)=2$. Neither $\mathbb{F}_3$ nor
  $\mathbb{F}_9$ has an element of that order, so the quartic has no linear
  or quadratic factor and is irreducible. The centraliser of $A B$ is
  $\mathbb{F}_{81}^\times$, with surjective determinant norm onto $\mathbb{F}_3^\times$.
  Its full GL class $\mathcal{K}$ is therefore a single SL class. Since
  $C_0^2$ is SL invariant and contains $A B$, it contains $\mathcal{K}$.

  Every prescribed SL class of this Jordan type is a GL conjugate of
  $C_0$. Conjugate the whole inclusion by a matrix carrying $A$ into
  the desired class. Normality of SL in GL keeps the transported
  factors in that prescribed SL class, while the full class $\mathcal{K}$
  is unchanged. This proves $\mathcal{K} \subseteq C^2$ for every positive
  unipotent class. Negating both factors gives the same inclusion for
  every negative unipotent class.

  For the scalar words, $A^3=I$ already gives $I \in C_0^3$. Reuse $A$
  and put
  \[
  v=\begin{pmatrix}I & 0 \\ -I & I\end{pmatrix}.
\]
  Then $\det v=1$, $v^3=I$, and
  \[
  A v+I=\begin{pmatrix}I & I \\ -I & -I\end{pmatrix}, \quad (A v+I)^2=0.
\]
  In characteristic three this gives $(A v)^3=-I$, and the rotation identity
  (\Cref{lem-rotation-identity}) reads
  \[
  A(v A v^{-1})(v^2 A v^{-2})=(A v)^3 v^{-3}=-I.
\]
  Each factor belongs to $C_0$ because its conjugator has determinant one.
  Transporting this entire word by GL conjugation gives it in every positive
  unipotent SL class, and negating all three letters swaps the two scalar
  products to supply both scalars for the negative classes.

  Finally, $\operatorname{diag}(-I_2,I_2)$ has determinant one and conjugates
  $A$ to $A^{-1}$. Transport and negation prove inverse stability for
  every class under consideration. Appending $c c^{-1}$ to either
  scalar word puts $I,-I$ in $C^5$, and \Cref{lem-sl43-two-ingredients} finishes
  the proof.
\end{proof}

Although two of these shears have enough combined deviation rank to move the
whole four-space, the tight scalar-action lemma (\Cref{lem-tight-scalar-action})
shows that three letters are in fact minimal for $-I$.

\begin{remark}[{The unipotent scalar word is shortest}]
\label{rem-sl43-unipotent-minimality}
For a class $C$ in \Cref{lem-sl43-unipotent}, suppose $S_1 S_2=-I$
  with $S_i \in C$, and put $N_i=S_i-\varepsilon I$. Each $N_i$ has rank
  two and square zero. The product origin is $\varepsilon^2=1$, and
  $\rho_1(-I)=4=2+2$, so the product is tight.
  \Cref{lem-tight-scalar-action} forces
  \[
  N_i^2=\varepsilon(-1-1)N_i=\varepsilon N_i,
\]
  a contradiction. A single member of $C$ is noncentral, while the
  preceding construction reaches $-I$ in three positive letters.
\end{remark}

For the other four quadratic polynomials, a single cyclic product serves all
four seed types, and only the closure onto the scalars differs: split
involutions and square roots of $-I$ return to the seed class in two letters,
whereas order-eight seeds use the cyclic product to return in four, so that in
every case one final letter closes the word.

\begin{lemma}[{The other quadratic classes}]
\label{lem-sl43-other-quadratic}
Let $C$ be the class in $\mathrm{SL}_4(\mathbb{F}_3)$ with two equal invariant
  factors $p$, where
  \[
  p \in \{x^2+1,x^2-1,x^2-x-1,x^2+x-1\}.
\]
  If $\mathcal{E}$ is the full cyclic GL class with polynomial $x^4+1$,
  then $\mathcal{E} \subseteq C^2$ and $I,-I \in C^5$.
  Consequently $C^5=\mathrm{SL}_4(\mathbb{F}_3)$.
\end{lemma}

\begin{proof}
Write $p=x^2-\alpha x-\beta$. In the quadratic-pair step of the minimum-rank
  product proposition (\Cref{prop-minimum-rank-product}), two copies of this
  polynomial give $\gamma=\beta^2=1$. Choose $Z=\operatorname{Comp}(x^2-x-1)$. The
  characteristic polynomial of $Z^{-1}$ is $x^2+x-1$, a different irreducible
  quadratic, so the construction supplies a cyclic product with polynomial
  \[
  (x^2-x-1)(x^2+x-1)=x^4+1.
\]
  Since $C$ is a full GL class, conjugacy invariance gives
  $\mathcal{E} \subseteq C^2$.

To construct scalar words, we use scalar reversal. For a nonzero
    scalar $z$ and $m\ge 2$, removing or appending the final factor gives
    \[
  z I \in C^m \quad \iff  \quad z C^{-1} \subseteq C^{m-1}.
\]
    One factorisation gives one element of the reversed class, and
    conjugacy invariance yields the rest.

  Let $\Gamma$ denote the class with minimum polynomial $x^2+1$.
  For $\Gamma$ and the split involution class, $C=C^{-1}=-C$.
  Scalar reversal therefore gives both scalar cubes once
  $C \subseteq C^2$. We construct that inclusion for each class.
  For $\Gamma$, reuse $J=\begin{pmatrix}0 & 2 \\ 1 & 0\end{pmatrix}$ from above and set
  \[
  X_0=\begin{pmatrix}0 & I \\ -I & 0\end{pmatrix}, \quad
    Y_0=\operatorname{diag}(J,-J), \quad X_0 Y_0=\begin{pmatrix}0 & -J \\ -J & 0\end{pmatrix}.
\]
  All three matrices square to $-I$, and hence belong to $\Gamma$.
  Thus $\Gamma \subseteq \Gamma^2$. For the split involution class $C$,
  the identity
  \[
  \operatorname{diag}(1,1,-1,-1)\operatorname{diag}(1,-1,-1,1)
    =\operatorname{diag}(1,-1,1,-1)
\]
  similarly gives $C \subseteq C^2$. Thus $I,-I \in C^3$ in both
  cases, and an inverse pair extends both words to length five.

  Now let $\mathcal{P}$ be either order-eight class, with minimum polynomial
  $x^2-a x-1$, $a=\pm 1$. Its quadratic relation implies $u^4=-I$ for
  every $u \in \mathcal{P}$. Squaring therefore supplies
  $\Gamma \subseteq \mathcal{P}^2$. Reciprocation and negation interchange
  the two polynomials, so $\mathcal{P}^{-1}=-\mathcal{P}$; we also already have
  $\mathcal{E} \subseteq \mathcal{P}^2$.

  Take $e \in \mathcal{E}$ and put $p_0=I+e^2$ and $g=p_0 e^{-1}$.
  Since $e^4=-I$, we obtain
  \[
  p_0^2=I-p_0=-e^2, \quad g^2=p_0^2 e^{-2}=-I.
\]
  Hence $g \in \Gamma$, while $p_0$ has irreducible minimum polynomial
  $x^2+x-1$ and belongs to one order-eight class. The factorisation
  $p_0=g e$ puts that whole class in $\Gamma \mathcal{E}$. Since
  $-\Gamma=\Gamma$, this product is stable under negation and contains
  the other order-eight class too. For either choice of seed,
  \[
  \mathcal{P} \cup \mathcal{P}^{-1} \subseteq \Gamma \mathcal{E}
    \subseteq \mathcal{P}^4.
\]
  Scalar reversal gives $I \in \mathcal{P}^5$ from the inverse class and
  $-I \in \mathcal{P}^5$ from $-\mathcal{P}^{-1}=\mathcal{P}$.
  Apply \Cref{lem-sl43-two-ingredients} in all four cases.
\end{proof}

\begin{remark}[{Shorter scalar products}]
\label{rem-sl43-short-scalars}
The constructions also retain the shorter endpoint $-I \in C^2$ for
  each of the four non-unipotent quadratic classes. In $\Gamma$ this
  follows by squaring. For the split involution class and either
  order-eight class, it follows from $-C^{-1}=C$.
\end{remark}

\begin{proof}[{Proof of \Cref{prop-sl-small-field-coverage}(ii)}]
Let $C$ be any noncentral conjugacy class in
  $\mathrm{SL}_4(\mathbb{F}_3)$. If its invariant-factor degrees differ from
  $(2,2)$, apply \Cref{lem-sl43-nonquadratic}. Otherwise
  \Cref{lem-sl43-quadratic-types} gives six possibilities: the first two
  are covered by \Cref{lem-sl43-unipotent} and the remaining four by
  \Cref{lem-sl43-other-quadratic}. Thus $C^5=\mathrm{SL}_4(\mathbb{F}_3)$ in every
  case, including both scalar targets at that exact positive length.
\end{proof}

This special linear result also feeds directly into determinant-fibre covering
below: every noncentral full class of determinant one in $\mathrm{GL}_4(\mathbb{F}_3)$
contains a prescribed $\mathrm{SL}_4(\mathbb{F}_3)$-class, so its fifth power fills the
entire determinant-one fibre.

\subsection{Determinant fibres over \texorpdfstring{$\mathbb{F}_3$}{F3}}
\label{sec-fibre-f3}
Over $\mathbb{F}_3$, parity determines how the two scalars $\pm I$ sit in the
determinant fibres. In odd dimension they lie in opposite fibres, so each fibre
contains a single scalar that must be reached at the same length $n$ as the
nonscalars; in even dimension both scalars have determinant one, leaving the
determinant-minus-one fibre scalar-free so that covering its nonscalars fills
it outright.

\textbf{Restatement of \Cref{prop-gl-f3-fibre} (Determinant fibres over $\mathbb{F}_3$).}
  Let $n\ge 3$ and let $C$ be a noncentral full conjugacy class in
  $\mathrm{GL}_n(\mathbb{F}_3)$ of determinant $d$. Put
  \[
  k=\begin{cases}n & \textnormal{if } n \textnormal{ is odd,} \\ n+1 & \textnormal{if } n \textnormal{ is even.}\end{cases}
\]
  Then
  \[
  C^k=D_k(C)=\mathcal{D}(d^k)=\{X \in \mathrm{GL}_n(\mathbb{F}_3):\det X=d^k\}.
\]
  Here $\mu_n(\mathbb{F}_3)$ is trivial for odd $n$ and equals $\{1,-1\}$
  for even $n$.

\subsubsection{The determinant scalar}
\label{sec-f3-determinant-scalar}
A primary seed has a single irreducible factor in its characteristic
polynomial, so over the corresponding extension field its semisimple part acts
as a scalar. Running through the Frobenius conjugates of that scalar multiplies
them to a norm, matching the determinant, while suitable powers of the
unipotent part cancel. For seeds with more than one primary component, the
rotation identity (\Cref{lem-rotation-identity}) supplies a shorter route to the
same scalar.

\begin{lemma}[{The determinant scalar over $\mathbb{F}_3$}]
\label{lem-f3-det-scalar}
For $n\ge 3$ and every noncentral
  full class $C$ in $\mathrm{GL}_n(\mathbb{F}_3)$, with $d=\det C$,
  \[
  d I \in C^n.
\]
\end{lemma}

\begin{proof}
Choose $g \in C$. First suppose $g$ is primary, so its characteristic
  polynomial is a power of one irreducible polynomial of degree $a$.
  Write $n=a k$. The semisimple part of $g$ equips the underlying space
  with the structure of a $k$-dimensional space over
  $E=\mathbb{F}_{3^a}$, on which
  \[
  g=\alpha U, \quad \alpha \in E^\times,
\]
  with $U$ unipotent. Choose a Jordan basis over $E$; the entries of
  $U$ then lie in $\mathbb{F}_3$.

  Form a list $\alpha_1,\dots,\alpha_n$ containing $k$ copies of each
  $\alpha^{3^j}$, $0\le j<a$. Choose $n$ integer exponents
  $e_1,\dots,e_n$, none divisible by three, with sum zero. For even $n$
  take pairs $1,-1$; for odd $n\ge 3$, take one triple $1,1,-2$ and
  then pairs. Set
  \[
  g_i=\alpha_i U^{e_i} \quad (1\le i\le n).
\]

  Each $g_i$ is GL conjugate to $g$. Coordinatewise Frobenius is
  $\mathbb{F}_3$-linear, sends $\alpha$ to its Frobenius conjugates and fixes
  $U$. Moreover, raising a unipotent matrix to an exponent prime to
  three preserves its Jordan block sizes. Indeed, if $U=I+N$, then
  $U^e-I=N Q(N)$ for a polynomial with $Q(0)=e\neq 0$ in $\mathbb{F}_3$;
  this holds also for negative $e$ by the finite nilpotent expansion.
  Since $Q(N)$ is invertible and commutes with $N$, all powers of
  $U^e-I$ have the same kernel dimensions as the corresponding powers
  of $N$. Thus $U^e$ is conjugate to $U$ over $E$, which also gives
  conjugacy over $\mathbb{F}_3$.

  The $n$ factors commute, and their product is
  \[
  \prod_{i=1}^n g_i
    = (\prod_{j=0}^{a-1} \alpha^{3^j})^k U^{\sum_{i=1}^n e_i}
    =\operatorname{N}_{E / \mathbb{F}_3}(\alpha)^k I
    =d I.
\]
  For the last equality, the determinant of multiplication by $\alpha$
  on $E^k$ is its norm to the power $k$, while $U$ has determinant
  one. There are exactly $a k=n$ factors, each in $C$. This also
  includes $k=1$, when $U=I$.

  Now suppose $g$ is nonprimary. Its characteristic polynomial has two coprime
  nonconstant factors, so by \cite[Theorem 4.2]{BungerNielsen1999} (equivalently
  \cite[Proposition 4.8]{KnuppelNielsen2006}, with the determinant-sign correction
  recorded in the remark following \Cref{prop-cyclic-power}) it admits a
  factorisation $g=s p^{-1}$
  where $s$ is conjugate to $\operatorname{Comp}(x^n-d)$ and $p$ to
  $\operatorname{Comp}(x^n-1)$. Their determinants satisfy $\det(s) \det(p)^{-1}=d$.
  The companion polynomials give $s^n=d I$ and $p^n=I$, so the rotation identity
  (\Cref{lem-rotation-identity}) gives
  \[
  \prod_{i=0}^{n-1} p^i g p^{-i}=(g p)^n p^{-n}=d I,
\]
  again using exactly $n$ positive letters from the full class $C$.
\end{proof}

\subsubsection{Completing the fibres}

We can now assemble the special linear results inside each fibre. In odd
dimension, multiplying a seed by its determinant moves it into
$\mathrm{SL}_n(\mathbb{F}_3)$; scaling back recovers all nonscalars, and the
determinant-scalar lemma (\Cref{lem-f3-det-scalar}) supplies the single remaining
scalar. In even dimension, a seed of determinant one stays in
$\mathrm{SL}_n(\mathbb{F}_3)$ throughout, whereas a seed of determinant minus one visits
the scalar-free fibre at every odd length and finishes as soon as the cyclic
construction reaches the nonscalars.

\begin{proof}[{Proof of \Cref{prop-gl-f3-fibre}}]
Throughout, $d=\det C \in \{1,-1\}$ and
  $D_j(C)=\mathcal{D}(d^j)=\{X:\det X=d^j\}$.

  \textbf{Odd dimension.} Let $n\ge 3$ be odd. The full class $d C$ has determinant
  $d^{n+1}=1$. Choose an SL-conjugacy subclass of $d C$; since $n\neq 4$,
  nonscalar coverage over $\mathbb{F}_3$ (\Cref{prop-sl-f3-nonscalar-coverage}) puts
  every nonscalar in $\mathrm{SL}_n(\mathbb{F}_3)$ into the \mbox{$n$th} power of this subclass.
  For any nonscalar $X \in D_n(C)$, the matrix $d X$ has determinant one, so
  \[
  d X \in (d C)^n=d C^n
\]
  (using $d^n=d$), giving $X \in C^n$. A scalar $z I$ lies in $D_n(C)$
  precisely when $z^n=d^n$, which in odd dimension forces $z=d$, and the
  determinant-scalar lemma (\Cref{lem-f3-det-scalar}) gives $d I \in C^n$.
  At length $n$, both parts of the fibre are accounted for:
  \begin{align*}
  D_n(C) \setminus \{d I\} & \subseteq C^n && \textnormal{(nonscalar coverage)}, \\
    \{d I\} & \subseteq C^n && \textnormal{(determinant scalar)}.
\end{align*}
  Hence $C^n=D_n(C)$.

  \textbf{Even dimension, determinant one.} Suppose $n\ge 6$ is even and $d=1$. Choose an
  SL subclass of $C$ and apply \Cref{prop-sl-f3-nonscalar-coverage} to obtain
  every nonscalar of $\mathrm{SL}_n(\mathbb{F}_3)$ in $C^n$. The one-letter completion
  lemma (\Cref{lem-one-letter-completion}) then gives
  \[
  C^{n+1}=\mathrm{SL}_n(\mathbb{F}_3)=D_{n+1}(C),
\]
  including both scalar targets at exactly $n+1$ letters. For $n=4$, choose an
  actual SL-conjugacy subclass $C_0 \subseteq C$; positive coverage in
  $\mathrm{SL}_4(\mathbb{F}_3)$ (\Cref{prop-sl-small-field-coverage}(ii), proved in
  \Cref{sec-sl-f2-comparison}) gives
  \[
  \mathrm{SL}_4(\mathbb{F}_3)=C_0^5 \subseteq C^5
    \subseteq \mathrm{SL}_4(\mathbb{F}_3),
\]
  so the fibre is again full at the required length. For every even $n\ge 4$ and
  $d=1$, both parts have therefore reached length $n+1$:
  \begin{align*}
  D_{n+1}(C) \setminus \{I,-I\} & \subseteq C^{n+1} && \textnormal{(SL coverage or $n=4$),}\\
    \{I,-I\} & \subseteq C^{n+1} && \textnormal{(completion or $n=4$).}
\end{align*}

  \textbf{Even dimension, determinant minus one.} Let $n\ge 4$ be even and $d=-1$. Every
  scalar matrix has determinant one, so $D_{n+1}(C)$ contains no scalars. A
  rank-one class fills this fibre directly by the rank-one product criterion
  (\Cref{thm-prescribed-rank-one}), since the length $n+1$ exceeds every deviation
  rank.

  Otherwise $\delta(C)\ge 2$. Put
  \[
  j=\left\lceil \frac{n-1}{\delta(C)} \right\rceil, \quad L=2j+1.
\]
  The minimum-rank product proposition (\Cref{prop-minimum-rank-product}) supplies
  a full cyclic class $\mathcal{K} \subseteq C^j$, and nonscalar coverage from a
  cyclic power (\Cref{prop-cyclic-power}) puts every nonscalar of $D_L(C)$ in
  $C^L$. Both results apply over $\mathbb{F}_3$ at $n=4$ as well, without appealing
  to the excluded dimension of \Cref{prop-sl-f3-nonscalar-coverage}. Since $n$ is
  even,
  \[
  L\le 2 \left\lceil \frac{n-1}{2} \right\rceil+1=n+1.
\]
  Moreover $L$ is odd, so $d^L=-1$ and $D_L(C)$ is scalar-free, giving
  $C^L=D_L(C)$ immediately. Full-fibre persistence
  (\Cref{lem-fibre-persistence}) then yields $C^{n+1}=D_{n+1}(C)$ (allowing the
  determinant $(-1)^r$ of the filled fibres $D_r(C)$ to alternate along the
  way; for $n=4$ the same estimate gives $j\le 2$ and $L\le 5$). At the final
  length, the nonscalar and scalar obligations are
  \begin{align*}
  D_{n+1}(C) & \subseteq C^{n+1} && \textnormal{(cyclic coverage and persistence)}, \\
    D_{n+1}(C) \cap \{I,-I\} & = \emptyset && \textnormal{(no scalar targets)}.
\end{align*}

  These three cases exhaust all noncentral full classes and all $n\ge 3$.
\end{proof}

Thus the small fields match the length bound of the larger fields: over
$\mathbb{F}_2$ even the identity has an $n$-letter positive word, while over
$\mathbb{F}_3$ the extra letter in even dimension suffices both for prescribed
special linear classes in dimension four and for full determinant fibres in
every dimension $n\ge 3$. Together these are the exact-length inputs used in the
special linear diameter theorem (\Cref{thm-sl-diameter}) and the determinant-fibre
covering theorem (\Cref{thm-gl-fibre-covering}).

\bibliographystyle{plain}
\bibliography{references}

\end{document}